\documentclass{article}
\usepackage[utf8]{inputenc}
\usepackage{eufrak}
\usepackage[breaklinks=true, pdftex, 
pdfborder={0 0 0}, colorlinks=true, linkcolor=blue, citecolor=blue, urlcolor=blue, 
bookmarks=false, pdfpagelabels=false]{hyperref}
\usepackage{fullpage}
\usepackage{times}
\usepackage{ifthen}   
\usepackage{graphicx,graphics}
\usepackage{setspace}
\usepackage{amssymb,amsmath,amsthm}
\usepackage{dsfont}
\usepackage{color}
\usepackage{wrapfig}
\usepackage{multirow}
\usepackage{multicol}
\usepackage{pdfpages}
\usepackage{longtable} 	
\usepackage{booktabs}
\usepackage{datetime}
\usepackage{xspace}
\usepackage{paralist} 
\usepackage{enumitem} 
\usepackage[it,small,belowskip=-15pt]{caption}
\usepackage[top=1in,bottom=1in,left=0.95in,right=0.95in]{geometry} 
\usepackage[T1]{fontenc}
\usepackage{doi}
\usepackage[numbers,compress,sort]{natbib} 
\usepackage{algorithmic}
\usepackage[noresetcount,vlined,boxed]{algorithm2e} 
\let\oldnl\nl
\newcommand{\nonl}{\renewcommand{\nl}{\let\nl\oldnl}}
\usepackage{caption}
\usepackage{float}   
\usepackage{subcaption}
\usepackage[framemethod=tikz]{mdframed}

\usepackage{comment}
\usepackage{cleveref} 
\usepackage{fontawesome5}

\usepackage{listings}
\usepackage{upgreek}
\usepackage[bbgreekl]{mathbbol}
\DeclareSymbolFontAlphabet{\mathbb}{AMSb}
\DeclareSymbolFontAlphabet{\mathbbl}{bbold}

\usepackage{lineno}  

\definecolor{codegreen}{rgb}{0,0.6,0}
\definecolor{codegray}{rgb}{0.5,0.5,0.5}
\definecolor{codepurple}{rgb}{0.58,0,0.82}
\definecolor{backcolour}{rgb}{0.95,0.95,0.92}
\definecolor{ForestGreen}{RGB}{34,139,34}

\usepackage{xcolor}
\usepackage{listings}

\usepackage{xparse}
\usepackage{adjustbox}
\usepackage{tikz}
\usetikzlibrary{fit, arrows, calc, positioning}
\usetikzlibrary{arrows, automata}
\usetikzlibrary{positioning}
\usetikzlibrary{shadows}
\tikzset{box/.style={draw,rectangle,rounded  corners=0pt, fill=blue!15, minimum  width=3cm, minimum  height=2cm}}
\tikzset{loz/.style={draw,diamond, aspect=1.5, fill=blue!15}}

\tikzstyle{rnd}=[circle,fill=blue!25,minimum  width=5em]
\tikzset{
	quote/.style={{|[width=1mm]}-{|[width=1mm]}}
}
\usepackage{stackengine}   
\newcommand{\ubar}[1]{\stackunder[1.0pt]{$#1$}{\rule{.8ex}{.075ex}}}    
\newcommand{\ubarlong}[1]{\stackunder[1.0pt]{$#1$}{\rule{1.2ex}{.075ex}}}
\NewDocumentCommand{\codeword}{v}{%
	\texttt{\textcolor{blue}{#1}}%
}

\lstdefinestyle{mystyle}{
	backgroundcolor=\color{backcolour},   
	commentstyle=\color{codegreen},
	keywordstyle=\color{magenta},
	numberstyle=\tiny\color{codegray},
	stringstyle=\color{codepurple},
	basicstyle=\ttfamily\footnotesize,
	breakatwhitespace=false,         
	breaklines=true,                 
	captionpos=b,                    
	keepspaces=true,                 
	numbers=left,                    
	numbersep=5pt,                  
	showspaces=false,                
	showstringspaces=false,
	showtabs=false,                  
	tabsize=2
}

\usepackage{titlesec}

\usepackage{mdframed}  
\usepackage{colortbl}
\usepackage{mathrsfs} 
\definecolor{highlight}{HTML}{DFEFF9}
\mdfdefinestyle{highlightbox}{leftmargin=0pt,rightmargin=0pt,%
	innerleftmargin=3pt,innerrightmargin=2pt,%
	innertopmargin=4pt,innerbottommargin=4pt,%
	backgroundcolor=highlight,linecolor=white}
\definecolor{algm}{HTML}{D4EFDF} 
\mdfdefinestyle{highlightalg}{leftmargin=0pt,rightmargin=0pt,%
	innerleftmargin=3pt,innerrightmargin=2pt,%
	innertopmargin=4pt,innerbottommargin=4pt,%
	backgroundcolor=algm,linecolor=white}
\definecolor{soft}{HTML}{F2D7D5} 
\mdfdefinestyle{highlightstruct}{leftmargin=0pt,rightmargin=0pt,%
	innerleftmargin=3pt,innerrightmargin=2pt,%
	innertopmargin=4pt,innerbottommargin=4pt,%
	backgroundcolor=soft,linecolor=white}

\usepackage{soul}

\makeatletter
\def\namedlabel#1#2{\begingroup
	#2%
	\def\@currentlabel{#2}%
	\phantomsection\label{#1}\endgroup
}
\makeatother

\newcommand{\randmat}[1]{{\scalebox{1}{${{\bm{\mathcal{#1}}}}$}}}

\newcommand{\Reals}{\mathbb{R}} 
\newcommand{\R}{\Reals} 

\DeclareMathOperator*{\argmin}{arg\,min}

\newcommand{\kappaFL}{\kappa_{{\rm FL}}}
\newcommand{\kappaas}{\kappa_{{\rm as}}}

\newcommand{\kappabmh}{\kappa_{{\rm mh}}}
\newcommand{\kappaf}{\kappa_{\rm f}}

\newcommand{\ceil}[1] {\lceil #1 \rceil}

\newcommand{\E}[2][]{\operatorname{\mathbb{E}}_{#1}\left[ #2\right]} 
 
\usepackage{bm} 

\newcommand{\eb}{\bm{e}}

\newcommand{\ub}{\bm{u}}

\newcommand{\vb}{\bm{v}}

\newcommand{\yb}{\bm{y}}

\newtheorem{definition}{Definition}[section]
\newtheorem{theorem}{Theorem}[section]
\newtheorem{lemma}{Lemma}[section]
\newtheorem{proposition}{Proposition}[section]
\newtheorem{assumption}{Assumption}[section]
\newtheorem{corollary}{Corollary}[section]

\newtheorem{remark}{Remark}[section]

\usepackage{tikz}
\usetikzlibrary{positioning}
\usetikzlibrary{calc}
\usetikzlibrary{shapes.geometric}
\usetikzlibrary{topaths}
\usetikzlibrary{external}
\usetikzlibrary{fadings,decorations.pathreplacing} 

\usepackage{pgfplots}
\usetikzlibrary{shapes,decorations.pathmorphing,patterns}
\pgfplotsset{compat=newest}
\usepgfplotslibrary{fillbetween}

\graphicspath{{figures/}}

\usepackage{hyperref}
\hypersetup{
unicode = false,
pdftoolbar = true,
pdfmenubar = true,
pdffitwindow = true,
pdftitle = {Direct search for stochastic optimization in random subspaces},
pdfauthor = {Dzahini},
pdfsubject = {RandomMatrix},
pdfnewwindow = true,
pdfkeywords = {Random, Matrix},
colorlinks = true,
linkcolor = blue,
citecolor = blue,
filecolor = black,
urlcolor = blue,
breaklinks = true
}

\newcommand{\Z}{\mathbb{Z}}
\newcommand{\Esp}{\mathbb{E}}

\newcommand{\N}{\mathbb{N}}

\newcommand{\pr}{\mathbb{P}}

\newcommand{\Hk}{\hat{\mathsf{H}}_k}

\newcommand{\ds}{\displaystyle}

\newcommand{\dmax}{\alpha_{\max}}

\newcommand{\ijk}{\mathds{1}_{\mathscr{J}_k}}

\newcommand{\ijkc}{\mathds{1}_{\bar{\mathscr{J}}_k}}
 
\newcommand{\iik}{\mathds{1}_{\mathscr{I}_k}}

\newcommand{\iikc}{\mathds{1}_{\bar{\mathscr{I}}_k}}

\newcommand{\abs}[1]{\left\lvert#1\right\rvert}
\newcommand{\parenth}[1]{\left(#1\right)}
\newcommand{\intbracket}[1]{\left[\!\left[#1\right]\!\right]}
\newcommand{\accolade}[1]{\left\lbrace#1\right\rbrace}
\newcommand{\accoladekinN}[1]{{\left\lbrace#1\right\rbrace}_{k\in\mathbb{N}}}

\newcommand{\prob}[1]{\mathbb{P}\left(#1\right)}
\newcommand{\var}[1]{\mathbb{V}\left[#1\right]}

\newcommand{\norme}[1]{\left\lVert#1\right\rVert}
\newcommand{\normeshort}[1]{{|\hspace*{-0.3mm}|}#1{|\hspace*{-0.3mm}|}}

\newcommand{\normii}[1]{{\left\lVert#1\right\rVert}_{2}}

\newcommand{\sk}{\bm{s}_k}

\newcommand{\rd}{\mathbb{R}^d}

\newcommand{\rp}{\mathbb{R}^p}

\newcommand{\rdp}{\mathbb{R}^{d\times p}}

\newcommand{\rpp}{\mathbb{R}^{p\times p}}

\newcommand{\s}{\bm{s}}
\newcommand{\vi}{\bm{v}}

\newcommand{\Q}{\mathsf{Q}} 
\newcommand{\Qk}{\Q_k}

\newcommand{\Mhatk}{\hat{M}_k} 
\newcommand{\mhatk}{\hat{m}_k} 

\newcommand{\Mhatell}{\hat{M}_{\ell}}

\newcommand{\fok}{\bar{f}^{[0]}_k}

\newcommand{\Foell}{\bar{F}^{[0]}_{\ell}}

\newcommand{\fsk}{\bar{f}^{[s]}_k}
\newcommand{\Fsell}{\bar{F}^{[s]}_{\ell}}

\newcommand{\dkun}{{\delta_{k+1}}}
 
\newcommand{\szero}{\bm{s}^0} 
 
\newcommand{\sone}{\bm{s}^1}

\newcommand{\sii}{\bm{s}^i}

\newcommand{\X}{\bm{X}}

\newcommand{\xkun}{\bm{x}_{k+1}^0}
\newcommand{\x}{\bm{x}}
\newcommand{\y}{\bm{y}}

\newcommand{\Lg}{\kappa_{L \rm g}}

\newcommand{\gkhat}{\bm{\hat{g}}_k}
\newcommand{\gk}{\bm{\hat{g}}_k}
\newcommand{\fmqalgebra}{\mathscr{F}^{\hat{M}\cdot \mathcal{Q}}_{k-1}}
\newcommand{\falgebra}{\mathscr{F}_{k-1}}
\newcommand{\fqalgebra}{\mathscr{F}_{k-1}^{\mathcal{Q}}}
\newcommand{\ef}{\epsilon_{\rm f}}
\newcommand{\Dk}{\Delta_k}

\newcommand{\dk}{\delta_k}

\newcommand{\ilkhat}{\mathds{1}_{\mathscr{L}_k}}
\newcommand{\ilkhatc}{\mathds{1}_{\bar{\mathscr{L}}_k}}

\newcommand{\Phik}{{\Phi}_k}
\newcommand{\Phikplusone}{{\Phi}_{k+1}}
\newcommand{\Tepsilon}{T_\varepsilon}

\newcommand{\kmge}{\kappa_{{\rm mge}}} 
\newcommand{\kfcd}{\kappa_{{\rm fcd}}}

\newcommand{\kmfe}{\kappa_{{\rm mfe}}}
\newcommand{\kfde}{\kappa_{{\rm fde}}}

\newcommand{\Qmax}{\kappa_{Q_{\max}}}

\newcommand{\bl}{\color{black}}

\newcommand{\uk}
{\mathscr{U}_k}
\newcommand{\erandom}{E}

\newcommand{\ebarokcur}{\bar{e}^{[0]}_k}
\newcommand{\ebarsk}{\bar{e}^{[s]}_k}
\newcommand{\erandombar}{\bar{E}}

\newcommand{\xrandomvec}{\bm{X}}

\newcommand{\vrandomvec}{\bm{V}}

\newcommand{\zrandom}{Z}

\newcommand{\Fbar}{\bar{F}}

\newcommand{\Fbarik}{\bar{F}^i_k}
\newcommand{\Fbarok}{\bar{F}^0_k}
\newcommand{\Fbarokbrak}{\bar{F}^{[0]}_k}
\newcommand{\Fbarokcur}{\bar{F}^{[0]}_k}
\newcommand{\Fbarsk}{\bar{F}^{[s]}_k}

\newcommand{\fbarik}{\bar{f}^i_k}
\newcommand{\fbarok}{\bar{f}^0_k}
\newcommand{\fbarokcur}{\bar{f}^{[0]}_k}
\newcommand{\fbarsk}{\bar{f}^{[s]}_k}

\newcommand{\nokrandombrak}{{N^{0,s}_k}}
\newcommand{\nskrandom}{N^{0,s}_k}

\newcommand{\nokbrak}{n^{0,s}_k}
\newcommand{\nsk}{n^{0,s}_k}
\newcommand{\nskbrak}{n^{0,s}_k}

\newcommand{\xks}{\bm{x}_k^{s}}
\newcommand{\xzerok}{\bm{x}_k^{0}}
\newcommand{\QkRandom}{\randmat{Q}_k}
\newcommand{\AkEvent}{\mathscr{A}_k}
\newcommand{\YkSet}{\mathscr{Y}_k}

\newcommand{\iak}{\mathds{1}_{\mathscr{A}_k}}

\newcommand{\mhatkRandom}{\hat{M}_k}

\newcommand{\siirandom}{{\bm{S}}^i}
\newcommand{\srandom}{\bm{S}}
\newcommand{\ebhat}{\hat{\eb}}
\newcommand{\Xzerok}{\xrandomvec_k^{0}}
\newcommand{\ik}{\mathscr{I}_k} 

\newcommand{\jk}{\mathscr{J}_k}

\newcommand{\galgebra}{\mathscr{G}_{j-1}}
\newcommand{\gjalgebra}{\mathscr{G}_{j}}
\newcommand{\gnalgebra}{\mathscr{G}_{n}}

\newcommand{\xk}{\xzerok}
\newcommand{\mk}{\mhatk}

\newcommand{\ivk}{\mathds{1}_{\mathscr{V}_k}}
\newcommand{\ivkc}{\mathds{1}_{\bar{\mathscr{V}}_k}}
\newcommand{\nkrandom}{N_k}
\newcommand{\nk}{n_k}
\newcommand{\phik}{{\phi}_k}
\newcommand{\phikplusone}{{\phi}_{k+1}}
\newcommand{\brandom}{B}

\newcommand{\GkhatRandom}{{\hat{\bm{G}}}_k}
\newcommand{\Gkhat}{{\hat{\bm{G}}}_k}

\newcommand{\rrandomk}{R_k}

\newcommand{\krandom}{K}

\newcommand{\taurandom}{T}

\newcommand{\ickwk}{\mathds{1}_{\mathscr{C}_k\cap \mathscr{W}_k}}
\newcommand{\ick}{\mathds{1}_{\mathscr{C}_k}}
\newcommand{\iwk}{\mathds{1}_{\mathscr{W}_k}}
\newcommand{\vepsell}{\varepsilon_\ell}

\newcommand{\epsLanbda}{\epsilon_\lambda}

\newcommand{\Zkj}{Z_k^{(j)}}
\newcommand{\ZkPlusOnej}{Z_{k+1}^{(j)}}
\newcommand{\Xkj}{X_k^{(j)}}
\newcommand{\XkPlusOnej}{X_{k+1}^{(j)}}
\newcommand{\kkj}{\varkappa_k^{(j)}}
\newcommand{\kkPlusOnej}{\varkappa_{k+1}^{(j)}}
\newcommand{\GhatRandom}{{\hat{G}}}

\newcommand{\Qrandom}{\randmat{Q}}

\newcommand{\Hkhat}{\hat{\mathsf{H}}_k} 
\newcommand{\HkhatRandom}{\hat{\randmat{H}}_k}
\newcommand{\HhatRandomComp}{\hat{\mathcal{H}}}
\newcommand{\HtildeRandomComp}{\tilde{\mathcal{H}}}
\newcommand{\QellRandom}{\randmat{Q}_{\ell}}
\newcommand{\QkTranspose}{\mathsf{Q}_k^{\top}} 
\newcommand{\QkRandomTranspose}{\randmat{Q}_k^{\top}}
\newcommand{\ssmath}{{\textnormal{\ss}}}
\newcommand{\halgebra}{\mathscr{H}}
\newcommand{\Nkbar}{{{\stackunder[1.0pt]{$N$}{\rule{1.3ex}{.075ex}}}_k}}
\newcommand{\Xki}{\X_k^i}
\newcommand{\Zkli}{Z_{k,\ell}^i}
\newcommand{\Dkli}{D_{k,\ell}^i}
\newcommand{\Dkni}{D_{k,n}^i}
\newcommand{\Dkls}{D_{k,\ell}^s}
\newcommand{\Dkns}{D_{k,n}^s}
\newcommand{\Ekli}{E_{k,\ell}^i}
\newcommand{\Eklo}{E_{k,\ell}^0}
\newcommand{\chifde}{\chi_{\rm fde}}
\newcommand{\gqalgebra}[1]{\mathscr{G}^{\mathcal{Q}}_{#1}}
\newcommand{\gmqalgebra}[1]{\mathscr{G}^{\hat{m}\cdot\mathcal{Q}}_{#1}}
\newcommand{\fkhat}{\hat{f}_k}
\newcommand{\kappah}{\kappabmh}
\newcommand{\Dscr}{\mathscr{D}}
\newcommand{\Phikun}{\Phikplusone} 
\newcommand{\gkalgebra}{\mathscr{G}_{k-1}}
\newcommand{\Wk}{\mathcal{W}_k}

\newcommand{\Gscr}{\mathscr{G}}
\newcommand{\Vscr}{\mathscr{V}}
\newcommand{\Nscr}{\mathscr{N}}
\newcommand{\Escr}{\mathscr{E}}

\newcommand{\Stepsize}{A}  

\newcommand{\Stepsizek}{\Stepsize_k}
\newcommand{\Cdelta}{C_\Stepsize}
\newcommand{\tb}{\bm{t}}
\newcommand{\Bp}{\mathbb{B}_p}
\newcommand{\deltmax}{\delta_{\max}}

\hypersetup{pdftitle={Sample complexity and scalability of adaptive-sampling trust-region methods for zeroth-order stochastic optimization --- Part I:   common random numbers meet supermartingales},pdfauthor={K. J. Dzahini, S. Shashaani, C. Cartis},pdfsubject={},pdfkeywords={}}
\usepackage{orcidlink}
\definecolor{SSgreen}{rgb}{0.0,0.45,0.0}
\newcommand{\SSadd}[1]{{\color{black}#1}}

\providecommand{\Dbar}{\bar{D}}
\providecommand{\uk}{\mathscr{U}_k}
\providecommand{\gjalgebra}[1][j]{\mathscr{G}_{#1}}
\definecolor{SSred}{rgb}{0.75,0.0,0.0}
\providecommand{\iuk}{\mathds{1}_{\uk}}
\providecommand{\iukc}{\mathds{1}_{\uk^c}}
\renewcommand{\AkEvent}{\mathscr{W}_k}
\renewcommand{\iak}{\mathds{1}_{\mathscr{W}_k}}

\definecolor{CCpurple}{rgb}{0.50,0.00,0.60
}
\newcommand{\SSaddnew}[1]{{\color{black}#1}}

\newcommand{\ASTRAF}{\textsf{ASTRA-F}}
\newcommand{\ASTRAR}{\textsf{ASTRA-R}}
\title{\SSaddnew{Sample complexity and scalability of adaptive-sampling trust-region methods for zeroth-order stochastic optimization --- Part I:   common random numbers meet supermartingales} 
\thanks{\SSaddnew{The first two authors contributed equally to this work.}}}

\author{
	\href{mailto:kdzahini@anl.gov}{K. J. Dzahini \orcidlink{0000-0002-1515-4251}
    }
    \thanks{Argonne National Laboratory, 9700 S. Cass Avenue, Lemont, IL 60439, USA (\href{https://www.anl.gov/profile/kwassi-joseph-dzahini}{www.anl.gov/profile/kwassi-joseph-dzahini}).
	}
	\and
	\href{mailto:sshasha2@ncsu.edu}{S. Shashaani \orcidlink{0000-0001-8515-5877}
    }
    \thanks{Fitts Dept. of Industrial \& System Eng., North Carolina State University, Raleigh, NC 27695, USA (\href{https://shashaani.wordpress.ncsu.edu/}{https://shashaani.wordpress.ncsu.edu/}).
	}
    \and
    \href{mailto:cartis@maths.ox.ac.uk}{C. Cartis \orcidlink{0000-0002-0963-5550}
    }
    \thanks{Mathematical Institute, University of Oxford, Oxford OX2 6GG, UK (\href{https://www.maths.ox.ac.uk/people/coralia.cartis}{https://www.maths.ox.ac.uk/people/coralia.cartis}).
    }
}
\date{\today}
\begin{document}


\maketitle


\noindent

{\bf Abstract:} We analyze adaptive sampling stochastic model-based trust-region methods that build their models in random subspaces and use common random numbers (CRN) within an iteration, and obtain, within a STORM-type analysis, the first $\tilde{\mathcal{O}}(\varepsilon^{-4})$ bound on the number of function evaluations needed to reach $\varepsilon$-stationarity, both almost surely and in expectation. The subspace model reduces the per-iteration computational cost, which allows scalability. The analysis relies on a general framework for the complexity of adaptive stochastic processes that carries the almost-sure analysis of ASTRO-DF into the STORM setting and adds expected bounds by a direct counting argument, without the renewal–reward argument of Blanchet et al. (2019).  The bounds are obtained for two variants of a stochastic model-based trust-region method that differ  in how the function evaluations used to measure progress are generated. The first algorithm, \ASTRAF, proceeds in a STORM-like manner: the trial point and the incumbent are re-evaluated on an independent set of replications, so that the estimates entering the acceptance test are independent of the model, as every existing STORM-type analysis requires their accuracy to be assessed conditionally on the model. The second, \ASTRAR, allows a practical approach: the trial point is evaluated on the replications already drawn to build the model, which avoids re-evaluating the incumbent, needs one adaptive sample size per iteration instead of two, and reduces the randomness entering the iteration. To our knowledge, no STORM-type analysis has been able to accommodate the latter approach as the accuracy of an estimate at a point chosen from the samples cannot be assessed conditionally on the model; we overcome this with a probabilistic chaining argument.

$ $\\
{\bf Keywords} \SSaddnew{stochastic optimization; derivative-free optimization; trust-region methods; adaptive sampling; common random numbers; random subspaces; supermartingales}

$ $\\
{\bf Mathematics Subject Classification} 90C15 $\cdot$ 90C30 $\cdot$ 90C56


\section{Introduction}
We consider unconstrained stochastic optimization problems of the form 
\begin{equation}\label{probl1}
\underset{\x\in\rd}{\min} f(\x), \qquad\mbox{with}\quad f(\x)=\Esp_{\Theta}\left[F(\x; \Theta)\right],
\end{equation}
where $\Theta$ is a random variable whose distribution is possibly unknown, the values of the smooth, nonconvex and lower-bounded function $f:\rd \rightarrow \R$ are available only through the stochastically noisy version $F(\cdot;\Theta)$ of~$f$, and $\Esp_{\Theta}$ denotes the {\it expectation} with respect to $\Theta$. For example when $f(\x)$ is observable with noise via a Monte Carlo simulation, the random variable $F(\x;\Theta)$ can be decomposed as $F(\x;\Theta):=f(\x)+\erandom(\x;\Theta)$, that is, a sum of the true function value $f(\x)$ and a stochastic error $\erandom(\x;\Theta)$. Problem~\eqref{probl1} arises, for example, in statistical machine learning~\cite{bottou2018optimization} where~$\Theta$ is a random variable describing a given sample set while $F(\cdot;\Theta)$ represents the composition of a loss function and a prediction function. 
 
\SSadd{\paragraph{Contributions and positioning.} 
Our main result is the first $\tilde{\mathcal{O}}(\varepsilon^{-4})$ sample complexity bound, \SSaddnew{within a STORM-type supermartingale analysis,} for stochastic optimization of non-convex functions with trust-region methods that use common random numbers and random subspaces. The subspace model keeps the number of design points per iteration at $2p+1$ with $p$ independent of~$d$, which  allows the methods to scale. The sample complexity bound is obtained for two variants of one algorithm, \ASTRAF\ (Algorithm~\ref{algoF}) and \ASTRAR\ (Algorithm~\ref{algoR}), which differ  in how the function evaluations used to measure progress are generated. \ASTRAF\ is the variant  closer to STORM: the incumbent and the trial point are re-evaluated on replications independent of those that built the model. That independence is crucial in the existing supermartingale analyses~\cite{blanchet2016convergence,chen2018stochastic,DzaWildSub2022} where the trial estimates are assumed to be probabilistically accurate conditioned on the $\sigma$-algebra generated by the model, which the trial point depends on. \ASTRAR\ is the more practical variant that one would actually implement: it evaluates the trial point on the replications already drawn, avoids re-evaluating the incumbent and reduces the randomness entering the iteration\SSaddnew{, at a lower per-iteration cost}. It has been impossible to accommodate such a reuse in STORM-type analyses up to this point for exactly the reason above. In this paper, we manage to remove the independence requirement altogether rather than weakening it using a different probabilistic tool, namely, by replacing an accuracy statement at the trial point with one that holds uniformly over the trust region, which is indifferent to how the trial point was chosen. The uniform statement is obtained by a probabilistic chaining argument (Lemma~\ref{Lemma:Chained}), and it holds at the adaptive, data-dependent sample size.

Both algorithms are analyzed within the supermartingale framework of Section~\ref{sec:SupermartFram}, which is stated in a form that can be reused for the complexity analysis of other algorithms. \SSaddnew{The almost-sure part of that framework carries the complexity analysis of ASTRO-DF in~\cite{ha2025complexity} into a STORM-type setting, in which the accuracy events hold at every iteration with conditional probability exceeding~{\bl $1/2$} (rather than eventually along each sample path). In contrast to the STORM analysis of~\cite{blanchet2016convergence}, the bound on the expected number of iterations is obtained by a direct count over a bounded horizon, without the renewal--reward argument, and the almost-sure iteration complexity is derived from it. The expected sample complexity then follows from 
a nonasymptotic bound on the conditional expectation of the adaptive sample size (Lemma~\ref{lemma:Nk}). 
The comparison with~\cite{blanchet2016convergence} is made in Section~\ref{sec:SupermartFram}.} 
We compare and contrast the two analyses: the deterministic part and the supermartingale part are shared, the only difference being the probabilistic tool that certifies the accuracy of the trial estimate (a pointwise stopped Bernstein inequality for \ASTRAF, the chained bound for \ASTRAR), and the complexity bounds in both are the same. The sample complexity bounds are proven both almost surely and in expectation: the expected number of iterations is $\mathcal O(\varepsilon^{-2})$, and \SSaddnew{ the expected number of evaluations is $\mathcal O(\varepsilon^{-4})$ up to a logarithmic factor, and $\mathcal O(\varepsilon^{-4})$ with a constant sample-size inflation (Corollary~\ref{cor:constantInflation}).} 
We hope the uniform-accuracy device and complexity framework open the way to analyses of other algorithms that reuse samples across the steps of an iteration.

Adaptive sampling with common random numbers is known to reduce the sample size per iteration from order $\Dk^{-4}$ to order $\Dk^{-2}$, and hence the sample complexity to order $\varepsilon^{-4}$; see~\cite{ha2025complexity}. 
The same reduction is obtained in~\cite{rinaldi2024stochastic,ding2025sequential} under a weak tail-bound condition on the estimates, without a sample complexity bound. 
The present work differs from~\cite{rinaldi2024stochastic,ding2025sequential} in three respects. There, the sample sizes are prescribed deterministic functions of $\Dk$ whose constants must be set from a moment bound on the noise; here they are stopping times driven by the observed sample variance, and 
no noise constant needs to be known. There, the reduced sample size is obtained from a tail condition on the estimate of the reduction, with the trial point evaluated separately from the model; here the Bernstein condition is placed on a measure of paired increment error (see $D$ in~\eqref{eq:DIncrement}), and in \ASTRAR\ the trial point shares the model's replications, which a pointwise tail argument cannot accommodate. The probabilistic tool is also different, namely, they use a weak tail bound at a fixed sample size there, while we use a stopped (time-uniform) Bernstein inequality and its chained extension. 
} \SSaddnew{Supermartingale analyses of STORM-type methods~\cite{blanchet2016convergence,chen2018stochastic,DzaWildSub2022} and the high-probability analyses of~\cite{cao2024first} bound the number of iterations, not the number of oracle calls. The total sample complexity of STORM was obtained recently in~\cite{JinSchXieAOA25}, from a lower bound on the step-size parameter holding with high probability, for sample sizes prescribed by the step size and a known variance bound. Here the sample sizes are stopping times driven by the observed variance, the trial point may share the model's replications, and the work is bounded without any lower bound on the radius. 
}

\SSaddnew{\paragraph{Companion paper.} This is Part~I of two. Part~II analyses the reuse variant \ASTRAR\ along sample paths: the accuracy event of Section~\ref{sec:uniformAcc} and the well-alignment event of the sketch are shown, by the Borel--Cantelli lemma, to fail only finitely often, and from that random index on the argument is the deterministic trust-region proof. Where the present paper requires the accuracy and alignment events to hold at every iteration with conditional probability bounded away from one half, Part~II requires only that their failure probabilities be summable, which it achieves by letting the subspace dimension grow logarithmically with the iteration count. It proves that $\normii{\nabla f(\Xzerok)}\to0$ almost surely; that the iteration and work complexities are of orders $\varepsilon^{-2}$ and $\varepsilon^{-4}$ almost surely along deterministic sequences $\varepsilon_\ell\downarrow0$; that the random index after which every iteration is accurate has finite moments, which gives a finite-time bound holding with high probability; and that at a fixed subspace dimension $p$ its argument gives $\varepsilon^{-2-\mathcal O(1/p)}$, so that the growth of $p_k$ is the price of dropping the per-iteration probability bound. It also tracks the dependence of the deterministic constants on $d$ and $p$ for the usual sketch ensembles, which is not pursued here. The two parts share the accuracy analysis of Section~\ref{sec:uniformAcc}; the per-iteration bound used here is what yields the expected complexity of Section~\ref{sec:complexity}, which the pathwise route does not give.}

\subsubsection*{Notation}
Throughout the manuscript, the following notations are used, unless expressly stated otherwise. Deterministic scalars are lowercase (e.g., $z, \theta$) and random scalars are uppercase (e.g., $\zrandom, \Theta$). Deterministic vectors are lowercase boldface (e.g., $\x\in\rd, d\geq 2$) and random vectors are uppercase boldface (e.g., $\xrandomvec$). Deterministic matrices are sans-serif uppercase (e.g., $\Q\in\rdp$) and random matrices are calligraphic uppercase boldface (e.g., $\Qrandom\in\rdp$). Script letters denote sets (e.g., $\mathscr{F}$). All random quantities live on the same probability space $(\Omega,\mathscr{F},\pr)$. 
Realizations of given random quantities $\zrandom$, $\xrandomvec$ and $\Qrandom$ are respectively denoted by $z:=\zrandom(\omega)$, $\x:=\xrandomvec(\omega)$ and $\Q:=\Qrandom(\omega)$, for some fixed $\omega\in \Omega$.
For two $\sigma$-algebras $\mathscr G$ and $\mathscr H$ on $\Omega$, $\mathscr G\vee \mathscr H:=\sigma(\mathscr G\cup \mathscr H)$ denotes the smallest $\sigma$-algebra generated by all events belonging to either $\mathscr G$ or $\mathscr H$, and therefore containing both $\mathscr G$ and~$\mathscr H$.
An increasing sequence $\accolade{\mathscr{G}_k}$ of $\sigma$-algebras of~$\mathscr{F}$ is called a {\it filtration}. $\accolade{\hat{\eb}^i}_{i=1}^p$ is the standard basis of $\rp$ and $\intbracket{a,b}:=[a,b]\cap\Z$ with $a\leq b$.

\section{General supermartingale-based framework for complexity analyses}\label{sec:SupermartFram}

{\bl
This section develops a general framework for analyzing a family of stopping times associated with a positive stochastic process~$\accoladekinN{\Stepsizek}$. The objective is to establish conditions under which these stopping times and certain accumulated nonnegative quantities admit bounds, both almost surely and in expectation, independently of the particular mechanism generating the process. We first exploit the interplay between the probabilistic behavior of successive values of~$\Stepsizek$ and the summability of their squares to derive stopping time bounds. We then extend the analysis to cumulative sums involving inverse powers of~$\Stepsizek$, which require additional moment control. The resulting framework will subsequently be applied to the algorithms introduced in Section~\ref{sec:Section3} by verifying that their associated stochastic processes satisfy Assumptions~\ref{Assumption:general1} and~\ref{Assumption:cr}.  Because the analysis is formulated only in terms of generic stochastic processes and stopping times, it is intended to be applicable more broadly to other adaptive stochastic optimization schemes, including direct-search~\cite{audet2019stomads,dzahini2020expected,DzaWildSubDir24,rinaldi2024stochastic}, line-search methods~\cite{paquette2018stochastic}, and adaptive optimization algorithms with stochastic oracles~\cite{JinSchXieAOA25}. 

We begin by recalling the definition of a stopping time for a discrete-time stochastic process.
}
\SSaddnew{\begin{definition}
A random variable $T$ with values in $\N\cup\{\infty\}$ is a stopping time with respect to $\{\gkalgebra\}_{k\in\N}$ if $\{T\le k\}\in\gkalgebra$ for every $k\in\N$.
\end{definition}}

{\bl The following assumptions \footnote{ \SSaddnew{ How Assumption~\ref{Assumption:general1} is satisfied by the algorithms of Section~\ref{sec:Section3} is shown in Remark~\ref{Remark:AssumptVerif}, Lemma~\ref{VilleLemmaResult} and the paragraph below~\eqref{stoppingTime}.}} inspired by existing works of~\cite{audet2019stomads,BaScVi2014,blanchet2016convergence,dzahini2020expected,DzaWildSubDir24,ha2025complexity,paquette2018stochastic,rinaldi2024stochastic,Shashaani2016Thesis,shashaani2018astro}, characterize the probabilistic behavior of the quantities $\accolade{\gkalgebra}_{k\in\N}$, $\accoladekinN{\Stepsizek}$, and $\accolade{\Tepsilon}_{\varepsilon>0}$ introduced below and required in the subsequent analysis.}

\SSaddnew{\begin{assumption}\label{Assumption:general1}
$\{\gkalgebra\}_{k\in\N}$ is a filtration of $\mathscr{F}${\bl , with respect to which} $\{\Tepsilon\}_{\varepsilon>0}$ is a family of stopping times; $\{\Stepsizek\}_{k\in\N}$ is a sequence of $\gkalgebra$-measurable random variables with $\Stepsize_0=\alpha_0$ deterministic; and $\{\Dscr_k\}_{k\in\N}$ is a sequence of events with $\Dscr_k\in\mathscr G_k$. Moreover, there exist constants $k_0\in\N$, $\pi_\star\in(1/2,1)$, $\dmax>0$, $\gamma>1$, $\kappa_a>0$ and $\varepsilon_0>0$ such that
\begin{itemize}
\item[(i)] for all $k\geq k_0$, almost surely, 
\begin{equation}\label{eq3:AssumpGen}
\prob{\Dscr_k|\gkalgebra}\geq\pi_\star;
\end{equation}
\item[(ii)] for all $k\in\N$, 
\begin{equation}\label{eq4:AssumpGen}
0<\Stepsizek\leq\dmax\qquad\text{and}\qquad\Stepsize_{k+1}\geq\gamma^{-1}\Stepsizek;
\end{equation}
\item[(iii)] for every $\varepsilon\in(0,\varepsilon_0]$, where $\varepsilon_0\le\min\{1,\gamma^{-k_0}\alpha_0/\kappa_a\}$,\footnote{Since~\eqref{eq5:AssumpGen} remains true when $\varepsilon_0$ is decreased, this restriction entails no loss of generality. It is used {\bl to establish \SSaddnew{Theorem~\ref{theor:GenIterCompl}}} through $\varepsilon\le1$, and in Lemma~\ref{lemma:br} to ensure $\Stepsize_{k_0}\ge\kappa_a\bar\varepsilon_j$ for every admissible~$j$.} and every $k\geq k_0$,
\begin{equation}\label{eq5:AssumpGen}
\mathscr{E}_{\kappa_a,\varepsilon}:=\accolade{k<T_{\varepsilon/2},\ \Stepsizek<\kappa_a\varepsilon}\cap\Dscr_k\subseteq\accolade{\Stepsize_{k+1}=\gamma\Stepsizek};
\end{equation}
\item[(iv)] $S_\Stepsize:=\sum_{k=0}^{\infty}\Stepsizek^2$ satisfies 
\begin{equation}\label{eq:SDeltaFinite}
\E{S_\Stepsize}<\infty \quad{\bl \text{which implies that}}\quad S_\Stepsize<\infty\quad\text{almost surely};
\end{equation}
\item[(v)] 
$\accolade{\Tepsilon}_{\varepsilon>0}$ is monotone in the sense that 
\begin{equation}\label{eq:MonotoneTimes}
0<\varepsilon_1\leq \varepsilon_2\qquad\Longrightarrow\qquad T_{\varepsilon_1}\geq T_{\varepsilon_2}\qquad\text{almost surely}.
\end{equation}
\end{itemize}
\end{assumption}}

{\bl 
Throughout this section, we note that the stopping times~$\Tepsilon$ are not assumed a priori to be finite; their almost-sure finiteness will follow from \SSaddnew{Theorem~\ref{theor:GenIterCompl}} under the conditions of Assumption~\ref{Assumption:general1}.
The latter describes the behavior of~$\accoladekinN{\Stepsizek}$ through probabilistically favorable events~$\Dscr_k$, whose conditional probabilities are uniformly bounded below by $\pi_\star\in (1/2, 1)$. These events need not be independent, and their probabilities may vary with~$k$. The process~$\accoladekinN{\Stepsizek}$ is assumed positive and uniformly bounded above, while its decrease between consecutive indices is limited by~$\gamma^{-1}$. The key requirement is that, before~$T_{\varepsilon/2}$, whenever~$\Stepsizek$ falls below~$\kappa_a\varepsilon$, the occurrence of~$\Dscr_k$ forces $\Stepsize_{k+1}=\gamma\Stepsizek$. Thus as in~\cite{blanchet2016convergence} (see, e.g.,~\cite[Page~95]{blanchet2016convergence}), sufficiently small values of~$\Stepsizek$ are associated with a tendency to increase, without requiring this increase at every index. The expected square summability of~$\accoladekinN{\Stepsizek}$ provides complementary control over the frequency of sufficiently large values of~$\Stepsizek$, while the monotonicity of~$\accolade{\Tepsilon}_{\varepsilon>0}$ allows bounds established at selected thresholds to be extended to arbitrary deterministic sequences of thresholds tending to zero.}
\subsection{Square-summability and stopping time bound}\label{sec:Section2pt1}
{\bl The goal of this subsection is first to establish the almost-sure finiteness of the stopping times by means of Theorem~\ref{theor:GenIterCompl}, and then to bound~$\Tepsilon$ in terms of~$\varepsilon$ and $S_\Stepsize$ satisfying the condition~\emph{(iv)} of Assumption~\ref{Assumption:general1}. This condition directly controls the number of indices preceding~$\Tepsilon$ for which~$\Stepsizek$ exceeds a fixed multiple of~$\varepsilon$ as will be seen in the proof of Theorem~\ref{theor:IterComplex:Gen}, but not those for which it is small, since their contributions to~$S_\Stepsize$ may be arbitrarily small. The latter are controlled via Lemma~\ref{VilleLemmaResult:Gen} by exploiting the events~$\Dscr_k$ and their effect on successive values of~$\Stepsizek$. 
} 

\SSaddnew{\begin{theorem}\label{theor:GenIterCompl}
Let the conditions {\bl {(i)}, {(ii)} and {(v)}} of Assumption~\ref{Assumption:general1} hold. Let $\{a_\varepsilon\}_{\varepsilon\in(0,\varepsilon_0]}$ be positive thresholds, nondecreasing in $\varepsilon$, such that for every $\varepsilon\in(0,\varepsilon_0]$ and every $k\ge k_0$,
\begin{equation}\label{eq:sept}
\accolade{k<T_{\varepsilon/2},\,\Stepsizek<a_\varepsilon}\cap\Dscr_k\subseteq\{\Stepsize_{k+1}=\gamma\Stepsizek\},
\end{equation}
and let $h:(0,\dmax]\to(0,\infty)$ be nondecreasing with $\E{S_{\Stepsize,h}}<\infty$, where $S_{\Stepsize,h}:=\sum_{k\ge0}h(\Stepsizek)$. Then, for every $\varepsilon\in(0,\varepsilon_0]$,
\begin{equation}\label{eq:ExpIterGen}
\E{\Tepsilon}\le\E{T_{\varepsilon/2}}\le k_0+\frac{\kappa_L}{2\pi_\star-1}+\frac{2\pi_\star}{2\pi_\star-1}\,\frac{\E{S_{\Stepsize,h}}}{h(a_\varepsilon)}<\infty,\qquad\kappa_L:=k_0+\log_\gamma\!\Bigl(\frac{\dmax}{\alpha_0}\Bigr),
\end{equation}
and $\Tepsilon<\infty$ almost surely for every $\varepsilon>0$. In particular, under Assumption~\ref{Assumption:general1}, the choices $a_\varepsilon=\kappa_a\varepsilon$ and $h(t)=t^2$ give, for every $\varepsilon\in(0,\varepsilon_0]$,
\begin{equation}\label{eq:ExpIter}
\E{\Tepsilon}\le\E{T_{\varepsilon/2}}\le\kappa_T\,\varepsilon^{-2},\qquad\kappa_T:=k_0+\frac{\kappa_L+2\pi_\star\kappa_a^{-2}\E{S_\Stepsize}}{2\pi_\star-1}.
\end{equation}
\end{theorem}
The proof is given in Appendix~\ref{app:B-theor:GenIterCompl}.}

\SSaddnew{\begin{remark}[relation to~\cite{blanchet2016convergence}]\label{rem:Blanchet}
Theorem~\ref{theor:GenIterCompl} is the counterpart of the expected stopping-time bound of~\cite{blanchet2016convergence}, which is also stated through a nondecreasing function of the step-size parameter. Condition~\eqref{eq:sept} is the forcing property of~\cite[Assumption~1(ii)]{blanchet2016convergence}, that the parameter tends to grow whenever it is small, with the favorable event made explicit and its conditional probability bounded below by~\eqref{eq3:AssumpGen}. Where~\cite{blanchet2016convergence} assumes a drift condition on a potential function and bounds $\E{\Tepsilon}$ by a renewal-reward argument, we assume only $\E{S_{\Stepsize,h}}<\infty$, which the drift provides in Section~\ref{sec:convAnalysis} (see~\eqref{fourpttwelve}) but which is weaker, and count the indices directly over a bounded horizon. The same count, carried out pathwise, gives the almost-sure bounds of Theorems~\ref{theor:IterComplex:Gen} and~\ref{Theorem:WorkComplexity:Gen} and the inverse-square sums behind the sample complexity bounds, which~\cite{blanchet2016convergence} does not provide; this is why the result is proved here rather than imported. The remaining results of this section are stated for $h(t)=t^2$ and $a_\varepsilon=\kappa_a\varepsilon$, the case the algorithms of Section~\ref{sec:Section3} require.
\end{remark}}

{\bl  
\SSaddnew{Lemma~\ref{VilleLemmaResult:Gen} repeats the counting argument of the proof of Theorem~\ref{theor:GenIterCompl} pathwise, along dyadic thresholds.} 
Indeed, before~$T_{\varepsilon/2}$, whenever~$\Stepsizek<\kappa_a\varepsilon$, the occurrence of~$\Dscr_k$ forces $\Stepsize_{k+1}=\gamma\Stepsizek$, whereas $\Stepsize_{k+1}\geq\gamma^{-1}\Stepsizek$ always holds. Thus, on the logarithmic scale, favorable outcomes at small values yield increments of~$+1$, while all other increments are bounded below by~$-1$. Since $A_k$ is uniformly bounded above, this controls the number of indices with small values in terms of those with larger values, up to the fluctuations of the number of favorable events around its conditional mean, which are handled by a martingale. Lemma~\ref{VilleLemmaResult:Gen} formalizes this argument along dyadic thresholds and extends the resulting bound to arbitrary deterministic sequences $\accolade{\vepsell}_{\ell\geq 1}$ tending to zero.}
\begin{lemma}\label{VilleLemmaResult:Gen}
Let Assumption~\ref{Assumption:general1} hold.
Let $\accolade{\vepsell}_{\ell\geq 1}$ be any deterministic positive sequence such that $\vepsell\downarrow 0$ as $\ell\to\infty$. Then 
we have almost surely,
\begin{equation*}
\underset{\ell\to\infty}{\limsup}\ \vepsell^2\abs{\accolade{0\leq k<\taurandom_{\vepsell}:\Stepsizek<\kappa_a\vepsell}}\leq \frac{S_{\Stepsize}}{(2\pi_\star-1)\kappa_a^2}.
\end{equation*}
\end{lemma}

\SSadd{The proof is given in Appendix~\ref{app:B-VilleLemmaResult:Gen}.}

{\bl Lemma~\ref{VilleLemmaResult:Gen} controls asymptotically the number of indices~$k<T_{\varepsilon_\ell}$ for which $\Stepsizek<\kappa_a\varepsilon_\ell$. The remaining indices are controlled by $S_\Stepsize=\sum_{k=0}^{\infty}\Stepsizek^2$, since each contributes at least $\kappa_a^2\varepsilon_\ell^2$ to this sum. Combining these two bounds yields the main result of this subsection stated next, valid eventually almost surely along any deterministic sequence of $\vepsell\downarrow 0$, with an integrable random coefficient expressed in terms of~$S_\Stepsize$, and independent of~$\accolade{\vepsell}_{\ell\geq 1}$.}
\begin{theorem}\label{theor:IterComplex:Gen}
Let Assumption~\ref{Assumption:general1} hold. 
Consider the almost surely finite random variable $R_\Stepsize:=1+\left(1+\frac{1}{2\SSaddnew{\pi_\star}-1}\right)\frac{S_\Stepsize}{\kappa_a^2}$, with $\E{R_\Stepsize}<\infty$.
For every deterministic sequence of $\vepsell\to 0$, 
\begin{equation}\label{complRes1:Gen}
\prob{\taurandom_{\vepsell}\vepsell^2>R_\Stepsize\ \ i.o.}=0.
\end{equation}
In other words, along every deterministic threshold sequence of $\vepsell\downarrow 0$, $\taurandom_{\vepsell}=\mathcal{O}(\vepsell^{-2})$ with an integrable random constant.
\end{theorem}

\SSadd{The proof is given in Appendix~\ref{app:B-theor:IterComplex:Gen}.}

\subsection{Inverse square sums and reciprocal moment control}
{\bl We now extend the preceding stopping time analysis to cumulative sums involving~$\Stepsizek^{-2}$. Although Theorem~\ref{theor:IterComplex:Gen} provides an eventual almost-sure bound on~$T_{\vepsell}$, and hence on the number of indices preceding it, it does not control the accumulated inverse squares of~$\Stepsizek$, which may become arbitrarily large when~$\Stepsizek$ is small. Consequently, the counting argument of Section~\ref{sec:Section2pt1} is insufficient. To address this difficulty, we establish moment bounds for appropriately normalized inverse squares of~$\Stepsizek$ under the following additional assumption.}
\begin{assumption}\label{Assumption:cr}
Let $\epsLanbda\in (0,1)$. There exists $r>\frac{1}{1-\epsLanbda}$ such that\footnote{A discussion on how this assumption is satisfied in the proposed algorithmic framework, is provided later in Remark~\ref{Remark:Parameters1}.} $c_r:=\pi_\star\gamma^{-2r}+(1-\pi_\star)\gamma^{2r}<1$.
\end{assumption}
{\bl Assumption~\ref{Assumption:cr} strengthens the probabilistic requirement on the favorable events~$\Dscr_k$ to control inverse moments of~$\Stepsizek$. Indeed for~$j\in\N$ sufficiently large, whenever~$k<T_{2^{-j-1}}$ and~$\Stepsizek<\kappa_a2^{-j}$, the occurrence of~$\Dscr_k$ forces~$\Stepsize_{k+1}=\gamma\Stepsizek$, reducing~$\Stepsize_k^{-2r}$ by the factor~$\gamma^{-2r}$, whereas its increase is otherwise bounded by~$\gamma^{2r}$. The condition~$c_r<1$ thus ensures via Lemma~\ref{lemma:br}-\eqref{ineq:br} a contraction in the corresponding conditional moment bound, up to an additive constant. Moreover, the restriction~$r(1-\epsilon_\lambda)>1$ later guarantees in Lemma~\ref{lemma:InverseRadius} -\eqref{eq:newarbound} the summability of probability bounds obtained through Markov's inequality. These properties allow Lemma~\ref{lemma:br} to establish uniform moment bounds for the normalized inverse-square process~$\Zkj$ defined below, and its truncated version~$\Xkj$.}
\begin{lemma}\label{lemma:br}
Let Assumptions~\ref{Assumption:general1}-(i)-(iii) and~\ref{Assumption:cr} hold. Consider for all $j\in \N$ the random processes defined for all $k\geq k_0$ by $\Zkj:=\max\accolade{1,\left(\frac{\kappa_a \SSadd{\bar\varepsilon_j}}{\Stepsizek}\right)^2}$ and $\Xkj:=\Zkj\mathds{1}_{\accolade{k<\taurandom_{\SSadd{\bar\varepsilon_j/2}}}}$ \SSadd{using $\bar\varepsilon_j=2^{-j}$}.
Then the following hold for 
\SSaddnew{every $j$ with $\bar\varepsilon_j\le\varepsilon_0$
 (we call such $j$ admissible) 
}
with $c_r:=\pi_\star\gamma^{-2r}+(1-\pi_\star)\gamma^{2r}$, $b_r:=1+\gamma^{2r}$ and $a_r:=1+\frac{b_r}{1-c_r}$.
\begin{equation}\label{ineq:br}
\SSaddnew{  \E{\left.\left(\XkPlusOnej\right)^r\right\rvert\gkalgebra}\leq c_r\left(\Xkj\right)^r+b_r\mathds{1}_{\accolade{k<\taurandom_{\bar\varepsilon_j/2}}},}
\end{equation}

\begin{equation}\label{ineq:ar}
\underset{k\geq k_0}{\sup}\E{\left(\Xkj\right)^r}\leq a_r,
\end{equation}
and for $l> k_0$,
\begin{equation}\label{ineq:Jensen}
\E{\left(\sum_{k=k_0}^{l-1}\Xkj\right)^r}\leq a_r l^r.
\end{equation}
\end{lemma}

\SSadd{The proof is given in Appendix~\ref{app:B-lemma:br}.}

{\bl Lemma~\ref{lemma:br} established uniform moment bounds for~$\Xkj$ and its cumulative values over deterministic index ranges. In the next lemma, we apply Markov's inequality to these sums, obtaining in~\eqref{eq:newarbound} probability bounds of order~$j^{-r(1-\epsilon_\lambda)}$, which are summable under Assumption~\ref{Assumption:cr}. Combining the resulting almost-sure estimates with the stopping time bound of Theorem~\ref{theor:IterComplex:Gen} allows the deterministic index ranges to be extended to~$T_{2^{-j-1}}$. Finally, using the definition of~$\Xkj$, Lemma~\ref{lemma:InverseRadius} establishes an eventual almost sure bound of order~$2^{4j}j^{1-\epsilon_\lambda}$ on the accumulated inverse squares of~$\Stepsizek$.}
\begin{lemma}\label{lemma:InverseRadius}
Let Assumptions~\ref{Assumption:general1} and~\ref{Assumption:cr} hold. Consider the random variable $\Cdelta:=\kappa_a^{-2}\left(4\ceil{R_\Stepsize}+5\right)<\infty$ almost surely.
Then, for all sufficiently large~$j\in\N$, \SSadd{and $\bar \varepsilon_j=2^{-j}$}
\begin{equation*}
\sum_{k=k_0}^{\taurandom_{\SSadd{\bar \varepsilon_j/2}}-1}\Stepsizek^{-2}\leq \Cdelta 2^{4j}j^{1-\epsilon_\lambda}.
\end{equation*}
\end{lemma}

\SSadd{The proof is given in Appendix~\ref{app:B-lemma:InverseRadius}.}

{\bl We now consider a general sequence $\accoladekinN{\Wk}$ of nonnegative random variables, assumed to be eventually bounded, almost surely, by a finite random multiple of $\lambda_k\Stepsizek^{-2}$, with~$\lambda_k$ specified below. This condition allows the extension of the preceding analysis to the accumulated quantities $W_\varepsilon:=\sum_{k=0}^{\taurandom_\varepsilon-1}\Wk$, without specifying how~$\Wk$ is generated. Combining Theorem~\ref{theor:IterComplex:Gen} with Lemma~\ref{lemma:InverseRadius} then yields the eventual almost-sure bound of order $\varepsilon^{-4}(\log(1/\varepsilon))^2$ for the accumulated quantities~$W_\varepsilon$, established in the following theorem, where the logarithmic exponent results from combining the growth of~$\lambda_k$ with the bound of Lemma~\ref{lemma:InverseRadius}.}
\begin{theorem}\label{Theorem:WorkComplexity:Gen}
Let Assumptions~\ref{Assumption:general1} and~\ref{Assumption:cr} hold. 
\SSadd{Let $R_\Stepsize$ and $\Cdelta$ be defined as in Theorem~\ref{theor:IterComplex:Gen} and Lemma~\ref{lemma:InverseRadius}.} 
Let $\accoladekinN{\lambda_k}$ be the deterministic sequence defined for all $k\in\N$ by
\begin{equation*}
\lambda_k=\ubar{\lambda}(\log(k+2))^{1+\epsilon_\lambda},\qquad \epsLanbda\in(0,1),\ \, \SSaddnew{\ubar{\lambda}> 0},
\end{equation*}
and let $\accoladekinN{\Wk}$ be a sequence of nonnegative random variables such that for some almost surely finite random quantities $\ubarlong{K}_W>0$ and $N_W>0$, 
\begin{equation}\label{eq:WkBound1}
\Wk\leq N_W\lambda_k\Stepsizek^{-2}\qquad\text{for all}\ k\geq \ubarlong{K}_W,\ \text{almost surely}.
\end{equation}
Define for all $\varepsilon>0$,
\begin{equation*}
W_\varepsilon:=\sum_{k=0}^{\taurandom_\varepsilon-1}\Wk.
\end{equation*}
Then for every deterministic sequence of $\vepsell\downarrow 0$, there exists an almost surely finite random variable $R_W>0$ depending on $\ubarlong{K}_W, N_W, C_\Stepsize, \epsLanbda$\SSaddnew{,$R_A, \underline\lambda$} and $k_0$, such that 
\begin{equation*}
\prob{W_{\vepsell}>R_W\vepsell^{-4}\left(\log(1/\vepsell)\right)^2\ \, i.o.}=0.
\end{equation*}
Equivalently, $W_{\vepsell}=\mathcal{O}\left(\vepsell^{-4}\left(\log(1/\vepsell)\right)^2\right)$ eventually almost surely along every deterministic sequence of $\vepsell\downarrow 0$.
\end{theorem}

\SSadd{The proof is given in Appendix~\ref{app:B-Theorem:WorkComplexity:Gen}.}

{\bl
The preceding results provide almost-sure control of the accumulated quantities~$W_\varepsilon$. We now derive their counterparts in expectation. To this end, Lemma~\ref{lemma:ExpInverseRadius} first establishes deterministic bounds on the expected inverse-square sums~$\sum_{k<T_\varepsilon}\Stepsizek^{-2}$ and, under the logarithmic weighting~$\lambda_k$, on~$\sum_{k<T_\varepsilon}\lambda_k\Stepsizek^{-2}$. These estimates then allow a conditional bound of the form $\E{\Wk|\gkalgebra}\leq \kappa_W\lambda_k\Stepsizek^{-2}$ to be converted directly into an expected bound for~$W_\varepsilon$, yielding Theorem~\ref{Theorem:ExpectedWork:Gen}.
}
\SSaddnew{\begin{lemma}\label{lemma:ExpInverseRadius}
Let Assumptions~\ref{Assumption:general1} and~\ref{Assumption:cr} hold and let $\kappa_T$ be as in {\bl \SSaddnew{Theorem~\ref{theor:GenIterCompl}}}. Then for every $\varepsilon\in(0,\varepsilon_0/2]$:
\begin{enumerate}
\item[\textup{(a)}] $\displaystyle\E{\sum_{k=0}^{\taurandom_\varepsilon-1}\Stepsize_k^{-2}}\le\kappa_\Stepsize\,\varepsilon^{-4}$, \ with $\displaystyle\kappa_\Stepsize:=\frac{k_0\gamma^{2k_0}}{\alpha_0^2}+\frac{1+b_r\kappa_T}{(1-c_r)\kappa_a^2}$;
\item[\textup{(b)}] if $\lambda_k=\ubar{\lambda}(\log(k+2))^{1+\epsLanbda}$ and $\E{u(S_\Stepsize)}<\infty$, where $u(t):=t\,(\log(t+3))^{1+\epsLanbda}$, then
\[
\E{\sum_{k=0}^{\taurandom_\varepsilon-1}\lambda_k\Stepsize_k^{-2}}\le\kappa_{\Stepsize,\lambda}\,\varepsilon^{-4}\bigl(1+\log(1/\varepsilon)\bigr)^{1+\epsLanbda},
\]
where $\kappa_{\Stepsize,\lambda}$ is a deterministic constant depending only on the constants of Assumptions~\ref{Assumption:general1} and~\ref{Assumption:cr}, on 
$\epsLanbda$, $\ubar{\lambda}$, and on $\E{S_\Stepsize}$ and $\E{u(S_\Stepsize)}$.
\end{enumerate}
\end{lemma}
The proof is given in Appendix~\ref{app:B-lemma:ExpInverseRadius}.
}

\SSaddnew{\begin{theorem}\label{Theorem:ExpectedWork:Gen}
Let Assumptions~\ref{Assumption:general1} and~\ref{Assumption:cr} hold and let $\{\Wk\}_{k\in\N}$ be nonnegative random variables with $\E{\Wk|\gkalgebra}\le\kappa_W\lambda_k\Stepsize_k^{-2}$ almost surely for every $k$, for a deterministic $\kappa_W$, and let $W_\varepsilon:=\sum_{k=0}^{\taurandom_\varepsilon-1}\Wk$. If $\lambda_k=\ubar\lambda(\log(k+2))^{1+\epsLanbda}$ and $\E{u(S_\Stepsize)}<\infty$, then for every $\varepsilon\in(0,\varepsilon_0/2]$,
\[
\E{W_\varepsilon}\le\kappa_W\,\kappa_{\Stepsize,\lambda}\,\varepsilon^{-4}\bigl(1+\log(1/\varepsilon)\bigr)^{1+\epsLanbda}.
\]
If instead $\lambda_k\equiv\ubar\lambda$, then $\E{W_\varepsilon}\le\kappa_W\,\ubar\lambda\,\kappa_\Stepsize\,\varepsilon^{-4}$ for every $\varepsilon\in(0,\varepsilon_0/2]$, without the moment condition on $S_\Stepsize$.
\end{theorem}
\begin{proof}
Since $\{k<\taurandom_\varepsilon\}\in\gkalgebra$ and $\Wk\ge0$, Tonelli's theorem and the tower property give $$\E{W_\varepsilon}=\sum_{k\ge0}\E{\mathds 1_{\{k<\taurandom_\varepsilon\}}\E{\Wk|\gkalgebra}}\le\kappa_W\,\E{\sum_{k<\taurandom_\varepsilon}\lambda_k\Stepsize_k^{-2}},$$ and Lemma~\ref{lemma:ExpInverseRadius}(b), respectively~(a), bounds the right-hand side.
\end{proof}}

\SSaddnew{While the integrability of $R_A$ in Theorem~\ref{theor:IterComplex:Gen} follows from Assumption~\ref{Assumption:general1}-\emph{(iv)}, Theorem~\ref{Theorem:WorkComplexity:Gen} guarantees only the almost-sure finiteness of $R_W$. Theorem~\ref{Theorem:ExpectedWork:Gen} has deterministic constants instead, at the price of a conditional moment bound on {\bl $\Wk$} and, under the logarithmic schedule, of the moment condition $\E{u(S_A)}<\infty$.}

\SSaddnew{\begin{remark}[summary of section]
The almost-sure results of this section, Lemma~\ref{VilleLemmaResult:Gen}, Theorems~\ref{theor:IterComplex:Gen} and~\ref{Theorem:WorkComplexity:Gen} and the inverse-square bounds of Lemmas~\ref{lemma:br}--\ref{lemma:InverseRadius}, carry the complexity analysis of ASTRO-DF in~\cite{ha2025complexity} into the present setting: there the accuracy events hold eventually along each sample path, here the favorable events $\Dscr_k$ hold at every index with conditional probability at least $\pi_\star>\frac12$, and the indices before $T_\varepsilon$ are again counted by the size of $A_k$, with $S_A$ controlling the large values and the forcing property controlling the small ones. The expected results here are different from those  of~\cite{blanchet2016convergence} in that  Theorem~\ref{theor:GenIterCompl} counts over a bounded horizon $T_\varepsilon\wedge\bar k$, where the conditional probability of $\Dscr_k$ and $\E{S_A}<\infty$ give $\E{T_\varepsilon}=\mathcal O(\varepsilon^{-2})$ in a few lines, with no drift condition on a potential and no exponential supermartingale. The almost-sure count of Lemma~\ref{VilleLemmaResult:Gen} is then a consequence, because the fluctuation of the number of favorable events around its conditional mean is a bounded-increment martingale whose variance is at most $\E{T_\varepsilon}$. For the accumulated quantities, Lemma~\ref{lemma:ExpInverseRadius} sums a conditional moment contraction for the normalized inverse squared $A_k$ up to the stopping time, and Theorem~\ref{Theorem:ExpectedWork:Gen} converts a conditional moment bound on {\bl $\Wk$} into $\E{W_\varepsilon}=\mathcal O(\varepsilon^{-4})$ up to a logarithm. Remark~\ref{rem:Blanchet} compares the assumptions with those of~\cite{blanchet2016convergence}.
\end{remark}}

{\bl The applicability of the main results derived in this section to a particular method depends on verifying the assumptions governing~$\accoladekinN{\Stepsizek}$, the stopping times, and the cumulative sums. These verifications will be carried out for the algorithms analyzed in the following sections.}
\section{Algorithmic framework and resulting stochastic process}\label{sec:Section3}
{\bl 
This section introduces the algorithmic framework and the stochastic process generated by the two variants analyzed in the remainder of the paper. Section~\ref{sec:theMethods} presents the common adaptive sampling trust-region framework in random subspaces, specifies the model decrease and update mechanisms, and distinguishes the proposed \ASTRAF\ and \ASTRAR\ algorithms through their respective implementations of the estimate computation step. Section~\ref{sec:estiModels} then makes the estimators, models, adaptive sample sizes, and filtrations precise, thereby formalizing the associated stochastic process and establishing the almost-sure finiteness of the adaptive sample sizes.
}
\subsection{The adaptive sampling stochastic trust-region method}\label{sec:theMethods}
Given a current iterate $\xzerok$ at iteration~$k$, an affine subspace $\YkSet:=\accolade{\xzerok+\Qk\s:\s\in\rp}\subseteq\rd$ is randomly chosen through the range of a matrix $\Qk\in\rdp$ with random selected entries. Then making use of realizations of the noisy function $F(\cdot;\Theta)$, a subspace local model $\mhatk:\rp\to\R$ is built on $\YkSet$ to approximate $f(\x)$ in $\mathscr{B}(\xzerok,\dk;\Qk):=\accolade{\xzerok+\Qk\s\in\YkSet:\normii{\s}\leq\dk}$. The quadratic subspace model employed in this manuscript is given by 
\[f(\xzerok+\Qk\s)\approx\mhatk(\s):=\fbarok+\gkhat^{\top}\s+\frac{1}{2}\s^{\top}\Hkhat\s, \]
where $\fbarok\approx f(\xzerok)$ while $\gkhat\in\rp$ and $\Hkhat\in\rpp$ are low-dimensional gradient and Hessian, respectively. By using the solution (or {\it trial step}) $\sk\approx\argmin\accolade{\mhatk(\s):\s\in\rp,\normii{\s}\leq \dk}$ resulting from an approximate minimization of~$\mhatk$ in $\mathscr{B}(\bm{0},\dk)$, a {\it tentative step} $\Qk\sk\in\rd$ is produced. Inspired by~\cite{CRsubspace2021,chen2018stochastic,DzaWildSub2022}, the trial step is required to provide a sufficient decrease in the model by satisfying the following condition, so-called standard fraction of Cauchy decrease.
 
\begin{assumption}\label{CauchyDecreaseAssump}
A trial step $\sk\in\rp$ is computed at each iteration~$k$ so that 
\begin{equation}\label{CauchyDecreaseEq}
\mhatk(\bm{0})-\mhatk(\sk)\geq \frac{\kfcd}{2}\normii{\gkhat}\min\accolade{\dk,\frac{\normii{\gkhat}}{\max(\normeshort{\Hkhat},1)}}
\end{equation}
for some constant $\kfcd\in (0,1]$.
\end{assumption}
\SSaddnew{The trial step is computed from the model and the radius by a deterministic rule: $\sk$ is a fixed Borel-measurable function of the model coefficients $(\fbarok,\gkhat,\Hkhat)$ and of $\dk$. The oracle $(\x,\theta)\mapsto F(\x;\theta)$ is jointly measurable. These two properties make each iteration a measurable function of the quantities drawn, which Section~\ref{sec:estiModels} records in the appropriate $\sigma$-algebras.} 

Estimates $\fbarokcur\approx f(\xzerok)$ and $\fbarsk\approx f(\xks)$, where $\xks:=\xzerok+\Qk\sk$ is a trial point, are constructed using realizations of the noisy function $F(\cdot,\Theta)$ and the random matrix $\QkRandom\sim\mathbb{Q}_k$. 
\SSadd{Step~\textbf{[3]} of Algorithm~\ref{algoAdapTR} is left generic and is instantiated in two different ways, leading to the two algorithms analyzed side by side in this manuscript. {\bl Inspired by~\cite{chen2018stochastic,DzaWildSub2022},} Algorithm~\ref{algoF}, {\bl referred to as} \ASTRAF, recomputes both estimates using an independent copy of the replication stream, as required by a pointwise analysis of the trial estimate. {\bl Inspired by~\cite{ha2025complexity,shashaani2018astro},} Algorithm~\ref{algoR}, {\bl referred to as} \ASTRAR, instead sets $\fbarokcur=\fbarok$ and evaluates the trial point using the replications already drawn. All subsequent statements concern Algorithms~\ref{algoF} and~\ref{algoR}; Algorithm~\ref{algoAdapTR} serves only to specify the components they share.}
 
\begin{algorithm}[htb]
\caption{Adaptive sampling stochastic trust-region algorithm in random subspaces}
\label{algoAdapTR} 
\textbf{[0] Initialization}\\
\hspace*{5mm}Choose expansion and shrinkage constant $\gamma>1$, 
criticality threshold $\mu>0$, maximum and initial radii \\ 
\hspace*{5mm}$\delta_{\max}>\delta_0>0$, success threshold $\eta\in (0,1)$, sample size inflation sequence $\accoladekinN{\lambda_k}$, starting point \\ 
\hspace*{5mm}$\x_0^0\in \rd$, adaptive sampling constant $\sigma_0, \kappaas>0$, well alignment parameters $\SSaddnew{\alpha_Q,\beta_Q}\in (0,1)$, and subspace\\ 
\hspace*{5mm}dimension $p\in\intbracket{1,d}$.\\ 
\hspace*{5mm}Set the iteration counter $k \gets 0$ and $q=2p$.\\
\textbf{[1] Construction of subspace local model}\\
\hspace*{5mm}Generate $\Qk\in\rdp$, a realization of a $WAM(\SSaddnew{\alpha_Q,\beta_Q})$ using a distribution $\mathbb{Q}_k$.\\
\hspace*{5mm}Select a design set $\mathscr{S}_k:=\accolade{\szero_k,\sone_k,\dots,\s^q_k}\subset\rp\cap\mathscr{B}(\bm{0},\dk)$, with $\szero_k=\bm{0}$, $\s^j_k=\dk\ebhat^j$ and $\s^{p+j}_k=-\dk\ebhat^j$,\\ \hspace*{5mm}$j\in\intbracket{1,p}$, and compute estimates $\fbarik(\sii_k,\nk)\approx f(\xzerok+\Qk\sii_k)$ with a sample size $\nk(\sigma_0,\kappaas)$, $i\in\intbracket{0,q}$,\\ 
\hspace*{5mm}to build a subspace model $\mhatk:\rp\to\R$ following Section~\ref{sec:estiModels}.\\
\textbf{[2] Step calculation}\\
\hspace*{5mm}Compute $\sk\approx \argmin\accolade{\mhatk(\s):\s\in\rp, \normii{\s}\leq \dk}$ and set $\xks=\xzerok+\Qk\sk$. \\
\textbf{[3] Estimates computation}\\
\hspace*{5mm}\SSadd{Obtain estimates $\fbarokcur\approx f(\xzerok)$ and $\fbarsk\approx f(\xks)$.}\\
\textbf{[4] Updates}\\
\hspace*{5mm}Compute $\hat{\rho}_k:=\frac{\SSaddnew{\fbarsk-\fbarokcur}
}{\mhatk(\sk)-\mhatk(\bm{0})}$ \SSaddnew{with $\hat\rho_k:=0$ if $\mhatk(\sk)=\mhatk(\bm 0)$}.  \\
\hspace*{5mm}If $\hat{\rho}_k\SSaddnew{\ge}\eta$ and $\mu\normii{\gkhat}\geq \dk$ (\textbf{success}): set $\xkun=\xks$ and $\dkun=\min\accolade{\gamma\dk,\SSaddnew{\delta_{\max}}}$.\\
\hspace*{5mm}Otherwise (\textbf{failure}): set $\xkun=\xzerok$ and $\dkun=\gamma^{-1}\dk$.\\
\hspace*{5mm}Update the iteration counter $k\gets k+1$ and go to \textbf{[1]}.
\end{algorithm}
\SSadd{\begin{algorithm}[htb]
\caption{\ASTRAF: Algorithm~\ref{algoAdapTR} with Step~\textbf{[3]} implemented on a fresh replication stream}\label{algoF}
{\bl Keep} Steps \textbf{[0]}, \textbf{[1]}, \textbf{[2]} and \textbf{[4]} {\bl of} Algorithm~\ref{algoAdapTR}, {\bl and update Step \textbf{[3]} as follows}.\\
\textbf{[3] Estimates computation (fresh stream)}\\
\hspace*{5mm}Draw an independent copy $\Theta^{0,s}_k$ of the replication stream. Obtain $\SSaddnew{\fbarokcur=}\fbarokcur\left(\bm{0},\nokbrak\right)$ and\\
\hspace*{5mm}$\SSaddnew{\fbarsk=}\fbarsk\left(\sk,\nskbrak\right)$ using {\bl a realization $\nsk$ of} the sample size $\nskrandom$ {\bl given by}~\eqref{n0sDefEq}.
\end{algorithm}}

\SSadd{\begin{algorithm}[htb]
\caption{\ASTRAR: Algorithm~\ref{algoAdapTR} with Step~\textbf{[3]} implemented by reuse}\label{algoR}
{\bl Keep} Steps \textbf{[0]}, \textbf{[1]}, \textbf{[2]} and \textbf{[4]} as in Algorithm~\ref{algoAdapTR}, {\bl and update Step \textbf{[3]} as follows}.\\
\textbf{[3] Estimates computation (reuse)}\\
\hspace*{5mm}Obtain $\SSaddnew{\fbarsk=}\fbarsk\left(\sk,\nk\right)$ on the samples $\Theta_{k,1},\dots,\Theta_{k,\nk}$ already drawn in Step \textbf{[1]}, and set\\
\hspace*{5mm}$\fbarokcur=\fbarok$, {\bl where $\nk$ is a realization of the sample size~$\nkrandom$ given by~\eqref{nikDefEq}.}
\end{algorithm}}

\SSaddnew{\paragraph{Roadmap of analyzing the algorithms.} Section~\ref{sec:estiModels} makes the sample size and the estimated quantities precise and shows that the adaptive sample sizes devised are almost surely finite. Then, the guarantees for Algorithms~\ref{algoF} and~\ref{algoR} can be obtained as follows: 
\begin{enumerate}
    \item We show, in Section~\ref{sec:Prob}, that an iteration is ``good'', where the model is accurate in its gradients and the estimates entering the acceptance test are accurate in their difference, with a probability close to one conditioned on the past. The two algorithms differ in proving this result based on how Step~\textbf{[3]} is implemented. 
    \item We combine, in Section~\ref{sec:convAnalysis}, the deterministic consequences of a good iteration (success when the radius is small relative to the gradient, otherwise a bounded increase of $f$ due to stochasticity), with a supermartingale argument on a Lyapunov function that mixes $f$ and $\Dk^2$, obtains $\E{\sum_k\Dk^2}<\infty$ and $\normii{\nabla f(\Xzerok)}\to0$ almost surely. 
    \item To prove the total sample complexity of both algorithms, in Section~\ref{sec:complexity}, we verify the assumptions of Section~\ref{sec:SupermartFram} that turns the results in 1 and 2 above into iteration bounds of order $\varepsilon^{-2}$, and combine them with moment bounds on the adaptive sample sizes to bound the number of oracle calls by order $\varepsilon^{-4}$ up to a logarithmic factor, almost surely and in expectation.
\end{enumerate}
}

\subsection{Estimates, model construction and stochastic process generated by \SSadd{Algorithms~\ref{algoF} and~\ref{algoR}}}\label{sec:estiModels}
 
At each iteration~$k$ of \SSadd{Algorithm~\ref{algoF} or~\ref{algoR}}, the local model~$\mhatk$ and estimates $\fbarokcur\approx f(\xzerok),\fbarsk\approx f(\xks)$  constructed using realizations of the noisy function $F(\cdot,\Theta)$ and the random matrix $\QkRandom\sim\mathbb{Q}_k$ as mentioned above, can be considered as realizations of random quantities $\Mhatk$, $\Fbarokcur$ and $\Fbarsk$, respectively, 
whose behavior influences each iteration of the algorithm. The latter therefore results in a stochastic process which is shown by the present work to have desirable convergence properties as long as the sequences $\accoladekinN{\Mhatk}$ and $\accoladekinN{\Fbarokcur,\Fbarsk}$ satisfy certain conditions, conditioned on the past history of the algorithm. 
 
Throughout the manuscript, 
\SSaddnew{the models are built on interpolation sets $\mathscr{S}_k$ of \textbf{Step~[1]} with $q:=2p$.} 
Let \SSaddnew{$\Theta_k:=(\Theta_{k,\ell})_{\ell\ge1}$ be a sequence of independent copies of $\Theta$.} 
To construct the model $\mhatk$,  estimates $\fbarik(\sii_k,\nk)\approx f(\xzerok+\Qk\sii_k)$ $i\in\intbracket{0,q}$, are obtained as realizations of $\Fbarik(\siirandom_k,\nkrandom;\Theta_k)$, where
\begin{equation}
\Fbarik(\siirandom_k,n;\Theta_k):=\frac{1}{n}\sum_{\ell=1}^{n}F(\X_k^i;\Theta_{k,\ell}), \quad i\in\intbracket{0,q},
\end{equation}
where $\X_k^i:=\Xzerok+\QkRandom\siirandom_k$. Consider the {\it slowly varying} sequence $\accoladekinN{\lambda_k}$ (see, e.g.,~\cite[Definition~2.1]{ha2025complexity}) defined for some constants $\SSaddnew{\ubar{\lambda}> 0}$ and $\epsilon_\lambda\in (0,1)$ by $\lambda_k=\ubar{\lambda}(\log(k+2))^{1+\epsilon_\lambda}$, for all~$k$, and let $\sigma_0,\kappaas>0$ be constants. Then the sample size $\nkrandom$ is chosen according to 
\begin{equation}\label{nikDefEq}
\nkrandom = \min\accolade{n\in\N_+:\frac{\sigma_0\vee\hat{\sigma}_F^{\max}(n;\Theta_k)}{\sqrt{n}}\leq \kappaas\frac{\Dk}{\sqrt{\lambda_k}}},
\end{equation}
where $\hat{\sigma}_F^{\max}(n;\Theta_k):=\max_{i\in\intbracket{0,q}} \hat{\sigma}_F(\siirandom_k,n;\Theta_k)$, with $\hat{\sigma}_F(\siirandom_k,n;\Theta_k)$ defined by the variance estimator given by
\begin{equation}\label{sigmaHatFsqEq}
\hat{\sigma}_F^2(\siirandom_k,n;\Theta_k):= \frac{1}{n}\sum_{\ell=1}^{n}\left(F(\Xzerok+\QkRandom\siirandom_k;\Theta_{k,\ell})-\Fbarik(\siirandom_k,n;\Theta_k)\right)^2, \quad i\in\intbracket{0,q}.
\end{equation}
Likewise, \SSadd{in \ASTRAF,} $\fbarokcur$ and $\fbarsk$ are respectively obtained as realizations of random estimates
\begin{equation}
\Fbarokcur(\bm{0},\nokrandombrak;\Theta_k^{0,s})=\frac{1}{\nokrandombrak}\sum_{\ell=1}^{\nokrandombrak}F(\Xzerok;\Theta_{k,\ell}^{0,s})\ \text{and}\  \Fbarsk(\srandom_k,\nskrandom;\Theta_k^{0,s})=\frac{1}{\nskrandom}\sum_{\ell=1}^{\nskrandom}F(\X_k^s;\Theta_{k,\ell}^{0,s}),
\end{equation}
where $\X_k^s:=\Xzerok+\QkRandom\srandom_k$ and the $\Theta_{k,\ell}^{0,s}$ are independent random samples of an independent copy $\Theta_k^{0,s}$ of $\Theta_k$. 
\SSadd{
The sample size of the fresh \SSaddnew{independent copy $\Theta^{0,s}_k$} is
} 
\begin{equation}\label{n0sDefEq}
\nskrandom= \min\accolade{n\in\N_+:\frac{\max\accolade{\sigma_0,\hat{\sigma}_F^{[0]}(n;\Theta_k), \hat{\sigma}_F^{[s]}(n;\Theta_k)} }{\sqrt{n}}\leq \kappaas\frac{\Dk}{\sqrt{\lambda_k}}},
\end{equation}
with 
\begin{equation}\label{sigm0Fdef}
\left(\hat{\sigma}_F^{[0]}(n;\Theta_k)\right)^2:= \frac{1}{n}\sum_{\ell=1}^{n}\left(F(\Xzerok;\Theta_{k,\ell}^{0,s})-\Fbarokcur(\bm{0},n;\Theta_k^{0,s})\right)^2
\end{equation}
and 
\begin{equation}\label{sigmsFdef}
\left(\hat{\sigma}_F^{[s]}(n;\Theta_k)\right)^2:= \frac{1}{n}\sum_{\ell=1}^{n}\left(F(\X_k^s;\Theta_{k,\ell}^{0,s})-\Fbarsk(\srandom_k,n;\Theta_k^{0,s})\right)^2.
\end{equation}
 
Thus, \SSadd{Algorithms~\ref{algoF} and~\ref{algoR}} use random models defined by 
\begin{equation}\label{eq:randoModel}
\Mhatk(\s)=\Fbar^0_k(\bm{0},N_k;\Theta_k)+\GkhatRandom^\top\s+\frac{1}{2}\s^\top\HkhatRandom\s,
\end{equation}
where the $i$th coordinate of $\GkhatRandom\in\rp$, obtained by central finite-difference approximation, is given by 
\begin{equation}\label{eq:Ghatki}
\GhatRandom_{k,i}:=\frac{\Fbarik(\siirandom_k,N_k;\Theta_k)-\Fbar^{p+i}_k(\srandom^{p+i},N_k;\Theta_k)}{2\Dk}, \qquad i\in\intbracket{1,p},
\end{equation}
while $\HkhatRandom=\operatorname{diag}\left(\HhatRandomComp_{k,1},\dots,\HhatRandomComp_{k,p}\right)$, with $\HhatRandomComp_{k,i}:=\min\accolade{\kappah, \max\accolade{-\kappah,\HtildeRandomComp_{k,i}}}$ for some constant $\kappah\geq 1$, where 
\begin{equation*}
\HtildeRandomComp_{k,i}:=\frac{\Fbarik(\siirandom_k,N_k;\Theta_k)-2\Fbar^0_k(\bm{0},N_k;\Theta_k)+\Fbar^{p+i}_k(\srandom^{p+i}_k,N_k;\Theta_k)}{\Dk^2},
\end{equation*}
so that every realization $\Hkhat$ of $\HkhatRandom$ satisfies 
\begin{equation}\label{eq:HessianBound}
\normeshort{\Hkhat}\leq\kappah
\end{equation} \SSaddnew{in the spectral norm.}
 
Let $(\QkRandom)_{ij}$, with $(i,j)\in \intbracket{1,d}\times\intbracket{1,p}=:\mathbb{I}_{d,p}$, denote the entries of $\QkRandom$. The aforementioned conditioning on the past is formalized, inspired by~\cite{chen2018stochastic,DzaWildSub2022}, by considering the filtrations $\accoladekinN{\falgebra}$, $\accoladekinN{\fqalgebra}$ and $\accoladekinN{\fmqalgebra}$ defined~by
\begin{equation*}
\begin{array}{rcl}
\falgebra &:=& \sigma\left(\Foell, \Fsell, \Mhatell,  (\QellRandom)_{ij} \mbox{ with}\  \ell \in \intbracket{0,k-1} \mbox{ \&}\ (i,j)\in \mathbb{I}_{d,p} \right)\\
\fqalgebra &:=& \sigma\left(\Foell, \Fsell, \Mhatell, (\QellRandom)_{ij}, (\QkRandom)_{ij} \mbox{ with}\  \ell \in \intbracket{0,k-1} \mbox{ \&}\ (i,j)\in \mathbb{I}_{d,p} \right)\\
\fmqalgebra &:=& {\sigma\left(\Foell, \Fsell, \Mhatell, (\QellRandom)_{ij}, (\QkRandom)_{ij}, \ \Mhatk \mbox{ with}\  \ell \in \intbracket{0,k-1} \mbox{ \&}\ (i,j)\in \mathbb{I}_{d,p}  \right),}
\end{array} 
\end{equation*}
where $\mathscr{F}_{-1}:=\sigma(\x^0_0)$ for {\it completeness}~\cite{billingsley1995ThirdEdition}. Thus, as in~\cite{DzaWildSub2022}, it holds by construction that $\falgebra\subset\fqalgebra\subset\fmqalgebra$. 
\SSaddnew{Since the trial step is a deterministic Borel-measurable function of the model coefficients and of the radius (Section~\ref{sec:theMethods}), and $\Dk$ is $\falgebra$-measurable, $\srandom_k$ and $\X^s_k$ are $\fmqalgebra$-measurable. Throughout, 
$\Theta_k$ is independent of $\fqalgebra$ and $\Theta^{0,s}_k$ independent of $\fmqalgebra$, which is required by Lemma~\ref{Lemma:24} at each of its uses.} 
 
\SSadd{In \ASTRAR\ (Algorithm~\ref{algoR}), no new randomness is generated after the model is built: the estimate at the trial point is
\begin{equation}\label{eq:FbarsReuse}
\Fbarsk(\srandom_k,\nkrandom;\Theta_k)=\frac{1}{\nkrandom}\sum_{\ell=1}^{\nkrandom}F(\X_k^s;\Theta_{k,\ell}),
\end{equation}
computed from the same $\Theta_{k,1},\dots,\Theta_{k,\nkrandom}$ used in~\eqref{nikDefEq}, and $\Fbarokcur:=\Fbar^0_k(\bm{0},\nkrandom;\Theta_k)$. Consequently the second sampling rule~\eqref{n0sDefEq} and the estimators~\eqref{sigm0Fdef}--\eqref{sigmsFdef} are not used. 
}
 
\SSadd{For \ASTRAR, the relevant conditioning is by the within-iteration filtration
\begin{equation}\label{eq:GkFiltration}
\gqalgebra{k,n}:=\fqalgebra\vee\sigma(\Theta_{k,1},\Theta_{k,2},\dots,\Theta_{k,n}),\qquad \gqalgebra{k,0}:=\fqalgebra,
\end{equation}
\SSaddnew{also used in the proof of Theorem~\ref{Theorem:22}}} together with {\bl $\falgebra\subseteq\fqalgebra$. 
For \ASTRAR, the third $\sigma$-algebra $\fmqalgebra$ is not needed in the analysis. Although $\Mhatk$, $\srandom_k$, $\X^s_k$ and the estimate~\eqref{eq:FbarsReuse} are $\gqalgebra{k,\nkrandom}$-measurable, the reused replications are generally not independent of $\fmqalgebra$, since they were also used to construct the model and select the trial point. Consequently, the conditional centering required for the pointwise trial-estimate accuracy bound of Theorem~\ref{Theorem:22} is no longer guaranteed. This motivates the uniform accuracy result developed in Section~\ref{sec:uniformAcc}, which also covers the sample-dependent trial point. The next lemma and assumption provide the basic probabilistic ingredients needed to analyze the adaptive sampling rules. In particular, Lemma~\ref{Lemma:24} justifies conditional expectations at random points selected before an independent replication is drawn, while Assumption~\ref{Assumpt:31} provides the required integrability; together they are used in Proposition~\ref{prop:finiteNk} to establish the stopping-time property and almost-sure finiteness of the adaptive sample sizes.}

\begin{lemma}\label{Lemma:24}
Let $Z$ and~$\Theta$ be random maps into $(\mathcal{Z},\mathscr{Z})$ and $(\mathcal{T},\mathscr{T})$, respectively. Let~$\halgebra$ be a $\sigma$-algebra and~$\psi$ be a measurable real-valued function on $\left(\mathcal{Z}\times\mathcal{T},\mathscr{Z}\otimes\mathscr{T}\right)$. Assume that~$Z$ is $\halgebra$-measurable, and $\sigma(\Theta)$ and~$\halgebra$ are independent. If $\E{\abs{\psi(Z,\Theta)}}<\infty$ (\cite[Theorem~2.7-($\ell$)]{bhattacharya2007basic}), or if~$\psi$ is either nonnegative or bounded (\cite[Proposition~A.2.5]{LamLap2008secEd}), then $\E{\psi(Z,\Theta)|\halgebra}=\Esp_{\Theta}\left[\psi(Z,\Theta)\right]=\zeta(Z)$ almost surely, where $\zeta(z):=\E{\psi(z,\Theta)}$ for all $z\in\mathcal{Z}$.
\end{lemma}
 
\begin{assumption}\label{Assumpt:31}
For all $\x\in \rd$, 
$\sigma^2_F(\x):=\var{F(\x;\Theta)}\leq v_F<\infty$ for some constant~$v_F>0$.
\end{assumption}
 \SSaddnew{For fixed $k$ the criteria depend only on the replications drawn so far, so $\nkrandom$ and $\nskrandom$ are stopping times for the within-iteration filtrations. Their finiteness, proven in the next proposition, is needed for $W_\varepsilon$ of~\eqref{eq:Weps} to be finite.}
SSadd{\begin{proposition}\label{prop:finiteNk}
Let Assumption~\ref{Assumpt:31} hold and let $\nkrandom$ be given by~\eqref{nikDefEq}. Then $\nkrandom$ is a stopping time with respect to $\accolade{\gqalgebra{k,n}}_{n\geq 0}$,
\begin{equation}\label{eq:NkFloorLB}
\nkrandom\geq \Nkbar:=\left\lceil\frac{\sigma_0^2\lambda_k}{\kappaas^2\Dk^2}\right\rceil\qquad\text{almost surely},
\end{equation}
and $\prob{\bigcap_{k\in\N}\accolade{\nkrandom<\infty}}=1$. For \ASTRAF\ (Algorithm~\ref{algoF}), $\prob{\bigcap_{k\in\N}\accolade{\nskrandom<\infty}}=1$, {\bl with~$\nskrandom$ satisfying~\eqref{eq:NkFloorLB}}.
\end{proposition}}
\SSadd{\begin{proof}
The stopping-time property and the floor are immediate: for each $n$, the criterion in~\eqref{nikDefEq} depends only on $\fqalgebra$-measurable quantities and on $\Theta_{k,1},\dots,\Theta_{k,n}$, and $\sigma_0\vee\hat\sigma^{\max}_F\ge\sigma_0$. Finiteness is the finiteness of the Chow--Robbins stopping rule~\cite{ChowRobbins1965,Mukhopadhyay1978}, here applied conditionally on the design, which is random; the proof, via the conditional strong law of~\cite{MaNoZiCondSLLN2005}, is in Appendix~\ref{app:finitenessOriginal}.
\end{proof}}

\SSaddnew{\begin{remark}[reuse versus a fresh stream]\label{rem:reuse}
\ASTRAR\ is the variant one would implement. Per iteration it makes $(q+2)\nkrandom$ oracle calls against $(q+1)\nkrandom+2\nskrandom$ for \ASTRAF, it runs one sequential sampling stage instead of two, it never re-evaluates the incumbent, and the estimate of $f(\Xzerok)$ in the acceptance test is the one the model was built from, so the common random numbers are kept across the whole iteration. In the analysis, the single accuracy event of Corollary~\ref{cor:UniformEvent} replaces the pair of model and estimate events of \ASTRAF: the third $\sigma$-algebra $\fmqalgebra$ is not needed, the lower bound on the conditional probability of a good iteration is the accuracy level itself rather than its square, and the probability-dependent constants of Section~\ref{sec:complexity} and the admissible sketch failure probability $\beta_Q$ improve accordingly (Remark~\ref{rem:theorem32R}). The price is a selection bias in the trial estimate, which is evaluated on the replications that selected it: no pointwise accuracy statement is available there, and Section~\ref{sec:uniformAcc} replaces it by one holding uniformly over the trust region, under one further assumption, continuity of the sample paths (Assumption~\ref{Assumption:cont}). \ASTRAF\ needs no such assumption, assesses the predicted decrease on replications independent of the model, and admits the pointwise analysis of~\cite{chen2018stochastic,DzaWildSub2022}. Both variants tolerate occasionally inaccurate models and estimates and reach the same rates.
\end{remark}}

\section{Probabilistic well alignment, weak accuracy, and {\bl subspace accuracy}}\label{sec:Prob}

{\bl \bl The convergence and complexity analyses of Sections~\ref{sec:convAnalysis} and~\ref{sec:complexity} \SSaddnew{use three properties of an iteration. First, the sketch must be well aligned (Assumption~\ref{wellAssump})}; this concerns $\QkRandom$ alone and is conditioned on $\falgebra$. Second, 
the model must be \emph{subspace-accurate} (Definition~\ref{def:semilinear}): its gradient approximates that of the reduced objective $\s\mapsto f(\xzerok+\Qk\s)$ to within $\frac12\kmge\Dk$ on the trust region. This corresponds to the gradient component of ``fully linear'' models and controls 
\SSaddnew{the error in the predicted decrease (Lemma~\ref{modFuncErrLemma}), and no accuracy of the model values is needed because the acceptance test uses only differences. Third,} the trial estimates must be \emph{weakly accurate} (Definition~\ref{def:weakAcc}): the difference of their estimation errors is bounded by $\kfde\Dk^2$. \SSaddnew{The second and third properties jointly will deem a ``good'' iteration and }this section establishes \SSaddnew{a lower bound on that probability} 
for Algorithms~\ref{algoF} and~\ref{algoR}.
	
The model-construction analysis is common to both algorithms. 
Under CRN, a Bernstein moment condition on \emph{paired  increments} \SSaddnew{of errors} (Assumption~\ref{Assumption:33}) 
\SSaddnew{makes} sample sizes of order $\lambda_k\Dk^{-2}$ \SSaddnew{sufficient}, rather than the order $\Dk^{-4}$ scaling typically associated with independent sampling. 
Since $\nkrandom$ is a stopping time, \SSaddnew{the concentration bound must hold at a stopping time (data-dependent) sample size; }Lemma~\ref{Lemma:52} provides \SSaddnew{that through }an exponential supermartingale and Ville's inequality. Applied at the $2p$ design points, \SSaddnew{which were selected before the replications are drawn, it gives weak accuracy of their differences with the centering design point} 
(Theorem~\ref{Theorem:22}), from which model subspace-accuracy follows (Theorem~\ref{mgeTheorem}) \SSaddnew{ for both algorithms}.
	
The algorithms differ in trial-point evaluation. In \ASTRAF, \SSaddnew{ the trial point is evaluated on $\Theta^{0,s}_k$, }
a fresh replication stream independent of model construction and trial-point selection, \SSaddnew{which} allows 
the pointwise accuracy bound to hold conditionally on $\fmqalgebra$ (Theorem~\ref{Theorem:22} and Corollary~\ref{theor31bis} \SSaddnew{in Section~\ref{sec:pointwiseAcc}}). In \ASTRAR, the model-construction samples are reused \SSaddnew{and the trial point is evaluated on the replications that selected it, so its estimate is not conditionally unbiased given the model, which is why existing STORM-type analyses do not cover this scheme as} the corresponding pointwise conditional argument is generally unavailable.  
Instead, Theorem~\ref{Theorem:Uniform} in Section~\ref{sec:uniformAcc} \SSaddnew{ replaces the pointwise bound by one that holds uniformly over the trust region, hence at the trial point however it was chosen.} 
}
\SSaddnew{Corollary~\ref{theor31bis} for \ASTRAF\ and Corollary~\ref{cor:UniformEvent} for \ASTRAR\ play the same role downstream: each bounds from below the conditional probability of a good iteration. What reuse costs and what it buys is recorded in Remark~\ref{rem:reuse}.}

In order to control differences of errors (computed with the same random seed) used in our analysis, we make the following Bernstein regularity assumption on the paired CRN increment {\bl $D(\x,\y;\Theta)$ 
given by}
\begin{equation}\label{eq:DIncrement}
D(\x,\y;\Theta):=(F(\y;\Theta)-f(\y))-(F(\x;\Theta)-f(\x)),
\end{equation}
satisfying $\E{D(\x,\y;\Theta)}=0$ and which is integrable thanks to Assumption~\ref{Assumpt:31}. 
{\bl The following assumptions provide the geometric and probabilistic regularity needed later for the accuracy analysis. Assumption~\ref{Assumption:Qmax} uniformly controls the action of the random subspace matrices, while Assumption~\ref{Assumption:33} imposes a Bernstein-type moment condition directly on the paired CRN error increments, rather than on the individual stochastic errors as in \cite[Assumption~6]{ha2025complexity}. This more targeted condition captures the cancellation induced by CRN and scales with the distance between the paired points; together with Assumption~\ref{Assumption:Qmax}, it allows the concentration bounds below to be expressed in terms of reduced-space displacements.} 

\begin{assumption}[Inspired by~\cite{CRsubspace2021,DzaWildSub2022,DzaWildSubDir24}]\label{Assumption:Qmax}
There exists a constant $\Qmax>0$ such that any matrix realization $\Qk$ used by \SSadd{Algorithm~\ref{algoF} or~\ref{algoR}} satisfies $\norme{\Qk}\leq\Qmax$ for all $k\in\N$. 
\end{assumption}
\begin{assumption}\label{Assumption:33}
There exist constants $b_D, \sigma_D>0$ such that for every $\x,\y\in\rd$ with $\normii{\y-\x}\leq \Qmax\deltmax$, and every integer $\ssmath\geq 2$
\begin{equation*}
\E{\abs{D(\x,\y;\Theta)}^\ssmath}\leq \frac{\ssmath!}{2}(b_D\normii{\y-\x})^{\ssmath-2}\sigma_D^2\normii{\y-\x}^2.
\end{equation*}
\end{assumption} 
{\bl The following result gives readily verifiable sufficient conditions for Assumption~\ref{Assumption:33}, thereby showing that its Bernstein-type control of paired error increments is satisfied under a natural stochastic Lipschitz regularity condition.}
\begin{proposition}\label{prop:BernsteinSuff}
Suppose that for all $\x,\y\in\rd$,
\begin{equation}\label{kappaFLipschitz1}
\abs{F(\y;\Theta)-F(\x;\Theta)}\leq \kappaFL(\Theta)\normii{\y-\x}\qquad\text{almost surely,}
\end{equation}
where the integrable random variable $\kappaFL(\Theta)\geq 0$ satisfies\footnote{We note that this assumption is 
\SSaddnew{satisfied by many light-tailed distributions} 
since the Bernstein-type moment condition $\E{\abs{Z}^\ssmath}\leq \frac{\ssmath!}{2}b^{\ssmath-2}\sigma^2$ is satisfied by a large class of variables~$Z$, including bounded, sub-Gaussian, and more generally sub-exponential and sub-gamma random variables.} 
\begin{equation}\label{Bernstein:cond2}
\E{\left(\kappaFL(\Theta)+\E{\kappaFL(\Theta)}\right)^\ssmath}\leq \frac{\ssmath!}{2}b_{\rm FL}^{\ssmath-2}\sigma_{\rm FL}^2,
\end{equation}
for constants $b_{\rm FL}, \sigma_{\rm FL}>0$ and \SSaddnew{all} integer\SSaddnew{s} $\ssmath\geq 2$. Then Assumption~\ref{Assumption:33} holds with $b_D=b_{\rm FL}$ and $\sigma_D=\sigma_{\rm FL}$.
\end{proposition}
 
\begin{proof}
\SSadd{By~\eqref{kappaFLipschitz1} and Jensen's inequality, $\abs{D(\x,\y;\Theta)}\leq\parenth{\kappaFL(\Theta)+\E{\kappaFL(\Theta)}}\normii{\y-\x}$, as in~\cite[Proposition~2.6]{rinaldi2024stochastic}, where the Lipschitz constant is deterministic. Raising both sides to the power $\ssmath$, taking expectations and applying~\eqref{Bernstein:cond2} gives the claim.}
\end{proof}
{\bl The next two lemmas provide the time-uniform concentration tools needed to handle the adaptive, data-dependent sample sizes used below. Lemma~\ref{VilleLemma} supplies the maximal inequality underlying the argument, while Lemma~\ref{Lemma:52} gives the conditional Bernstein bound applied directly in Section~\ref{sec:pointwiseAcc} for pointwise weak accuracy and, in Section~\ref{sec:uniformAcc}, to the individual edges of the chaining construction.}
\begin{lemma}\label{VilleLemma}(Ville's maximal inequality  for nonnegative supermartingales~\cite[Lemma~1]{HowRamMcASek2020}).
If $\accolade{L_t}_{t\in \mathscr{T}\cup\accolade{0}}$ is a nonnegative supermartingale with respect to the filtration $\accolade{\mathscr{G}_t}_{t\in \mathscr{T}\cup\accolade{0}}$, then for any $a>0$, 
\begin{equation*}
\prob{\exists t\in\mathscr{T}:L_t\geq a\, |\,\mathscr{G}_0}\leq \frac{L_0}{a}\qquad\text{almost surely.}
\end{equation*}
\end{lemma}

The following conditional time-uniform Bernstein inequality is a specialization of the scalar sub-gamma line-crossing bounds of Howard et al.~\cite{HowRamMcASek2020}; see in particular~\cite[Corollary~4(b)]{HowRamMcASek2020}. We state and prove the result in the form needed here because the Bernstein parameters and the lower sample-size index are allowed to be $\mathscr{G}_0$-measurable random variables. The proof makes this conditional extension explicit by constructing the corresponding exponential supermartingale and applying Lemma~\ref{VilleLemma} after a $\mathscr{G}_0$-measurable normalization.

\begin{lemma}\label{Lemma:52}
Let 
$\accolade{Z_j}_{j\geq 1}$ \SSaddnew{be adapted to filtration $\accolade{\gjalgebra}_{j\in\N}$, that is, $Z_j$ is $\gjalgebra$-measurable, and} satisfy $\E{Z_j|\galgebra}=0$ almost surely. Assume that for every integer $\ssmath\geq 2$, 
\begin{equation}\label{MomentCond20}
\E{\abs{Z_j}^\ssmath|\galgebra}\leq \frac{\ssmath !}{2}B^{\ssmath-2}V^2\qquad\text{almost surely,}
\end{equation}
where $B,V>0$ are $\mathscr{G}_0$-measurable random variables. Let~$A>0$ be a $\mathscr{G}_0$-measurable random variable, and~$\aleph$ be a positive, integer-valued $\mathscr{G}_0$-measurable random variable. Then 
\begin{equation}\label{eq:Bernstein1}
\prob{\left.\underset{n\geq \aleph}{\sup}\, \frac{1}{n}\abs{\sum_{j=1}^{n}Z_j}\SSaddnew{>}
A\right\rvert\mathscr{G}_0}\leq 2\exp\left(-\frac{\aleph A^2}{2(V^2+BA)}\right)\qquad\text{almost surely.}
\end{equation}
\end{lemma}

\begin{proof}
Set $S_n:=\sum_{j=1}^{n}Z_j$ and $S_0=0$. We first show that 
\begin{equation}\label{eq:36}
\prob{\left.\exists n\geq \aleph:\frac{S_n}{n}\geq A\right\rvert\mathscr{G}_0}\leq \exp\left(-\frac{\aleph A^2}{2(V^2+BA)}\right)\qquad\text{almost surely.}
\end{equation}
For every~$z\in\R$ and the $\mathscr{G}_0$-measurable random variable $\Upsilon:=\frac{A}{V^2+BA}$ satisfying $0<\Upsilon<B^{-1}$ almost surely, 
\begin{equation*}
e^{\Upsilon z}=1+\Upsilon z+\sum_{\ssmath=2}^{\infty}\frac{(\Upsilon z)^\ssmath}{\ssmath !}\leq 1+\Upsilon z+\sum_{\ssmath=2}^{\infty}\frac{\Upsilon^\ssmath \abs{z}^\ssmath}{\ssmath !}.
\end{equation*}
Because $\Upsilon$ is $\galgebra$-measurable, taking conditional expectations, using $\E{Z_j|\galgebra}=0$ almost surely, the moment condition~\eqref{MomentCond20}, and conditional monotone convergence for nonnegative \SSaddnew{random variables, with conditional expectations allowed to take the value $+\infty$ (see, e.g.,~\cite[Theorem~5.1]{kallenberg1997foundations}) }
gives almost surely,
\begin{equation}\label{eq:condMon1}
\begin{split}
\E{e^{\Upsilon Z_j}|\galgebra}&\leq 1+\sum_{\ssmath=2}^{\infty}\frac{\Upsilon^\ssmath}{\ssmath !}\E{\abs{Z_j}^\ssmath|\galgebra}\leq 1+\frac{\Upsilon^2V^2}{2}\sum_{\ssmath=2}^{\infty}(\Upsilon B)^{\ssmath-2}\\
&=1+\frac{\Upsilon^2 V^2}{2(1-\Upsilon B)}\leq \exp\left(\frac{\Upsilon^2 V^2}{2(1-\Upsilon B)}\right).
\end{split}
\end{equation}
\SSaddnew{All conditional expectations here are of variables that are conditionally integrable given $\galgebra$ ($Z_j$ by~\eqref{MomentCond20} with $\ssmath=2$), taken in the generalized sense of~\cite[Theorem~5.1]{kallenberg1997foundations}, under which linearity and the pull-out of the factor $\Upsilon$ hold.}
Define $\psi(\Upsilon):=\frac{\Upsilon^2 V^2}{2(1-\Upsilon B)}=\frac{A^2}{2(V^2+BA)}$ and $L_n:=\exp(\Upsilon S_n-n\psi(\Upsilon)), n\geq 0$, so that $L_0=1$ and $L_n$ is nonnegative and $\gnalgebra$-measurable. We note that $\frac{L_n}{L_{n-1}}=e^{\Upsilon Z_n}e^{-\psi(\Upsilon)}$, and since~$\psi(\Upsilon)$ is $\mathscr{G}_{n-1}$-measurable, then it follows from~\eqref{eq:condMon1} that 
\begin{equation*}
\E{\left.\frac{L_n}{L_{n-1}}\right\rvert\mathscr{G}_{n-1}}=\E{\left.e^{\Upsilon Z_n}e^{-\psi(\Upsilon)}\right\rvert\mathscr{G}_{n-1}}=e^{-\psi(\Upsilon)}\E{\left.e^{\Upsilon Z_n}\right\rvert\mathscr{G}_{n-1}}\leq 1,
\end{equation*}
which shows that $\E{L_n|\mathscr{G}_{n-1}}\leq L_{n-1}$. In particular, by induction, $\E{L_n}\leq\E{L_0}=1$, so that each $L_n$ is integrable. Hence $\accolade{L_n}_{n\geq 0}$ is a nonnegative supermartingale. Now define $D:=\Upsilon A-\psi(\Upsilon)=\frac{A^2}{2(V^2+BA)}$, and suppose that for some $n\geq \aleph$, $S_n\geq n A$. Then,
\begin{equation*}
L_n\geq \exp(n(\Upsilon A-\psi(\Upsilon)))=e^{nD}\geq e^{\aleph D},
\end{equation*}
hence
\begin{equation}\label{eq:reddot}
\accolade{\exists n\geq \aleph:S_n\geq nA}\subseteq\accolade{\exists n\geq 1:L_n\geq e^{\aleph D}}.
\end{equation}
Moreover, $e^{-\aleph D}\in [0,1]$ is $\mathscr{G}_0$-measurable, $\tilde{L}_n:=e^{-\aleph D}L_n$ is trivially integrable, satisfying 
\begin{equation*}
\E{\tilde{L}_n|\mathscr{G}_{n-1}}=e^{-\aleph D}\E{L_n|\mathscr{G}_{n-1}}\leq e^{-\aleph D}L_{n-1}=\tilde{L}_{n-1},
\end{equation*}
hence the normalized process $\accolade{\tilde{L}_n}_{n\geq 0}$ is a nonnegative supermartingale with $\tilde{L}_0=e^{-\aleph D}$. We note that $L_n\geq e^{\aleph D}$ if and only if $\tilde{L}_n\geq 1$, hence it follows from~\eqref{eq:reddot} that
\begin{equation}\label{eq:red2dot}
\accolade{\exists n\geq \aleph:S_n\geq nA}\subseteq\accolade{\exists n\geq 1:\tilde{L}_n\geq 1}.
\end{equation}
According to Lemma~\ref{VilleLemma}, almost surely, 
\begin{equation*}
\prob{\left. \exists n\geq 1:\tilde{L}_n\geq 1\right\rvert\mathscr{G}_0}\leq \tilde{L}_0=e^{-\aleph D}, 
\end{equation*}
which, combined with~\eqref{eq:red2dot}, yields~\eqref{eq:36}. 
 
The negative tail is handled by applying exactly the same arguments to $\accolade{-Z_j}_{j\geq 1}$. Indeed, $\E{-Z_j|\galgebra}=0$ and $\E{\abs{-Z_j}^\ssmath|\galgebra}=\E{\abs{Z_j}^\ssmath|\galgebra}$, so 
\begin{equation}\label{eq:37}
\prob{\left.\exists n\geq \aleph:\frac{S_n}{n}\leq -A\right\rvert\mathscr{G}_0}\leq \exp\left(-\frac{\aleph A^2}{2(V^2+BA)}\right)\qquad\text{almost surely.}
\end{equation}
We finally note that 
\begin{equation*}
\accolade{\underset{n\geq \aleph}{\sup}\, \frac{\abs{S_n}}{n}\SSaddnew{>}
 A}=\accolade{\exists n\geq \aleph: S_n\SSaddnew{>}
nA}\cup\accolade{\exists n\geq \aleph: S_n<
-nA}\SSaddnew{\subseteq\accolade{\exists n\geq \aleph: S_n\geq 
nA}\cup\accolade{\exists n\geq \aleph: S_n\leq 
-nA}},
\end{equation*}
so the conditional union bound together with~\eqref{eq:36} and~\eqref{eq:37} yields~\eqref{eq:Bernstein1}, which completes the proof.
\end{proof}

\SSadd{\subsection{Pointwise weak accuracy}\label{sec:pointwiseAcc}
{\bl Following the notions of $\epsilon_f$-accurate estimates~\cite{chen2018stochastic,DzaWildSub2022}, we introduce a weaker requirement (motivated by the paired-error balance condition of~\cite{ha2025complexity}) that controls only the difference of the two estimation errors at the scale~$\dk^2$. Theorem~\ref{Theorem:22} and Corollary~\ref{theor31bis} show that it holds with high conditional probability under Assumptions~\ref{Assumpt:31},~\ref{Assumption:Qmax}, and~\ref{Assumption:33}.}
 }

\SSadd{\begin{definition}\label{def:weakAcc}
Given $\kfde>0$ and $\Qk\sk\in\rd$, $\fbarokcur:=\fbarokcur\left(\bm{0},\nokbrak\right)$ and $\fbarsk:=\fbarsk\left(\sk,\nskbrak\right)$ are called $\kfde$-weakly accurate estimates of $f(\xzerok)$ and $f(\xzerok+\Qk\sk)$, respectively, for a given $\dk$ if 
\begin{equation}\label{weakDef1}
\abs{\ebarsk-\ebarokcur}=\abs{\parenth{\fbarsk-f(\xzerok+\Qk\sk)}-\parenth{\fbarokcur-f(\xzerok)}}\leq \kfde\dk^2,
\end{equation}
where $\ebarsk:=\fbarsk-f(\xzerok+\Qk\sk)$ and $\ebarokcur:=\fbarokcur-f(\xzerok)$.
\end{definition}

\begin{theorem}\label{Theorem:22}
Let Assumptions~\ref{Assumpt:31}, \ref{Assumption:Qmax} and~\ref{Assumption:33} hold, and consider the random errors $\erandombar^i_k(\nkrandom)$, $i\in\intbracket{0,q}$, defined by $\Fbarik(\srandom^i_k,\nkrandom;\Theta_k)=f(\Xki)+\erandombar^i_k(\nkrandom)$, \SSaddnew{with} $\Theta_k$ \SSaddnew{as in Section~\ref{sec:estiModels}.} 
For all~$\kfde>0$ and all $k\geq 0$,
\begin{equation}\label{eq:chifde}
\prob{\left.\abs{\erandombar^i_k(\nkrandom)-\erandombar^0_k(\nkrandom)}\leq \kfde\Dk^2\right\rvert\fqalgebra}\geq 1-2e^{-\chifde\lambda_k} \qquad\text{almost surely,}
\end{equation}
where $\chifde:=\frac{\sigma_0^2\kfde^2}{2\kappaas^2(\sigma_D^2\Qmax^2+b_D\Qmax\kfde\deltmax)}$. 
 
\SSaddnew{Likewise} for \ASTRAF, for all~$\kfde>0$ and all $k\geq 0$, the events $\jk:=\accolade{\abs{\erandombar^{[s]}_k\parenth{\nskrandom}-\erandombar^{[0]}_k\parenth{\nokrandombrak}}\leq \kfde\Dk^2}$ satisfy
\begin{equation}\label{jkEq1} 
\prob{\left.\jk\right\rvert\fmqalgebra}\geq 1-2e^{-\chifde\lambda_k} \qquad\text{almost surely},
\end{equation}
where the errors $\erandombar^{[s]}_k\parenth{\nskrandom}$ and $\erandombar^{[0]}_k\parenth{\nokrandombrak}$ are defined by $\erandombar^{[s]}_k\parenth{\nskrandom}=\Fbarsk\parenth{\srandom_k,\nskrandom;\Theta^{0,s}_k}-f(\xrandomvec^s_k)$, and $\erandombar^{[0]}_k\parenth{\nokrandombrak}=\Fbarokcur\parenth{\bm{0},\nokrandombrak;\Theta^{0,s}_k}-f(\Xzerok)$, respectively, with $\Theta^{0,s}_k$ being an independent copy of $\Theta_k$.
\end{theorem}}

\begin{proof}[Proof of Theorem~\ref{Theorem:22}]
For each $\ell\in \intbracket{1,n}$, let $\Ekli:=F(\Xki;\Theta_{k,\ell})-f(\Xki)$, $\Dkli:=\Ekli-\Eklo=D(\SSaddnew{\X_k^0,}\Xki
;\Theta_{k,\ell})$ and $\Dkls:=D(\X_k^\SSaddnew{0},\X_k^\SSaddnew{s};\Theta_{k,\ell}^{0,s})$, with~$D$ introduced in~\eqref{eq:DIncrement}, so that $\erandombar^i_k(n)-\erandombar^0_k(n)=\frac{1}{n}\sum_{\ell=1}^{n}\Dkli$. 
\SSaddnew{Recall $\gqalgebra{k,n}$ from~\eqref{eq:GkFiltration}.}
For $\ssmath\geq 2$, consider the deterministic moment function $\mu_\ssmath(\x,\yb):=\Esp_{\Theta}\left[\abs{D(\x,\y;\Theta)}^\ssmath\right]$. We note that the function $(\x,\yb,t)\mapsto\abs{D(\x,\yb,t)}^\ssmath$ is trivially nonnegative, $\Xki$ and $\Xzerok$ are $\gqalgebra{k,n-1}$-measurable, and $\sigma(\Theta_{k,n})$ and $\gqalgebra{k,n-1}$ are independent. Thus, it follows from Lemma~\ref{Lemma:24} that almost surely,
\begin{equation*}
\begin{split}
\E{\left.\abs{\Dkni}^\ssmath\right\rvert\gqalgebra{k,n-1}}&=\E{\left.\abs{D(\SSaddnew{\X_k^0,}\Xki
;\Theta_{k,n})}^\ssmath\right\rvert\gqalgebra{k,n-1}}=\mu_\ssmath(\SSaddnew{\X_k^0,}\Xki
)=\Esp_{\Theta}\left[\abs{D(\SSaddnew{\X_k^0,}\Xki
;\Theta_{k,n})}^\ssmath\right].
\end{split}
\end{equation*}
According to Assumption~\ref{Assumption:33}, there exist $b_D,\sigma_D>0$ such that 
\begin{equation*}
\Esp_{\Theta}\left[\abs{D(\SSaddnew{\X_k^0,}\Xki
;\Theta_{k,n})}^\ssmath\right]\leq \frac{\ssmath!}{2}(b_D\Qmax\Dk)^{\ssmath-2}(\sigma_D\Qmax\Dk)^2\qquad\text{almost surely},
\end{equation*}
and since $\E{\left.\Dkni\right\rvert\gqalgebra{k,n-1}}=\E{\left.D(\SSaddnew{\X_k^0,}\Xki
;\Theta_{k,n})\right\rvert\gqalgebra{k,n-1}}=0$, it follows from Lemma~\ref{Lemma:52} applied with $B=b_D\Qmax\Dk$, $V=\sigma_D\Qmax\Dk$, $\aleph=\Nkbar=\left\lceil\frac{\sigma_0^2\lambda_k}{\kappaas^2\Dk^2}\right\rceil$, $A=\kfde\Dk^2$ and $\mathscr{G}_0=\gqalgebra{k,0}$, that almost surely,
\begin{equation}\label{eq:supNkbar}
\begin{split}
\prob{\left.\underset{n\geq \Nkbar}{\sup}\, \frac{1}{n}\abs{\sum_{\ell=1}^{n}\Dkli}\SSaddnew{>}
\kfde\Dk^2\right\rvert\gqalgebra{k,0}}&\leq 2\exp\left(-\frac{\Nkbar\kfde^2\Dk^4}{2(\sigma_D^2\Qmax^2\Dk^2+b_D\Qmax\kfde\Dk^3)}\right)\leq 2\exp\left(-\chifde\lambda_k\right),
\end{split}
\end{equation}
where we note that all the random quantities $B,V,A$ and~$\aleph$ are $\gqalgebra{k,0}$-measurable given that $\Dk$ is $\fqalgebra$-measurable. Thus,~\eqref{eq:chifde} 
follows from~\eqref{eq:supNkbar} and the fact that
\begin{equation*}
\accolade{\abs{\erandombar^i_k(\nkrandom)-\erandombar^0_k(\nkrandom)}> \kfde\Dk^2}\subseteq\accolade{\underset{n\geq \Nkbar}{\sup}\, \frac{1}{n}\abs{\sum_{\ell=1}^{n}\Dkli}\SSaddnew{>}
\kfde\Dk^2}. 
\end{equation*}
The second result follows by repeating the same arguments, making use of the $\sigma$-algebras $\gmqalgebra{k,n}$ given by $\gmqalgebra{k,n}:=\fmqalgebra\vee\sigma(
\SSaddnew{\Theta^{0,s}_{k,1},\dots,\Theta^{0,s}_{k,n}})$ with $\gmqalgebra{k,0}:=\fmqalgebra$, and observing that $\Xzerok,\X_k^s$ are $\gmqalgebra{k,n-1}$-measurable, $\sigma(\Theta_{k,n}^{0,s})$ and $\gmqalgebra{k,n-1}$ are independent, $\E{\Dkns|\gmqalgebra{k,n-1}}=0$ almost surely, $N_k^{0,s}\geq \Nkbar$ and $\normeshort{\X_k^s-\Xzerok}_2\leq \Qmax\Dk$, and the proof is completed.
\end{proof}

{\bl In the spirit of the probabilistically $\epsilon_f$-accurate estimates of~\cite{chen2018stochastic,DzaWildSub2022}, the following corollary recasts the bounds of Theorem~\ref{Theorem:22} at the probability level used later in the analysis. We call the estimates good on~$\jk$ for \ASTRAF\ and, as introduced in Section~\ref{sec:uniformAcc}, on~$\uk$ for \ASTRAR, and bad otherwise. The model will similarly be called good on~$\ik$ once that event is introduced below.}

\SSadd{\begin{corollary}\label{theor31bis}
Let all the assumptions made in Theorem~\ref{Theorem:22} hold. Then for any $\ef>0$ fixed, there exists a deterministic integer $\kappaf\in\N$ sufficiently large, such that for all $k\geq \kappaf$,
\begin{equation}
\prob{\left.\abs{\erandombar^i_k(\nkrandom)-\erandombar^0_k(\nkrandom)}\leq \kfde\Dk^2\right\rvert\fqalgebra}\geq \beta_{k,\ef}:= 1-\frac{1}{k^{1+\ef}}\qquad\text{almost surely}.
\end{equation}
\SSaddnew{Likewise }
for \ASTRAF, for all~$k\geq \kappaf$, the events
\SSaddnew{$\jk$ of Theorem~\ref{Theorem:22}} 
satisfy,
$\prob{\left.\jk\right\rvert\fmqalgebra}\geq \beta_{k,\ef}$ almost surely.
\end{corollary}}
{\bl The proof is given in Appendix~\ref{app:pointwise}.}

{\bl The design-point statements in Theorem~\ref{Theorem:22} and Corollary~\ref{theor31bis} apply to both \ASTRAF\ and \ASTRAR\ and are used in Theorem~\ref{mgeTheorem} to establish model subspace-accuracy. Their trial-point statements, conditioned on $\fmqalgebra$, apply specifically to \ASTRAF\ and provide its probabilistic weak-accuracy guarantee. For \ASTRAR, the corresponding guarantee is obtained from the uniform event of Corollary~\ref{cor:UniformEvent} in Section~\ref{sec:uniformAcc}.}
 
\SSadd{\subsection{Uniform weak accuracy}\label{sec:uniformAcc}}
 
\SSadd{This subsection establishes, for \ASTRAR\ (Algorithm~\ref{algoR}), the analogue of Theorem~\ref{Theorem:22} and Corollary~\ref{theor31bis}. Everything downstream of Corollary~\ref{cor:UniformEvent} below is unchanged: Sections~\ref{sec:convAnalysis} and~\ref{sec:complexity} use the accuracy results only through the conditional lower bound on the good event, and that bound is supplied here in the same form.} \SSadd{The difficulty is that the trial point is now a function of the samples being averaged, so the summands of $\Fbarsk-\Fbarok$ are not a martingale difference sequence and Lemma~\ref{Lemma:52} cannot be invoked at $\srandom_k$. What is established instead is a bound valid simultaneously at every point of the trust region, which is indifferent to the mechanism selecting one of them. Two features of the present setting make this attainable, and both are already in place: the identity
\begin{equation}\label{eq:Dtelescope}
D(\x,\y;\Theta)-D(\x,\y';\Theta)=D(\y',\y;\Theta),
\end{equation}
immediate from~\eqref{eq:DIncrement}, which shows that increments of the field $\s\mapsto\Dbar_k(\s,n)$ are themselves paired increments and are therefore governed by Assumption~\ref{Assumption:33}; and the fact that the Bernstein parameters in Lemma~\ref{Lemma:52} are $\mathscr{G}_0$-measurable, so that the tilt is fixed before the first replication of the iteration is drawn.} \SSadd{Here and below, for $\normii{\s}\leq\Dk$ we write
\begin{equation}\label{eq:Dbarfield}
\Dbar_k(\s,n):=\frac{1}{n}\sum_{\ell=1}^{n}D\parenth{\Xzerok,\Xzerok+\QkRandom\s;\Theta_{k,\ell}},
\end{equation}
so that $\Dbar_k(\siirandom_k,\nkrandom)=\erandombar^i_k(\nkrandom)-\erandombar^0_k(\nkrandom)$ at a design point and $\Dbar_k(\srandom_k,\nkrandom)=\erandombar^{s}_k(\nkrandom)-\erandombar^0_k(\nkrandom)$ at the trial point.}

\SSaddnew{Lemma~\ref{Lemma:Chained} extends Lemma~\ref{Lemma:52} from one centred sequence to a random field indexed by the unit ball. {\bl Given~$\epsilon>0$ and two subsets $\mathscr{N}$ and $\mathscr{K}$ of a metric space, with $\mathscr{N}\subseteq\mathscr{K}$, the set $\mathscr{N}$ is called an $\epsilon$-net of $\mathscr{K}$ if every point of~$\mathscr{K}$ is at distance less than~$\epsilon$ from some point of $\mathscr{N}$ (see, e.g.,~\cite[Definition~4.2.1]{vershynin2018HDP} and~\cite[Definition~3.9]{DzaWildHashing2022}). In particular, for any $\epsilon\in(0,1)$ the Euclidean unit ball of $\rp$ admits a finite $\epsilon$-net of cardinality at most $(1+2/\epsilon)^p$ by the classical \emph{volumetric estimate}~\cite[Corollary~4.2.13]{vershynin2018HDP}. Chaining, as in Dudley's \emph{multiscale net argument} (see e.g.,~\cite[Section~8.1]{vershynin2018HDP}), controls a random field over a continuum by approximating each point through successive nets of mesh~$2^{-m}$, $m\geq 1$, and expressing its value relative to the origin as a telescoping sum of increments. 
Here,  \eqref{eq:Dtelescope} identifies these increments as paired error differences, so Assumption~\ref{Assumption:33} provides Bernstein bounds with parameters proportional to the distances between successive net points. Applying Lemma~\ref{Lemma:52} to the admissible edges at each level, with geometrically decreasing failure probabilities, yields simultaneous control over all levels and all sample sizes beyond the prescribed lower bound. Continuity then extends this control to the entire ball, giving the conditional time-uniform conclusion of Lemma~\ref{Lemma:Chained}.}
}
{\bl
\begin{lemma}\label{Lemma:Chained}(Conditional time-uniform chained Bernstein inequality). Let $\Bp:=\accolade{\tb\in\rp:\normii{\tb}\leq 1}$, let $\accolade{\gjalgebra}_{j\in\N}$ be a filtration, and let $\Delta>0$ be finite and $\Gscr_0$-measurable. Consider the adapted random field $\accolade{Z_j(\tb)}_{j\geq 1,\tb\in\Bp}$ satisfying, almost surely
\begin{enumerate}
\item[\textup{(i)}] $Z_j(\bm{0})=0$ and $\tb\mapsto Z_j(\tb)$ is continuous on $\Bp$;
\item[\textup{(ii)}] $\E{Z_j(\tb)|\mathscr{G}_{j-1}}=0$ for every $\tb\in\Bp$;
\item[\textup{(iii)}] for all $\tb,\tb'\in\Bp$ with $\normii{\tb-\tb'}\leq 1$ and every integer $\ssmath\geq 2$,
\[
\E{\left.\abs{Z_j(\tb)-Z_j(\tb')}^\ssmath\right\rvert\mathscr{G}_{j-1}}\leq \frac{\ssmath!}{2}\parenth{B\Delta\normii{\tb-\tb'}}^{\ssmath-2}\parenth{V\Delta\normii{\tb-\tb'}}^2,
\]
where $B,V>0$ are $\mathscr{G}_0$-measurable.
\end{enumerate}
If $\Vscr:=\accolade{\s\in\rp:\normii{\s}\leq \Delta}$ and $\widetilde{Z}_j(\s):=Z_j(\s/\Delta)$ for all $j\geq 1$ and all $\s\in\Vscr$, then for every $\vartheta>0$ and every integer-valued $\mathscr{G}_0$-measurable random variable $\aleph>0$, almost surely,
\begin{equation}\label{eq:chained2}
\prob{\left.\underset{n\geq \aleph}{\sup}\ \underset{\s\in\mathscr{V}}{\sup}\ \frac{1}{n}\abs{\sum_{j=1}^{n}\widetilde{Z}_j(\s)}\geq 15V\Delta\sqrt{\frac{p+\vartheta}{\aleph}}+72B\Delta\frac{p+\vartheta}{\aleph}\right\rvert\mathscr{G}_0}\leq e^{-\vartheta}.
\end{equation}
\end{lemma}

\begin{proof}
The result is proved conditionally on $\Gscr_0$ with respect to which $\Delta,B,V$ and~$\aleph$ are measurable and therefore fixed. 

\emph{Step 1 (Construction of the nets and admissible edges).} For every integer $m\geq 1$, choose a deterministic $2^{-m}$-net $\Nscr_m$ of $\Bp$ satisfying
\begin{equation}\label{eq:soixante}
\abs{\Nscr_m}\leq (1+2^{m+1})^p\leq (3\cdot 2^m)^p
\end{equation}
thanks to the standard volumetric estimate, and set $\Nscr_0:=\accolade{\bm{0}}$.
Let $\pi_m:\Bp\longrightarrow\Nscr_m$ denote the nearest-point selection map, so that $\normii{\tb-\pi_m(\tb)}\leq 2^{-m}$ and hence for $m\geq 2$, 
\begin{equation*}
\normii{\pi_m(\tb)-\pi_{m-1}(\tb)}\leq 2^{-m}+2^{-(m-1)}=3\cdot 2^{-m},
\end{equation*}
while for $m=1$,
\begin{equation*}
\normii{\pi_1(\tb)-\pi_0(\tb)}\leq 1\qquad\text{since } \pi_0(\tb)=\bm{0}\quad \text{and}\quad \pi_1(\tb)\in\Bp.
\end{equation*}
Consider the finite set of admissible edges defined for all $m\geq 1$ by 
\begin{equation*}
\Escr_m:=\accolade{(\ub,\vb)\in\Nscr_{m-1}\times \Nscr_m:\normii{\vb-\ub}\leq \min\accolade{3\cdot 2^{-m},1}}
\end{equation*}
which contains every consecutive pair $(\pi_{m-1}(\tb),\pi_m(\tb))$. We observe that $\abs{\Escr_m}\leq \abs{\Nscr_{m-1}}\abs{\Nscr_m}$ which, combined with the volumetric bounds from~\eqref{eq:soixante}, yields
\begin{equation}\label{eq:Un}
\log(4\abs{\Escr_m})+m\log 2\leq 4p(m+1)=\vartheta_m-\vartheta\qquad\Longrightarrow\qquad 4\abs{\Escr_m}\leq 2^{-m}e^{\vartheta_m}e^{-\vartheta},
\end{equation}
where $\vartheta_m:=\vartheta+4p(m+1)$.

\emph{Step 2 (Bernstein bound for individual edges).}
Set $S_n(\tb):=\sum_{j=1}^n Z_j(\tb)$ and for $e:=(\ub,\vb)\in\Escr_m$, consider $T_n^{\SSadd{S}}:=S_n(\vb)-S_n(\ub)=\sum_{j=1}^n\SSadd{T_j^Z}$, where the increments~$\SSadd{T_j^Z}$ are given by $\SSadd{T_j^Z}:=Z_j(\vb)-Z_j(\ub)$. Since $Z_j(\cdot)$ is adapted to $\accolade{\Gscr_j}$, then~$\SSadd{T_j^Z}$ is $\Gscr_{\SSaddnew{j}}$-measurable. Moreover, it follows from~(ii) that almost surely,
\begin{equation*}
\E{\SSadd{T_j^Z}|\Gscr_{j-1}}=\E{Z_j(\vb)|\Gscr_{j-1}}-\E{Z_j(\ub)|\Gscr_{j-1}}=0,
\end{equation*}
hence $\accolade{\SSadd{T_j^Z}}_{j\geq 1}$ is a martingale-difference sequence with respect to $\accolade{\Gscr_j}$. On the other hand, $\normii{\vb-\ub}\leq 1$ holds for the pair $(\ub,\vb)\in\Escr_m$, hence by defining $R_m:=3\cdot 2^{-m}\Delta$, it follows from~(iii) that 
\begin{equation*}
\E{\left.\abs{\SSadd{T_j^Z}}^\ssmath\right\rvert\Gscr_{j-1}}\leq \frac{\ssmath !}{2}(B\Delta\normii{\vb-\ub})^{\ssmath-2}(V\Delta\normii{\vb-\ub})^2\leq \frac{\ssmath !}{2}B^{\ssmath-2}V^2R_m^{\ssmath}.
\end{equation*}
Thus, the increments of $\accolade{T_n^{\SSadd{S}}}_{n\geq 0}$ satisfy the conditional Bernstein moment condition with parameters $B^e:=BR_m$ and $V^e:=VR_m$. Consequently, Lemma~\ref{Lemma:52} applied with $A_m:=R_m\left(V\sqrt{\frac{2\vartheta_m}{\aleph}}+2B\frac{\vartheta_m}{\aleph}\right)$, yields almost surely,
\begin{equation}\label{eq:expothetam1}
\prob{\left.\underset{n\geq\aleph}{\sup}\frac{\abs{T_n^{\SSadd{S}}}}{n}>A_m\right\rvert\Gscr_0}\leq 2\exp\left(-\frac{\aleph A_m^2}{2(V^2R_m^2+BR_mA_m)}\right)\leq 2e^{-\vartheta_m},
\end{equation}
where the second inequality follows from the fact that 
\begin{equation}\label{eq:bigstar1}
\frac{\aleph A_m^2}{2(V^2R_m^2+BR_mA_m)}=\vartheta_m+\frac{\vartheta_mBR_m\sqrt{\frac{2V^2R_m^2\vartheta_m}{\aleph}}}{V^2R_m^2+BR_mA_m}\geq \vartheta_m
\end{equation}
established next. Let $X_m:=\sqrt{\frac{2V^2R_m^2\vartheta_m}{\aleph}}$ and $Y_m:=\frac{2BR_m\vartheta_m}{\aleph}$ satisfying $X_m+Y_m=A_m$, $\aleph X_m^2-2V^2R_m^2\vartheta_m=0$ and $\aleph Y_m-2BR_m\vartheta_m=0$. Then,
\begin{align*}
\aleph A_m^2-2\vartheta_m(V^2R_m^2+BR_mA_m)&=\aleph(X_m+Y_m)^2-2\vartheta_mV^2R_m^2-2\vartheta_mBR_m(X_m+Y_m)\\
&=(\aleph X_m^2-2V^2R_m^2\vartheta_m)+(\aleph Y_m-2BR_m\vartheta_m)Y_m+(2\aleph Y_mX_m-2\vartheta_mBR_mX_m)\\
&=2(2BR_m\vartheta_m)X_m-2\vartheta_mBR_mX_m=2\vartheta_mBR_m\sqrt{\frac{2V^2R_m^2\vartheta_m}{\aleph}},
\end{align*}
which yields~\eqref{eq:bigstar1}.

\emph{Step 3 (Simultaneous control of all edges).}
Taking a union bound over all admissible edges and all levels $m\geq 1$ using~\eqref{eq:expothetam1}, and making use of~\eqref{eq:Un}, gives
\begin{equation}\label{eq:period1}
\begin{split}
\prob{\left.\exists m\geq 1, e\in\Escr_m,n\geq \aleph:\abs{T_n^{\SSadd{S}}}>nA_m\right\rvert\Gscr_0}\leq 2\sum_{m=1}^{\infty}\abs{\Escr_m}e^{-\vartheta_m}\leq e^{-\vartheta}\sum_{m=1}^{\infty}2^{-m}=e^{-\vartheta},
\end{split}
\end{equation}
where we used $\sum_{m=1}^{\infty}2^{-m}=1$, so that outside the event in~\eqref{eq:period1}, every admissible edge satisfies 
\begin{equation*}
\abs{T_n^{\SSadd{S}}}\leq nA_m\qquad\text{for all } n\geq \aleph\ \text{and } m\geq 1.
\end{equation*}

\emph{Step 4 (Passage from edges to the whole field).}
For every $\tb\in\Bp$ and $\phi\geq 1$, using $S_n(\pi_0(\tb))=S_n(\bm{0})=0$, we have 
\begin{equation*}
S_n(\pi_\phi(\tb))=\sum_{m=1}^\phi\left(S_n(\pi_m(\tb))-S_n(\pi_{m-1}(\tb))\right).
\end{equation*}
Thus, outside the event in~\eqref{eq:period1}, by considering the edges $e=(\pi_{m-1}(\tb),\pi_m(\tb))$ and their associated $T_n^{\SSadd{S}}$, we have for all $n\geq \aleph$ and $\phi\geq 1$,
\begin{equation*}
\abs{S_n(\pi_\phi(\tb))}\leq \sum_{m=1}^\phi\abs{T_n^{\SSadd{S}}}\leq n\sum_{m=1}^\phi A_m\qquad\Longleftrightarrow\qquad\frac{\abs{S_n(\pi_\phi(\tb))}}{n}\leq \sum_{m=1}^\phi A_m.
\end{equation*}
Since $S_n(\cdot)$ is continuous, then letting $\phi\to\infty$ and using $\pi_\phi(\tb)\underset{\phi\to\infty}{\longrightarrow}\tb$, yields
\begin{equation*}
\frac{\abs{S_n(\tb)}}{n}\leq \sum_{m=1}^\infty A_m\qquad\text{for every } \tb\in\Bp\ \text{and } n\geq \aleph.
\end{equation*}
It remains to estimate $\sum_{m=1}^\infty A_m$ which is achieved next. By using the values of $A_m$ and $\vartheta_m$, we obtain 
\begin{equation}\label{eq:star1}
\sum_{m=1}^\infty A_m=3V\Delta\sqrt{\frac{2}{\aleph}}\sum_{m=1}^\infty 2^{-m}\sqrt{\vartheta_m}+\frac{6B\Delta}{\aleph}\sum_{m=1}^\infty 2^{-m}\vartheta_m,
\end{equation}
with 
\begin{equation*}
\sum_{m=1}^\infty 2^{-m}\vartheta_m=\sum_{m=1}^\infty (\vartheta+4p(m+1))2^{-m}=\vartheta+4p\sum_{m=1}^\infty(m+1)2^{-m}=\vartheta+12p,
\end{equation*}
while by Cauchy--Schwarz,
\begin{equation*}
\sum_{m=1}^\infty 2^{-m}\sqrt{\vartheta_m}\leq \left(\sum_{m=1}^\infty 2^{-m}\right)^{1/2}\left(\sum_{m=1}^\infty 2^{-m}\vartheta_m\right)^{1/2}\leq \sqrt{\vartheta+12p},
\end{equation*}
where we used $\sum_{m=1}^{\infty}2^{-m}=1$ and $\sum_{m=1}^{\infty}m2^{-m}=2$. Therefore,~\eqref{eq:star1} yields
\begin{equation*}
\begin{split}
\sum_{m=1}^\infty A_m&\leq 3V\Delta\sqrt{\frac{2}{\aleph}}\sqrt{\vartheta+12p}+\frac{6B\Delta}{\aleph}(\vartheta+12p).
\end{split}
\end{equation*}
Thus, \SSaddnew{using $\vartheta+12p\le12(p+\vartheta)$, $3\sqrt{24}<15$ and $6\cdot12=72$,} observing that 
\begin{equation*}
\frac{\abs{S_n(\tb)}}{n}\SSaddnew{<} 15V\Delta\sqrt{\frac{p+\vartheta}{\aleph}}+72B\Delta\frac{p+\vartheta}{\aleph}\qquad\text{for every } \tb\in\Bp\ \text{and } n\geq \aleph,
\end{equation*}
holds outside the event in~\eqref{eq:period1}, yields
\begin{equation*}
\prob{\left.\underset{n\geq \aleph}{\sup}\ \underset{\tb\in\Bp}{\sup}\ \frac{1}{n}\abs{\sum_{j=1}^{n}Z_j(\tb)}\geq 15V\Delta\sqrt{\frac{p+\vartheta}{\aleph}}+72B\Delta\frac{p+\vartheta}{\aleph}\right\rvert\mathscr{G}_0}\leq e^{-\vartheta},
\end{equation*}
which implies~\eqref{eq:chained2}, and the proof is completed.
\end{proof}
}

{\bl To apply Lemma~\ref{Lemma:Chained} to the \ASTRAR\ error field, we require continuity of the sample paths so that the chained control extends from the successive nets to the entire trust region. Under this additional assumption, Theorem~\ref{Theorem:Uniform} yields the desired uniform weak-accuracy bound at the adaptive sample size.}
\SSaddnew{\begin{assumption}\label{Assumption:cont}
The sample paths $\x\mapsto F(\x;\Theta)$ are continuous almost surely.
\end{assumption}

\SSaddnew{Assumption~\ref{Assumption:cont} holds in particular when the sample paths are Lipschitz, which is the pathwise form (almost surely, for all $\x,\y$) of the sufficient condition~\eqref{kappaFLipschitz1} for Assumption~\ref{Assumption:33}; in that case reuse costs nothing beyond Assumption~\ref{Assumption:33}.}

}

\SSadd{\begin{theorem}[uniform weak accuracy]\label{Theorem:Uniform}
Let Assumptions~\ref{Assumpt:31},~\ref{Assumption:Qmax},~\ref{Assumption:33} and~\ref{Assumption:cont} hold, consider \ASTRAR\ (Algorithm~\ref{algoR}) 
and define
\begin{equation}\label{eq:rF}
r_F:=\min\accolade{1,\ \parenth{\frac{\sigma_0\kfde}{\SSaddnew{30}\kappaas\sigma_D\Qmax}}^2,\ \frac{\sigma_0^2\kfde}{\SSaddnew{144}\kappaas^2 b_D\Qmax\deltmax}},\qquad c_F:=\frac{r_F}{2}.
\end{equation}
Then for every $k\geq 0$ with $\lambda_k\geq 2p/r_F$, almost surely,
\begin{equation}\label{eq:uniformAcc}
\prob{\left.\underset{\normii{\s}\leq\Dk}{\sup}\ \abs{\Dbar_k(\s,\nkrandom)}> \kfde\Dk^2\right\rvert\fqalgebra}\leq e^{-c_F\lambda_k}.
\end{equation}
\end{theorem}}

{\bl
\begin{proof}
We apply Lemma~\ref{Lemma:Chained} with $\mathscr{G}_0=\fqalgebra$, $\gjalgebra=\gqalgebra{k,j}$, $\Delta=\Dk$, $\aleph=\Nkbar$, $\vartheta=c_F\lambda_k$, and consider, for every $\tb\in\Bp$, the normalized random field $Z_j(\tb):=D(\Xzerok,\Xzerok+\QkRandom\Dk\tb;\Theta_{k,j})$, so that $\widetilde{Z}_j(\s):=Z_j(\s/\Dk)=D(\Xzerok,\Xzerok+\QkRandom\s;\Theta_{k,j})$, with $\normii{\s}\le\Dk$, satisfies $\frac1n\sum_{j\le n}\widetilde{Z}_j(\s)=\Dbar_k(\s,n)$. The quantities $\Dk$, $\QkRandom$, $\Xzerok$ and $\Nkbar$ are $\fqalgebra$-measurable. 

Hypothesis~(i) holds because $Z_j(\bm0)=0$ by~\eqref{eq:DIncrement}, and $\tb\mapsto Z_j(\tb)$ is continuous by Assumption~\ref{Assumption:cont} and the continuity of $f$. \SSaddnew{For hypothesis (ii) fix $\tb\in\Bp$; for (iii) fix $\tb,\tb'\in\Bp$ with $\normii{\tb-\tb'}\le1$. Then }
\eqref{eq:Dtelescope} gives $Z_j(\tb)-Z_j(\tb')=D(\Xzerok+\QkRandom\Dk\tb',\Xzerok+\QkRandom\Dk\tb;\Theta_{k,j})$, a paired increment at displacement $$\normii{\QkRandom\Dk(\tb-\tb')}\leq\Qmax\Dk\normii{\tb-\tb'}\leq\Qmax\deltmax.$$ 

Since $\Xzerok$, $\QkRandom$ and $\Dk$ are $\gqalgebra{k,j-1}$-measurable and $\Theta_{k,j}$ is independent of $\gqalgebra{k,j-1}$, Lemma~\ref{Lemma:24} and Assumption~\ref{Assumption:33} give $\E{Z_j(\tb)|\gqalgebra{k,j-1}}=0$ and~(iii) with $B=b_D\Qmax$ and $V=\sigma_D\Qmax$. 

Since $\Nkbar\ge\sigma_0^2\lambda_k\kappaas^{-2}\Dk^{-2}$ and $\Dk\le\deltmax$,
\[
15V\Dk\sqrt{\tfrac{p+\vartheta}{\Nkbar}}+72B\Dk\tfrac{p+\vartheta}{\Nkbar}
\le\Dk^2\Bigl[\tfrac{15\kappaas\sigma_D\Qmax}{\sigma_0}\sqrt{\tfrac{p+\vartheta}{\lambda_k}}+\tfrac{72\kappaas^2b_D\Qmax\deltmax}{\sigma_0^2}\,\tfrac{p+\vartheta}{\lambda_k}\Bigr]\le\kfde\Dk^2 .
\]
Here $(p+\vartheta)/\lambda_k\le r_F$ by the hypothesis $\lambda_k\ge2p/r_F$, and each bracketed term is then at most $\kfde/2$ by~\eqref{eq:rF}. Since $\nkrandom\ge\Nkbar$, the event in~\eqref{eq:uniformAcc} is contained in the event of~\eqref{eq:chained2}, whose conditional probability is at most $e^{-\vartheta}=e^{-c_F\lambda_k}$.
\end{proof}}
 
\SSadd{\begin{remark}\label{rem:lambdaIndex}
The hypothesis $\lambda_k\geq 2p/r_F$ is a condition on an index, not on $\ubar{\lambda}$: since $\lambda_k\to\infty$ it holds beyond a deterministic $k_\lambda$ for every $\ubar{\lambda}>0$, and so does every other threshold of this kind below, including~\eqref{eq:zetastar}. With a constant schedule each such threshold would instead constrain an input through unknown constants, and could be missed.
\end{remark}}

{\bl The uniform bound of Theorem~\ref{Theorem:Uniform} can now be specialized simultaneously to the design points and the data-dependent trial point. The following result collects these consequences in the event~$\uk$ and expresses its conditional probability at the level used later in the analysis.}
\SSadd{\begin{corollary}\label{cor:UniformEvent}
Let the assumptions of Theorem~\ref{Theorem:Uniform} hold and set
\[
\uk:=\accolade{\underset{\normii{\s}\leq\Dk}{\sup}\ \abs{\Dbar_k(\s,\nkrandom)}\leq \kfde\Dk^2}.
\]
Then $\prob{\left.\uk\right\rvert\fqalgebra}\geq 1-e^{-c_F\lambda_k}$ almost surely \SSaddnew{for every $k$ with $\lambda_k\ge2p/r_F$}, and on $\uk$ both conclusions used in the sequel hold simultaneously: $\abs{\erandombar^i_k(\nkrandom)-\erandombar^0_k(\nkrandom)}\leq\kfde\Dk^2$ for every $i\in\intbracket{1,q}$, and $\abs{\erandombar^{s}_k(\nkrandom)-\erandombar^0_k(\nkrandom)}\leq\kfde\Dk^2$. Consequently, for any $\ef>0$ there is a deterministic $\kappaf\in\N\SSaddnew{, \kappaf\ge k_\lambda}$ with
\begin{equation}\label{eq:UniformBeta}
\prob{\left.\uk\right\rvert\fqalgebra}\geq \beta_{k,\ef}=1-\frac{1}{k^{1+\ef}}\qquad\text{for all }k\geq\kappaf,\ \text{almost surely}.
\end{equation}
\end{corollary}}
 
\SSadd{\begin{proof}
The probability bound is Theorem~\ref{Theorem:Uniform}. The two displayed conclusions are the specialisations of the event $\uk$ to $\s=\siirandom_k$, $i\in\intbracket{1,q}$, and to $\s=\srandom_k$, both admissible since $\normii{\siirandom_k}=\Dk$ and $\normii{\srandom_k}\leq\Dk$. Finally, for $\ef>0$ fixed, the argument of Corollary~\ref{theor31bis} applies verbatim with $\chifde$ replaced by $c_F$: since $(\log(k+2))^{1+\epsLanbda}/\log k\to\infty$, there is $\kappaf\in\N$ with $c_F\ubar{\lambda}(\log(k+2))^{1+\epsLanbda}\geq(1+\ef)\log k$ for all $k\geq\kappaf$, which gives~\eqref{eq:UniformBeta}.
\end{proof}}

\SSadd{\begin{remark}\label{rem:UniformReplaces}
The two conclusions of Corollary~\ref{cor:UniformEvent} are exactly the inputs that Theorem~\ref{mgeTheorem} and Lemma~\ref{lemma4p7} need from the estimates
, namely $(\kmge,\Qk)$-\SSaddnew{subspace-accuracy} of the model and accuracy of the trial estimates. 
Every later use of those results therefore goes through with the single event $\uk$ in place of the pair $\ik\cap\jk$, with all conditioning on $\fqalgebra$, and with $\fmqalgebra$ absent. 
\SSaddnew{Since one event replaces two, the lower bound $\pi_{\min}^2$ on the conditional probability of a good iteration, obtained in \eqref{tower1} as a product of two bounds, becomes $\pi_{\min}$. The consequences for the constants are recorded in Remark~\ref{rem:theorem32R}.}
\end{remark}}

\subsection{\SSaddnew{Subspace accuracy} and well-alignment}
{\bl A crucial model requirement for the analysis is the accurate approximation of the reduced gradient over the trust region. Unlike the subspace fully linear models of STARS~\cite[Definition~3.10]{DzaWildSub2022}, which control both model values and gradients, we retain only the gradient component, motivated also by the reduced model requirements considered in~\cite[Section~1.3.2(c)]{ha2025complexity}. This is sufficient here because the acceptance test depends on differences, and the resulting gradient control yields the required bound on the error in predicted reduction through Lemma~\ref{modFuncErrLemma}. Assumption~\ref{fminAndLipGradAssump} provides the smoothness needed for this analysis, Definition~\ref{def:semilinear} formalizes the resulting notion of subspace-accuracy, and Theorem~\ref{mgeTheorem} shows that the models constructed in Section~\ref{sec:estiModels} satisfy it with high conditional probability.}
\begin{assumption}\label{fminAndLipGradAssump}
There exists a constant $f_{\min}>-\infty$ such that $f(\x)\geq f_{\min}$ for all $\x\in\rd$. 
Moreover, $f$ is continuously differentiable with $\Lg$-Lipschitz continuous gradient, for some positive constant $\kappa_{Lg}<\infty$.
\end{assumption}
 
 \SSaddnew{\begin{definition}\label{def:semilinear}
Given $\Qk\in\rdp$ and $\xzerok\in\rd$, a function $\mhatk:\rp\to\R$ satisfying
\begin{equation}\label{semiLin1}
\normii{\QkTranspose\nabla f(\xzerok+\Qk\s)-\nabla\mhatk(\s)}\leq \frac{1}{2}\kmge\dk,\quad\text{for all}\ \s\in \mathscr{B}(\bm{0},\dk)\subset\rp
\end{equation}
for some constant $\kmge>0$ and $\dk>0$, is called a $\left(\kmge,\Qk\right)$-subspace-accurate model of~$f$ in $\mathscr{B}(\xzerok,\dk;\Qk)$.
\end{definition}
}
 
\begin{theorem}\label{mgeTheorem}
Let Assumptions~\ref{fminAndLipGradAssump}, \ref{Assumpt:31}, \ref{Assumption:Qmax} and~\ref{Assumption:33} hold, and consider the random models~$\mhatkRandom$ constructed in Section~\ref{sec:estiModels}. Introduce the events 
\begin{equation*}
\ik:=\underset{\normii{\s}\leq \Dk}{\bigcap}\accolade{\normii{\QkRandomTranspose\nabla f(\Xzerok+\QkRandom\s)-\nabla\mhatkRandom(\s)}\leq\frac{1}{2}\kmge\Dk},
\end{equation*}
where $\kmge:=2\left(\Lg\Qmax^2+\kappah+\sqrt{p}\left(\frac{1}{2}\Lg\Qmax^2+\kfde\right)\right)$. Then, for all~$k$ sufficiently large and any $\ef>0$ fixed,
\begin{equation}\label{mgeEquation1}
\prob{\ik\left\lvert\fqalgebra\right.}=\E{\left.\iik\right\rvert\fqalgebra}\geq 1-\frac{2p}{k^{1+\ef}}\quad \text{almost surely}.
\end{equation}
\end{theorem}

\begin{proof}
We consider the reduced objective $\fkhat:\s\mapsto\fkhat(\s):=f(\Xzerok+\QkRandom\s)$, $\s\in\rp$, and note that under Assumptions~\ref{Assumption:Qmax} and~\ref{fminAndLipGradAssump}, for all $\s_1,\s_2\in\rp$, 
\begin{equation*}
\begin{split}
\normii{\nabla\fkhat(\s_2)-\nabla\fkhat(\s_1)}&=\normii{\QkRandomTranspose\left(\nabla f(\Xzerok+\QkRandom\s_2)-\nabla f(\Xzerok+\QkRandom\s_1)\right)}\\
&\leq \Lg\norme{\QkRandom}\normii{\QkRandom(\s_2-\s_1)}\leq \Lg\Qmax^2\normii{\s_2-\s_1},
\end{split}
\end{equation*}
that is, $\fkhat$'s gradient is $\Lg\Qmax^2$-Lipschitz continuous. For each coordinate $i\in\intbracket{1,p}$, recalling the interpolation points $\siirandom_k=\Dk\ebhat^i$ and $\srandom^{p+i}_k=-\Dk\ebhat^i$, we consider the central finite-difference approximation of $\partial_i\fkhat(\bm{0})$ given by
\begin{equation}\label{eq:Gkjf}
G_{k,i}^f:=\frac{\fkhat(\Dk\ebhat^i)-\fkhat(-\Dk\ebhat^i)}{2\Dk}=\frac{1}{2\Dk}\int_{-\Dk}^{\Dk}\partial_i\fkhat(t\ebhat^i)dt,
\end{equation}
where the second equality follows from the fundamental theorem of calculus. Then,
\begin{equation}\label{eq:dot1}
\begin{split}
\abs{G_{k,i}^f-\partial_i\fkhat(\bm{0})}&\leq \frac{1}{2\Dk}\int_{-\Dk}^{\Dk}\abs{\partial_i\fkhat(t\ebhat^i)-\partial_i\fkhat(\bm{0})}dt\leq \frac{1}{2\Dk}\int_{-\Dk}^{\Dk}\abs{(\ebhat^i)^\top\left(\nabla\fkhat(t\ebhat^i)-\nabla\fkhat(\bm{0})\right)}dt\\
&\leq \frac{1}{2\Dk}\int_{-\Dk}^{\Dk}\Lg\Qmax^2\abs{t}dt=\frac{1}{2}\Lg\Qmax^2\Dk.
\end{split}
\end{equation}
On the other hand, observing that $\Fbarik(\srandom^i_k,\nkrandom;\Theta_k)=\fkhat(\Dk\ebhat^i)+\erandombar^i_k(\nkrandom)$ and $\Fbar^{p+i}_k(\srandom^{p+i},N_k;\Theta_k)=\fkhat(-\Dk\ebhat^i)+\erandombar^{p+i}_k(\nkrandom)$, it follows from~\eqref{eq:Ghatki} and~\eqref{eq:Gkjf} that 
\begin{equation*}
\GhatRandom_{k,i}-G_{k,i}^f=\frac{\erandombar^i_k(\nkrandom)-\erandombar^{p+i}_k(\nkrandom)}{2\Dk}.
\end{equation*}
 
Thus, on $\overset{2p}{\underset{j=1}{\bigcap}}\mathscr{E}^j_k$ with $\mathscr{E}^j_k:=\accolade{\abs{\erandombar^j_k(\nkrandom)-\erandombar^0_k(\nkrandom)}\leq \kfde\Dk^2}$, $j\in\intbracket{1,2p}$, it holds that for each $i\in\intbracket{1,p}$,
\begin{equation*}
\begin{split}
\abs{\GhatRandom_{k,i}-G_{k,i}^f}&\leq \frac{1}{2\Dk}\left(\abs{\erandombar^i_k(\nkrandom)-\erandombar^0_k(\nkrandom)}+\abs{\erandombar^{p+i}_k(\nkrandom)-\erandombar^0_k(\nkrandom)}\right)\leq \kfde\Dk,
\end{split}
\end{equation*}
which, combined with~\eqref{eq:dot1}, implies 
\begin{equation*}
\abs{\GhatRandom_{k,i}-\partial_i\fkhat(\bm{0})}\leq \abs{\GhatRandom_{k,i}-G_{k,i}^f}+\abs{G_{k,i}^f-\partial_i\fkhat(\bm{0})}\leq \left(\frac{1}{2}\Lg\Qmax^2+\kfde\right)\Dk.
\end{equation*}
Therefore, 
\begin{equation}\label{eq:GhatNablafhat}
\normii{\GhatRandom_k-\nabla\fkhat(\bm{0})}^\SSaddnew{2}=\sum_{i=1}^{p}\abs{\GhatRandom_{k,i}-\partial_i\fkhat(\bm{0})}^2\leq p\left(\frac{1}{2}\Lg\Qmax^2+\kfde\right)^2\Dk^2.
\end{equation}
Recalling~\eqref{eq:randoModel}, we have for all $\s\in\mathscr{B}(\bm{0},\Dk)$
\begin{equation*}
\begin{split}
\normii{\nabla\fkhat(\s)-\nabla\Mhatk(\s)}&\leq \normii{\nabla\fkhat(\s)-\nabla\fkhat(\bm{0})}+\normii{\GhatRandom_k-\nabla\fkhat(\bm{0})}+\normeshort{\HkhatRandom}\normii{\s}\\
&\leq \Lg\Qmax^2\Dk+\sqrt{p}\left(\frac{1}{2}\Lg\Qmax^2+\kfde\right)\Dk+\kappah\Dk=\frac{1}{2}\kmge\Dk,
\end{split}
\end{equation*}
where we used the $\Lg\Qmax^2$-Lipschitz continuity of $\nabla\fkhat$ combined with~\eqref{eq:GhatNablafhat} and~\eqref{eq:HessianBound}.
 
Observing that $\nabla\fkhat(\s)=\QkRandomTranspose\nabla f(\Xzerok+\QkRandom\s)$, we have thus proved that $\overset{2p}{\underset{i=1}{\bigcap}}\mathscr{E}^i_k\subseteq\ik$, with the $\mathscr{E}^i_k$ satisfying $\prob{\bar{\mathscr{E}}^i_k |\fqalgebra}\leq\frac{1}{k^{1+\ef}}$ for~$k$ sufficiently large and any $\ef>0$ thanks to Corollary~\ref{theor31bis}. Hence, almost surely
\begin{equation*}
\begin{split}
\prob{\left. \ik \right\rvert\fqalgebra}&\geq \prob{\left.\overset{2p}{\underset{i=1}{\bigcap}}\mathscr{E}^i_k \right\rvert\fqalgebra}\geq 1-\sum_{i=1}^{2p}\prob{\bar{\mathscr{E}}^i_k |\fqalgebra}\geq 1-\frac{2p}{k^{1+\ef}}, 
\end{split}
\end{equation*}
which completes the proof.
 
\end{proof}

 {\bl To connect reduced-space progress with full-space stationarity, the random subspace must retain a nonnegligible component of the full-space gradient. The corresponding standard well-alignment condition was formalized in \cite[Definition~2]{CRsubspace2021}, while \cite[Section~3.1, Definition~3.2]{DzaWildSub2022} introduced its conditional probabilistic version 
 adopted here. The following definition states the condition needed below.}
\begin{definition}\label{def:WAM}(Inspired by~\cite[Definition~3.2]{DzaWildSub2022})
A sequence $\accoladekinN{\QkRandom}$ of random matrices is called $(1-\beta_Q)$-probabilistically $\alpha_Q$-well aligned, for fixed $\alpha_Q,\beta_Q\in (0,1)$, if for any $\falgebra$-measurable random vector $\vrandomvec$ with realizations $\vi\in\rd$, the events $\AkEvent$ defined by $\AkEvent:=\accolade{\normii{\QkRandomTranspose\vrandomvec}\geq \alpha_Q\normii{\vrandomvec}}$ satisfy almost surely
\begin{equation}\label{WAMEquation}
\prob{\AkEvent|\falgebra}=\E{\iak|\falgebra}\geq 1-\beta_Q,
\end{equation}
where $\iak$ is the indicator function of $\AkEvent$, satisfying $\iak(\omega)=1$ if $\omega\in \AkEvent$ and $\iak(\omega)=0$ otherwise. Random matrices satisfying~\eqref{WAMEquation} are referred to as $WAM(\alpha_Q,\beta_Q)$.
\end{definition} 
{\bl Several random matrix ensembles satisfy Definition~\ref{def:WAM}. In particular, when $\QkRandom$ is generated independently of $\falgebra$, any Johnson--Lindenstrauss transform~\cite{JohnsonLind1984} satisfying $\normii{\QkRandom^\top \vb}\geq \alpha_Q\normii{\vb}$ with probability at least $1-\beta_Q$ for every fixed~$\vb$ yields~\eqref{WAMEquation}. This includes Gaussian matrices~\cite[Corollary~3.5]{DzaWildSub2022}, $s$-hashing matrices~\cite[Corollary~3.8]{DzaWildSub2022} based on the sparse constructions of~\cite{KaNel2014SparseLidenstrauss}, $s$-hashing-like matrices~\cite[Theorem~2.7]{DzaWildHashing2022}, and suitably scaled submatrices of Haar-distributed orthogonal matrices~\cite[Theorem~3.9]{DzaWildSub2022} (see also~\cite[Theorem~2.2]{DzaWildSubDir24}). The $s$-hashing, $s$-hashing-like, and Haar constructions also admit deterministic bounds on $\norme{\QkRandom}$ and hence satisfy Assumption~\ref{Assumption:Qmax}, whereas Gaussian matrices satisfy such bounds only with high probability. Since these constructions and their parameter choices have been studied in the aforementioned works, we focus here only on the properties required by Definition~\ref{def:WAM} and Assumption~\ref{Assumption:Qmax}.}

\section{Convergence analysis}\label{sec:convAnalysis}

{\bl 
The convergence analysis adapts the supermartingale argument of~\cite{DzaWildSub2022} to the accuracy mechanisms established in Section~\ref{sec:Prob}. Algorithms~\ref{algoF} and~\ref{algoR} share the same Lyapunov function~$\Phik$, trust-region radius dynamics, and deterministic ingredients, in particular Lemmas~\ref{lemmaIncr} and~\ref{lemma4p7}, which respectively control objective increases on accepted steps and force success when the radius is sufficiently small. Their analyses differ only in how model and estimate accuracy enter Theorem~\ref{convDeltaTheorem}: for \ASTRAF, the model and estimate events are distinct, leading to three subcases in the first case and an auxiliary $\fmqalgebra$-measurable event in the second; for \ASTRAR, both properties follow from the single event~$\uk$, reducing each case to two branches and eliminating~$\fmqalgebra$. Consequently, the product bound $\SSaddnew{\pi_{\min}^2}$ associated with the two separate accuracy events is replaced by the single-event bound $\prob{\uk|\fqalgebra}\geq\pi_{\min}$ for \ASTRAR. We present the proof for \ASTRAR\ and defer the \ASTRAF\ proof, which closely follows~\cite{DzaWildSub2022}, to Appendix~\ref{app:proofConvF}.
}
 
\subsection{Preliminary results}\label{sec:sec5pt1}
{\bl The following preliminary results provide the deterministic ingredients needed for the drift analysis of Theorem~\ref{convDeltaTheorem}. Lemmas~\ref{descentLemma}--\ref{lemma4p6} relate smoothness and model subspace-accuracy to the error in predicted reduction and to actual decrease in~$f$ when the trust-region radius is sufficiently small; Lemmas~\ref{lemma4p7}--\ref{lemmaIncr} then use weakly accurate estimates to characterize successful iterations, guarantee decrease on good accepted steps, and control the possible increase on erroneous ones. Lemma~\ref{lemmafourp7} records the measurability needed for the subsequent conditional case analysis.}
 
\begin{lemma}\label{descentLemma}
Let $f$ be differentiable with $\kappa_{Lg}$-Lipschitz continuous gradient, let $\s\in\rp$, $\xzerok\in\rd$ and $\Qk\in\rdp$. Under Assumption~\ref{Assumption:Qmax}, it holds that $\abs{f(\xk+\Qk\SSaddnew{\s})-f(\xk)-(\Qk\SSaddnew{\s})^\top\nabla f(\xk)}\leq \frac{1}{2}\kappa_{Lg}\Qmax^2\normii{\s}^2$.
\end{lemma}

 \begin{lemma}\label{modFuncErrLemma}
Let Assumptions~\ref{fminAndLipGradAssump} and~\ref{Assumption:Qmax} hold. Suppose that $\mhatk$ is a $(\kmge,\Qk)$-\SSaddnew{subspace-accurate} model of $f$ in $\mathscr{B}(\xzerok,\dk;\Qk)$. 
Let $\kmfe:=\kmge+\kappa_{Lg}\Qmax^2+\SSaddnew{\kappah}$. Then 
\begin{align}\label{modFuncErrEq}
\abs{(\mk(\sk)-f(\xk+\Qk\sk))-(\mk(\bm{0})-f(\xk))}\leq\frac{1}{2}\kmfe\dk^2.
\end{align}
\end{lemma}
\begin{proof} We can write 
\begin{eqnarray*}
\Delta^{\mk\cdot f}&:=&\abs{(\textcolor{black}{\mk(\sk)-\mk(\bm{0})})-(f(\xk+\Qk\sk)-f(\xk))} \\
&=& \abs{\textcolor{black}{\gk^\top\sk+\frac{1}{2}\sk^\top\Hk\sk}-\left(f(\xk+\Qk\sk)-f(\xk)-(\Qk\sk)^\top\nabla f(\xk)\right)-(\Qk\sk)^\top\nabla f(\xk)}\\
&\leq& \abs{(\gk-\QkTranspose\nabla f(\xk))^\top\sk}+\abs{\textcolor{black}{f(\xk+\Qk\sk)-f(\xk)-(\Qk\sk)^\top\nabla f(\xk)}}+\frac{1}{2}\abs{\sk^\top\Hk\sk}\\
&\leq& \frac{1}{2}\kmge\dk^2+\textcolor{black}{\frac{1}{2}\kappa_{Lg}\Qmax^2\dk^2}+\frac{1}{2}\SSaddnew{\kappah}\dk^2,
\end{eqnarray*}
where the last inequality follows from the \SSaddnew{subspace-accuracy} assumption, Lemma~\ref{descentLemma} and~\eqref{eq:HessianBound}, respectively, 
which completes the proof.
\end{proof}

\begin{lemma}\label{lemma4p5}(\cite[Lemma~4.1]{DzaWildSub2022})
Let all the assumptions made in Lemma~\ref{modFuncErrLemma} and Assumption~\ref{CauchyDecreaseAssump} hold. If 
\begin{align}\label{firstDeltaCond}
\dk\leq\min\accolade{\frac{1}{\SSaddnew{\kappah}},\frac{\kfcd}{2\kmfe}}\normii{\gk},
\end{align}
then the trial step $\sk$ leads to an improvement in $f$ such that
\begin{align}\label{decrLemma4p5}
f(\xk+\Qk\sk)-f(\xk)\leq -\frac{1}{4}\kfcd\normii{\gk}\dk.
\end{align}
\end{lemma}
\SSadd{The proof is given in Appendix~\ref{app:proofsSec5}.}
 
\begin{lemma}\label{lemma4p6}(\cite[Lemma~4.2]{DzaWildSub2022})
Let all the assumptions made in Lemma~\ref{lemma4p5} hold. If 
\begin{align}\label{secDeltaCond}
\dk\leq\min\accolade{\frac{1}{\SSaddnew{\kappah}+\frac{1}{2}\kmge},\frac{1}{\frac{2\kmfe}{\kfcd}+\frac{1}{2}\kmge}}\normii{\QkTranspose\nabla f(\xk)},
\end{align}
then the trial step $\sk$ leads to an improvement in $f$ such that
\begin{align}\label{seconImprov}
f(\xk+\SSaddnew{\Qk}\sk)-f(\xk)\leq -c_1\normii{\QkTranspose\nabla f(\xk)}\dk,
\end{align}
where $c_1=\frac{1}{2}\kfcd\max\accolade{\frac{\SSaddnew{\kappah}}{2\SSaddnew{\kappah}+\kmge},\frac{2\kmfe}{4\kmfe+\kmge\kfcd}}$.
\end{lemma}
\SSadd{The proof is given in Appendix~\ref{app:proofsSec5}.}
 
\begin{lemma}\label{lemma4p7}(\cite[Lemma~4.3]{DzaWildSub2022}) 
Let all the assumptions made in Lemma~\ref{lemma4p5} hold, and assume further that the estimates $\accolade{\fbarokcur,\fsk}$ are $\kfde$-weakly accurate. 
f 
\begin{align}\label{thirdDeltCond}
\dk\leq\min\accolade{\mu,\frac{1}{\SSaddnew{\kappah}},\frac{\kfcd(1-\eta)}{2(\kfde+\kmfe)}}\normii{\gk},
\end{align}
then the $k$th iteration is successful.
\end{lemma}
\SSadd{The proof is given in Appendix~\ref{app:proofsSec5}.}
 
\begin{lemma}\label{lemma4p8}(\cite[Lemma~4.4]{DzaWildSub2022})
\SSaddnew{Let Assumption~\ref{CauchyDecreaseAssump} and~\eqref{eq:HessianBound} hold, let the estimates be $\kfde$-weakly accurate and the iteration be successful.} Assume further that $\kfde<\frac{\eta}{2\mu}\kfcd\min\accolade{\frac{1}{\SSaddnew{\kappah}\mu},1}$.
Define $c_2=\frac{\eta}{2\mu}\kfcd\min\accolade{\frac{1}{\SSaddnew{\kappah}\mu},1}-\kfde>0$. Then the improvement in~$f$ is bounded as
\begin{align}\label{eq:fourpt5}
f(\xkun)-f(\xk)\leq -c_2\dk^2.
\end{align}
\end{lemma}
\SSadd{The proof is given in Appendix~\ref{app:proofsSec5}.}
 
\begin{lemma}\label{lemmaIncr}({\bl Inspired by} \cite[Lemma~4.5]{DzaWildSub2022})
Under Assumptions~\ref{fminAndLipGradAssump} and~\ref{Assumption:Qmax}, assume that \SSaddnew{Algorithm~\ref{algoF} or~\ref{algoR}}  
erroneously accepts a step $\sk$ yielding an increase in~$f$ when $\normii{\QkTranspose\nabla f(\xzerok)}\geq \zeta\dk$ for some constant $\zeta>0$, and define $c_3(\zeta)=1+\frac{3\kappa_{Lg}}{2\zeta}\Qmax^2$. Then the increase in~$f$ is such that
\begin{equation}\label{fourpt6}
f(\xzerok+\Qk\sk)-f(\xk)\leq c_3(\zeta)\normii{\QkTranspose\nabla f(\xzerok)}\dk.
\end{equation}
{\bl Likewise, if $\normii{\QkTranspose\nabla f(\xzerok)}< \zeta\dk$, then 
\begin{equation}\label{fourpt6bis}
f(\xzerok+\Qk\sk)-f(\xk)< c_3(\zeta)\zeta\dk^2.
\end{equation}
}
\end{lemma}
{\bl The proof is given in Appendix~\ref{app:proofsSec5}.}
 
\begin{lemma}\label{lemmafourp7}(\cite[Lemma~4.7]{DzaWildSub2022})
Let $\zeta\geq 0$ be a constant, and consider the event $\mathscr{V}_k:=\SSaddnew{\left\{\normii{\QkRandomTranspose\nabla f(\Xzerok)}\ge\zeta\Dk\right\}}$. 
Then $\SSaddnew{\normii{\QkRandomTranspose\nabla f(\Xzerok)}}$, $\ivk$ and $\ivkc$ are $\fqalgebra$-measurable.
\end{lemma}

\subsection{
\SSaddnew{Global} 
convergence {\bl of the sequence of trust region radii}}
{\bl We now establish the global convergence of the trust-region radii by combining the probabilistic subspace-accuracy and weak-accuracy guarantees of Section~\ref{sec:Prob} with the descent, success, and objective-increase bounds of Section~\ref{sec:sec5pt1}. Inspired by~\cite{audet2019stomads,blanchet2016convergence,chen2018stochastic,dzahini2020expected,dzahini2020constrained,DzaWildSubDir24,DzaWildSub2022,paquette2018stochastic,rinaldi2024stochastic}, Theorem~\ref{convDeltaTheorem} introduces a suitable Lyapunov function and shows that its conditional expected change is bounded above by a negative multiple of $\Dk^2$ for both \ASTRAF\ and \ASTRAR. This yields the summability of $\accoladekinN{\Dk^2}$ and, consequently, $\Dk\to 0$ almost surely.}
\begin{theorem}\label{convDeltaTheorem}\SSadd{(\ASTRAF\ and \ASTRAR)}
\SSaddnew{Let Assumptions \ref{CauchyDecreaseAssump}, \ref{Assumpt:31}, \ref{Assumption:Qmax}, \ref{Assumption:33}, \ref{fminAndLipGradAssump} and, for \ASTRAR, Assumption~\ref{Assumption:cont} hold, let $\kfde$ satisfy the condition of Lemma \ref{lemma4p8}, and let $c_1,c_2,c_3:=c_3(\zeta)$ be the constants of Lemmas \ref{lemma4p6}, \ref{lemma4p8} and \ref{lemmaIncr}}

for some fixed $\zeta>0$ satisfying 
\begin{equation}\label{fourpt8}
\zeta\geq \frac{1}{2}\kmge +\SSaddnew{\kappa_{\rm suc}\quad\text{ where }\quad\kappa_{\rm suc}:=\max\accolade{\frac{1}{\mu},\SSaddnew{\kappah},\frac{2(\kfde+\kmfe)}{\kfcd(1-\eta)}}.}
\end{equation}
Let $\nu\in (0,1)$ be chosen according to 
\begin{equation}\label{fourpt9}
\frac{\nu}{1-\nu}\geq \frac{8\gamma^2}{\min\accolade{c_2,\zeta c_1}}.
\end{equation}
\SSaddnew{Let $\pi_{\min}\in(0,1)$ and set $b_1:=(1-\nu)(1-\gamma^{-2})$ and $\varrho:=\frac12\pi_{\min}b_1>0$. For \ASTRAF, let $\pi_{\min}$ satisfy
\begin{equation}\label{fourpt10}
\frac{\SSaddnew{\pi_{\min}^2}-
\hat\gamma_r}{2(1-\SSaddnew{\pi_{\min}})}\geq \frac{c_3}{c_1}\qquad\text{and}\qquad \frac{\SSaddnew{\pi_{\min}}}{1-\SSaddnew{\pi_{\min}}}\geq \frac{2\nu\zeta\gamma^2c_3}{(1-\nu)(\gamma^2-1)}+2\gamma^2,
\end{equation}
for $r>1/(1-\epsLanbda)$ fixed,  $\hat\gamma_r:=\max\{\gamma^{2r}/(1+\gamma^{2r}),\frac12\}$,  
$\ef$ of Corollary~\ref{theor31bis}, and  $k_0:=\kappaf\vee\lceil(2p/(1-\pi_{\min}))^{1/(1+\ef)}\rceil$. 

For \ASTRAR, $\pi_{\min}$ is arbitrary and $k_0\ge1$ is such that, for all $k\ge k_0$, $\lambda_k\ge2p/r_F$ and\footnote{\SSaddnew{Since $\lambda_k\to\infty$, both conditions hold for all $k$ beyond a deterministic index, for every $\ubar\lambda>0$ (Remark~\ref{rem:lambdaIndex}).}}}
\begin{equation}\label{eq:zetastar}
\varpi_k:=e^{-c_F\lambda_k}\le\min\{\varpi_\star,1-\pi_{\min}\},\ \varpi_\star:=\min\Bigl\{\frac{c_1}{2(c_1+c_3)},\ \frac{b_1}{2(b_1+b_4)}\Bigr\},\ \ b_4:=\nu c_3\zeta+(1-\nu)(\gamma^2-1).
\end{equation}
Then the random function $\Phik:=\nu(f(\Xzerok)-f_{\min})+(1-\nu)\Dk^2$ satisfies
\begin{equation}\label{fourpt11}
\E{\Phikplusone-\Phik|\fqalgebra}\leq -\varrho\Dk^2\quad\text{for all}\ k\geq k_0,
\end{equation}
which implies that, 
\begin{equation}\label{fourpttwelve}
\E{\sum_{k=0}^{\infty}\Delta_k^2}\leq k_0\delta_{\max}^2+ \frac{1}{\varrho}\E{\Phi_{k_0}},
\end{equation}
hence
\begin{equation}\label{fourpt12}
\sum_{k=0}^{\infty}\Dk^2<\infty\quad\text{almost surely.}
\end{equation}
\end{theorem}
 
\SSadd{The proof for \ASTRAF, which follows the three-subcase argument of~\cite{DzaWildSub2022} closely, is given in Appendix~\ref{app:proofConvF}; the proof for \ASTRAR\ follows.}

 \SSaddnew{\paragraph{Preliminaries common to the two proofs.} Denote by $\phik$ the realizations of $\Phik$. On a successful iteration $\xkun=\xzerok+\Qk\sk$ and $\dkun=\min\accolade{\gamma\dk,\deltmax}$, so that
\begin{equation}\label{fourpt13}
\phikplusone-\phik\leq \nu\bigl(f(\xkun)-f(\xzerok)\bigr)+(1-\nu)(\gamma^2-1)\dk^2,
\end{equation}
while on an unsuccessful one $\xkun=\xzerok$ and $\dkun=\gamma^{-1}\dk$, so that
\begin{equation}\label{fourpt14}
\SSaddnew{\phikplusone-\phik\le-b_1\dk^2}
\end{equation}
Recall $\mathscr{V}_k=\accolade{
\SSaddnew{\normii{\QkRandomTranspose\nabla f(\Xzerok)}}\geq \zeta\Dk}$ from Lemma~\ref{lemmafourp7}, with $\zeta$ as in~\eqref{fourpt8}, and set, 
\begin{equation}\label{fourpt17}
\SSaddnew{b_{2,k}}:=-\nu c_1\normii{\QkTranspose\nabla f(\xzerok)}\dk+(1-\nu)(\gamma^2-1)\dk^2,\qquad
\SSaddnew{b_{3,k}}:=\nu c_3\normii{\QkTranspose\nabla f(\xzerok)}\dk+(1-\nu)(\gamma^2-1)\dk^2 ;
\end{equation}
on $\mathscr{V}_k$, $\SSaddnew{b_{2,k}}\le(-\nu c_1\zeta+(1-\nu)(\gamma^2-1))\dk^2<0$ by~\eqref{fourpt9}, while $\SSaddnew{b_{3,k}}>0$, and the random variables $\SSaddnew{\brandom_{2,k}},\SSaddnew{\brandom_{3,k}}$ are $\fqalgebra$-measurable by Lemma~\ref{lemmafourp7}. Two facts serve both algorithms: on $\mathscr{V}_k$, if the model is $(\kmge,\Qk)$-\SSaddnew{subspace-accurate} and the trial estimates are $\kfde$-weakly accurate, then the iteration is successful and $\phikplusone-\phik\le \SSaddnew{b_{2,k}}$ (Lemmas~\ref{lemma4p6} and~\ref{lemma4p7} with~\eqref{fourpt13}; the computation is Subcase~1(i) in Appendix~\ref{app:proofConvF}); and on $\mathscr{V}_k$, whatever the model and estimates, Lemma~\ref{lemmaIncr} with~\eqref{fourpt13} gives $\phikplusone-\phik\le \SSaddnew{b_{3,k}}$. Both proofs establish
\begin{align}
\E{\ivk(\Phikplusone-\Phik)|\fqalgebra}&\leq -2\ivk(1-\nu)(\gamma^2-1)\Dk^2\leq -\tfrac{1}{2}\ivk\SSaddnew{\pi_{\min}}(1-\nu)\bigl(1-\gamma^{-2}\bigr)\Dk^2,\label{fourpt15}\\
\E{\ivkc(\Phikplusone-\Phik)|\fqalgebra}&\leq -\tfrac{1}{2}\ivkc\SSaddnew{\pi_{\min}}(1-\nu)\bigl(1-\gamma^{-2}\bigr)\Dk^2,\label{fourpt16}
\end{align}
which add up to~\eqref{fourpt11}. 
\SSaddnew{Since each step has length at most $\Qmax\deltmax$, 
$\normii{\Xzerok-\x^0_0}\le k\Qmax\deltmax$, and the Lipschitz continuity of $\nabla f$ gives $$f(\Xzerok)-f_{\min}\le f(\x^0_0)-f_{\min}+k\Qmax\deltmax\normii{\nabla f(\x^0_0)}+\frac12\kappa_{Lg}(k\Qmax\deltmax)^2,$$
so $\Phi_k$ is bounded by a deterministic constant for each fixed $k$ and $\E{\Phi_{k_0}}<\infty$. Taking expectations in~\eqref{fourpt11} and summing over $k_0\le k<\SSaddnew{\bar k}$ gives $$\varrho\sum_{k=k_0}^{\bar k-1}\E{\Dk^2}\le\E{\Phi_{k_0}}-\E{\Phi_\SSaddnew{\bar k}}\le\E{\Phi_{k_0}};$$ letting $\SSaddnew{\bar k}\to\infty$ gives~\eqref{fourpttwelve}, and~\eqref{fourpt12} follows.}
}

\SSadd{\begin{proof}[Proof of Theorem~\ref{convDeltaTheorem} for \ASTRAR]
Throughout, $k\geq k_0$ and all conditioning is on $\fqalgebra$; the $\sigma$-algebra $\fmqalgebra$ is not used. Recall $$\mathscr{V}_k=\accolade{\normii{\QkRandomTranspose\nabla f(\Xzerok)}\geq\zeta\Dk}$$ from Lemma~\ref{lemmafourp7}, and recall~\eqref{fourpt13} and~\eqref{fourpt14}, which are statements about the radius update alone and are therefore unchanged.

\textbf{Case~1: $\ivk=1$.} On $\uk$ the model is $(\kmge,\Qk)$-\SSaddnew{subspace-accurate} and the trial estimates are accurate, by Corollary~\ref{cor:UniformEvent}; 
\SSaddnew{as recorded in the preliminaries, the iteration is successful and} hence
\begin{equation}\label{eq:reuse17}
\ivk\iuk\parenth{\Phikplusone-\Phik}\leq \ivk\iuk\SSaddnew{\brandom_{2,k}},\qquad \SSaddnew{\brandom_{2,k}}<0 .
\end{equation}
On $\uk^c$ nothing is known about the model or the estimates, but Lemma~\ref{lemmaIncr} bounds the increase in~$f$ produced by an accepted step without reference to either, 
SSaddnew{so the second fact of the preliminaries applies} and
\begin{equation}\label{eq:reuse19}
\ivk\iukc\parenth{\Phikplusone-\Phik}\leq \ivk\iukc\SSaddnew{\brandom_{3,k}},\qquad \SSaddnew{\brandom_{3,k}}>0 .
\end{equation}
Subcases~1(ii) and~1(iii) of the original proof are thereby covered at once, since both are instances of $\uk^c$. Because $\ivk$, $
\SSaddnew{\normii{\QkRandomTranspose\nabla f(\Xzerok)}}$, $\SSaddnew{\brandom_{2,k}}$ and $\SSaddnew{\brandom_{3,k}}$ are $\fqalgebra$-measurable while $\prob{\uk|\fqalgebra}\geq 1-\SSaddnew{\varpi_k}$ by Theorem~\ref{Theorem:Uniform}, taking conditional expectations in~\eqref{eq:reuse17}--\eqref{eq:reuse19} gives
{\bl
\begin{equation*}
\begin{split}
\E{\ivk\parenth{\Phikplusone-\Phik}|\fqalgebra}&\leq \ivk\left[(1-\SSaddnew{\varpi_k})\SSaddnew{\brandom_{2,k}}+\SSaddnew{\varpi_k\brandom_{3,k}}\right]\\
&=
\SSaddnew{-\ivk\nu\normii{\QkRandomTranspose\nabla f(\Xzerok)}\Dk}\left[(1-\SSaddnew{\varpi_k})c_1-\SSaddnew{\varpi_k}c_3\right]+\ivk(1-\nu)(\gamma^2-1)\Dk^2.
\end{split}
\end{equation*}
}

The first entry of~\eqref{eq:zetastar} gives $(1-\SSaddnew{\varpi_k})c_1-\SSaddnew{\varpi_k}c_3\geq c_1/2$, and $
\SSaddnew{\normii{\QkRandomTranspose\nabla f(\Xzerok)}}\geq\zeta\Dk$ on $\mathscr{V}_k$, so the right-hand side is at most $\ivk\left[-\frac{1}{2}\nu c_1\zeta+(1-\nu)(\gamma^2-1)\right]\Dk^2$, which by~\eqref{fourpt9} is at most $-2\ivk(1-\nu)(\gamma^2-1)\Dk^2$. This is~\eqref{fourpt15}. \SSaddnew{Since $(1-\nu)(\gamma^2-1)\Dk^2\le\frac18\nu c_1\|\QkRandomTranspose\nabla f(\Xzerok)\|\Dk$ on $\mathscr V_k$ by~\eqref{fourpt9}, the right-hand side is also at most $-\frac14\nu c_1\ivk\|\QkRandomTranspose\nabla f(\Xzerok)\|\Dk$, the form used in the proof of Lemma~\ref{l1deltakconv}.}

\textbf{Case~2: $\ivkc=1$.} No auxiliary event is required. On $\uk$, an unsuccessful iteration satisfies~\eqref{fourpt14}; a successful one satisfies $\mu\normii{\gkhat}\geq\Dk$ by Step~\textbf{[4]}, so Lemma~\ref{lemma4p8} applies and gives a decrease in~$f$ of at least $c_2\Dk^2$, whence $$\Phikplusone-\Phik\leq\parenth{-\nu c_2+(1-\nu)(\gamma^2-1)}\Dk^2.$$ Both are at most $-b_1\Dk^2$. 
On $\uk^c$, Lemma~\ref{lemmaIncr} together with $
\SSaddnew{\normii{\QkRandomTranspose\nabla f(\Xzerok)}}<\zeta\Dk$ bounds the increase by $c_3\zeta\Dk^2$, 
so $\Phikplusone-\Phik\leq b_4\Dk^2$. 
Taking conditional expectations,
\[
\E{\ivkc\parenth{\Phikplusone-\Phik}|\fqalgebra}\leq \ivkc\left[-b_1(1-\SSaddnew{\varpi_k})+b_4\SSaddnew{\varpi_k}\right]\Dk^2\leq -\frac{1}{2}\ivkc b_1\Dk^2,
\]
by the second entry of~\eqref{eq:zetastar}. Since \SSaddnew{$\pi_{\min}\le 1$} 
this is~\eqref{fourpt16} with the stated~$\varrho$.

Adding the two cases gives~\eqref{fourpt11}, and~\eqref{fourpttwelve}--\eqref{fourpt12} follow \SSaddnew{as noted in the preliminaries.}
\end{proof}}

\subsection{Liminf-type convergence}\label{sec:liminf}
{\bl The convergence result of Theorem~\ref{convDeltaTheorem} does not by itself imply first-order stationarity, since $\normii{\QkRandom^\top\nabla f(\Xzerok)}$ may be small even when $\normii{\nabla f(\Xzerok)}$ is bounded away from zero. Assumption~\ref{wellAssump} rules out this behavior with sufficiently high conditional probability by requiring the subspaces to be probabilistically well aligned. Combined with model subspace-accuracy and weak accuracy of the estimates, well alignment ensures that, whenever $\normii{\nabla f(\Xzerok)}$ remains bounded away from zero and $\Dk$ is sufficiently small, iterations on which all three properties hold force the trust-region radius to increase. Inspired by~\cite[Theorem~4.9]{DzaWildSub2022}, Theorem~\ref{theorem32} uses this mechanism to establish $\liminf_{k\to\infty}\normii{\nabla f(\Xzerok)}=0$ almost surely for both \ASTRAF\ and \ASTRAR.}
\begin{assumption}\label{wellAssump}{\bl (Inspired by~\cite[Assumption~2]{DzaWildSub2022}).}
The sequence $\accolade{\QkRandom}_{k\in\N}$ of matrices used by \SSadd{Algorithm~\ref{algoF} or~\ref{algoR}} is $(1-\beta_Q)$-probabilistically $\alpha_Q$-well aligned for some fixed $\alpha_Q\in (0,1)$ and $\beta_Q\in (0,1/2)$.
\end{assumption}
 
\begin{theorem}\label{theorem32}
Let all the assumptions made in Theorem~\ref{convDeltaTheorem} hold. Assume further that Assumption~\ref{wellAssump} holds with \SSaddnew{$\beta_Q<1-\hat\gamma_r/\pi_{\min}^2$ for \ASTRAF\ and $\beta_Q<1-\hat\gamma_r/\pi_{\min}$ for \ASTRAR.}
Then the sequence $\accoladekinN{\Xzerok}$  of random iterates generated by \SSadd{Algorithm~\ref{algoF} or~\ref{algoR}} satisfies  
\begin{equation*}
\underset{k\to\infty}{\liminf}\normii{\nabla f(\Xzerok)}=0\quad\text{almost surely.}
\end{equation*}
\end{theorem}

\SSadd{The proof is given in Appendix~\ref{app:proofsSec5}.}

\SSaddnew{The next two remarks check that the parameter conditions of Theorems~\ref{convDeltaTheorem} and~\ref{theorem32} are compatible, and record what changes for \ASTRAR.}
\begin{remark}\label{Remark:Parameters1} 
The parameter conditions in Theorems~\ref{convDeltaTheorem} and~\ref{theorem32} are mutually consistent. Indeed, for $r>\frac{1}{1-\epsLanbda}$ fixed in Theorem~\ref{convDeltaTheorem}, with $\epsLanbda\in (0,1)$, let $\hat{\gamma}_r:=\max\accolade{\frac{\gamma^{2r}}{1+\gamma^{2r}},\frac{1}{2}}$. The first condition $\frac{\SSaddnew{\pi_{\min}^2}-\hat{\gamma}_r }{2(1-\SSaddnew{\pi_{\min}})}\geq \frac{c_3}{c_1}$ in~\eqref{fourpt10}, with $c_1, c_3>0$, implies $\SSaddnew{\pi_{\min}^2}>\hat{\gamma}_r$. Hence, the upper bound $1-\frac{\hat{\gamma}_r}{\SSaddnew{\pi_{\min}^2}}>0$ imposed on~$\beta_Q$ in Theorem~\ref{theorem32} is nonvacuous. Moreover, any $0\leq \beta_Q<1-\frac{\hat{\gamma}_r}{\SSaddnew{\pi_{\min}^2}}$ gives $\pi_\star:=\SSaddnew{\pi_{\min}^2}(1-\beta_Q)>\hat{\gamma}_r\geq \frac{1}{2}$, which is the condition needed for the positive drift underlying the submartingale argument for $\accolade{\widetilde W_k}$ in the proof of Theorem~\ref{theorem32}. In addition, $\pi_\star>\frac{\gamma^{2r}}{1+\gamma^{2r}}$ is equivalent to $c_r:=\pi_\star\gamma^{-2r}+(1-\pi_\star)\gamma^{2r}<1$ (see Assumption~\ref{Assumption:cr}). Thus, the same choices of~$r, \SSaddnew{\pi_{\min}}$ and~$\beta_Q$ simultaneously ensure the probabilistic conditions required in both the convergence and complexity analyses. \SSaddnew{For \ASTRAR, $\pi_{\min}$ is free (Theorem~\ref{convDeltaTheorem}) and $\pi_{\min}^2$ is replaced by $\pi_{\min}$ (Remark~\ref{rem:theorem32R}); any $\pi_{\min}\in(\hat\gamma_r,1)$ and $0\le\beta_Q<1-\hat\gamma_r/\pi_{\min}$ give $\pi_\star:=\pi_{\min}(1-\beta_Q)>\hat\gamma_r$, and the same conclusions follow. 
}
\end{remark}

\SSaddnew{\begin{remark}\label{rem:theorem32R}
The proof of Theorem~\ref{theorem32} and Remark~\ref{Remark:AssumptVerif} use the good event $\mathscr{D}_k:=\AkEvent\cap\ik\cap\jk$ for both algorithms; for \ASTRAR, $\jk$ is the event of Definition~\ref{def:weakAcc} for the estimates~\eqref{eq:FbarsReuse} and $\Fbarokcur=\Fbar^0_k(\bm0,\nkrandom;\Theta_k)$. By Corollary~\ref{cor:UniformEvent} and the proof of Theorem~\ref{mgeTheorem}, $\uk\subseteq\ik\cap\jk$. Enlarging $k_0$ if necessary, $\prob{\uk|\fqalgebra}\ge\pi_{\min}$ for $k\ge k_0$, and since $\AkEvent$ is $\fqalgebra$-measurable, $$\prob{\mathscr{D}_k|\falgebra}\ge\prob{\AkEvent\cap\uk|\falgebra}\ge\pi_{\min}(1-\beta_Q).$$ Hence $\pi_{\min}^2$ is replaced by $\pi_{\min}$ throughout: the admissible range in Theorem~\ref{theorem32} becomes $\beta_Q<1-\hat\gamma_r/\pi_{\min}$, and $\pi_\star:=\pi_{\min}(1-\beta_Q)$ in Remark~\ref{Remark:Parameters1} and Assumption~\ref{Assumption:cr}. Thus, at the same accuracy level $\pi_{\min}$, the lower bound on the probability of a good iteration is $\pi_{\min}(1-\beta_Q)$ for \ASTRAR\ against $\pi_{\min}^2(1-\beta_Q)$ for \ASTRAF, because one accuracy event is required instead of two. 
At the same $\pi_{\min}$ and $\beta_Q$, the larger $\pi_\star$ yields a smaller $c_r$ and a larger $2\pi_\star-1$ (two places where $\pi_\star$ enters the constants $\kappa_T$, $\kappa_\Delta$ and $R_\Delta$ of Section~\ref{sec:complexity}; their other inputs, $k_0$ and $\E{S_\Delta}$, are not comparable across the two variants), and it admits 
a larger sketch failure probability, $\beta_Q<1-\hat\gamma_r/\pi_{\min}$ instead of $\beta_Q<1-\hat\gamma_r/\pi_{\min}^2$. 
\end{remark}}

\subsection{Lim-type convergence}
{\bl Theorem~\ref{theorem32} guarantees only that the gradient norms become arbitrarily small along a subsequence. To strengthen this $\liminf$-type result to convergence of the entire sequence, Lemma~\ref{l1deltakconv} shows that, for any fixed $\varepsilon>0$, the trust-region radii are summable over the iterations for which $\normii{\nabla f(\Xzerok)}>\varepsilon$. Since $\normii{\X_{k+1}^0-\Xzerok}\leq\Qmax\Dk$ and $\nabla f$ is Lipschitz continuous, every passage of the gradient norm from below $\varepsilon$ to above $2\varepsilon$ requires a fixed positive amount of cumulative trust-region radius. Lemma~\ref{l1deltakconv} therefore rules out infinitely many such passages, and Theorem~\ref{theor:limtypeconv} combines this fact with Theorem~\ref{theorem32} to obtain $\normii{\nabla f(\Xzerok)}\to 0$ almost surely.}

\begin{lemma}\label{l1deltakconv}
Let all the assumptions made in Theorem~\ref{theorem32} hold. Let $\accolade{\Xzerok}$ and $\accolade{\Dk}$ be the sequence of random iterates and trust-region radii generated by \SSadd{Algorithm~\ref{algoF} or~\ref{algoR}}. Fix $\SSaddnew{\varepsilon}>0$ and consider the random set $\mathscr{K}_\SSaddnew{\varepsilon}$ of integers~$k$ for which $\normii{\nabla f(\Xzerok)}>\SSaddnew{\varepsilon}$, that is, $\mathscr{K}_\SSaddnew{\varepsilon}:=\accolade{k\in \N:\normii{\nabla f(\Xzerok)}>\SSaddnew{\varepsilon}}$. Then 
\begin{equation*}
\sum_{k\in \mathscr{K}_\SSaddnew{\varepsilon}}\Dk <\infty\quad\text{almost surely.}
\end{equation*}
\end{lemma}

{\bl The proof presented in Appendix~\ref{app:proofsSec5} is inspired by that of~\cite[Lemma~4.17]{chen2018stochastic}, but corrects a gap in its summation argument. In that proof, 
a conditional decrease estimate of the form $\E{\Phi_{k+1}-\Phi_k|\mathscr{F}_{k-1}^{M\cdot F}}\leq-c\SSaddnew{\varepsilon}\Dk$ is  summed over a subset \SSaddnew{of the random index set $\mathscr K_\varepsilon$, whose }
indices need not be consecutive, so 
\SSaddnew{the sum of increments $\Phi_{k+1}-\Phi_k$} 
over that \SSaddnew{sub}set does not telescope. Consequently, the lower boundedness of~$\Phik$ alone does not imply the claimed summability of the corresponding radii. 
\SSaddnew{We instead introduce the $\falgebra$-measurable event $\mathscr C_k$, that $\normii{\nabla f(\Xzerok)}>\varepsilon$ while $\Dk$ is sufficiently small, prove~\eqref{fortyseven} for every deterministic $k$ and sum over all $k$, where the increments telescope and $\sum_k\Dk^2<\infty$ by Theorem~\ref{convDeltaTheorem}; since every index of $\mathscr K_\varepsilon$ eventually satisfies $\mathscr C_k$, the summability over $\mathscr K_\varepsilon$ follows.}

}
 
\begin{theorem}\label{theor:limtypeconv}
Let all the assumptions made in Theorem~\ref{theorem32} hold. Then the sequence $\accoladekinN{\Xzerok}$  of random iterates generated by \SSadd{Algorithm~\ref{algoF} or~\ref{algoR}} satisfies  
\begin{equation*}
\underset{k\to\infty}{\lim}\normii{\nabla f(\Xzerok)}=0\quad\text{almost surely.}
\end{equation*}
\end{theorem}
 
\begin{proof}
As in~\cite{chen2018stochastic} which used the following arguments, the proof is almost identical to that of the same result of~\cite[Theorem~4.3]{BaScVi2014}. Recalling the set $\mathscr{K}_\SSaddnew{\varepsilon}$ of Lemma~\ref{l1deltakconv}, the key idea is to show that if the result does not hold, then with positive probability, 
\begin{equation}\label{infDk}
\sum_{k\in \mathscr{K}_\SSaddnew{\varepsilon}}\Dk=\infty. 
\end{equation}
\SSaddnew{Since $\normii{\X^0_{k+1}-\Xzerok}\le\Qmax\Dk$, Lipschitz continuity of $\nabla f$ gives $\bigl|\normii{\nabla f(\X^0_{k+1})}-\normii{\nabla f(\Xzerok)}\bigr|\le\kappa_{Lg}\Qmax\Dk$. On a path with $\liminf_k\normii{\nabla f(\Xzerok)}=0$ and $\limsup_k\normii{\nabla f(\Xzerok)}>2\varepsilon$ there are infinitely many disjoint passages from above $2\varepsilon$ to below $\varepsilon$; each spans indices in $\mathscr K_\varepsilon$ whose radii sum to at least $\varepsilon/(\kappa_{Lg}\Qmax)$, so $\sum_{k\in\mathscr K_\varepsilon}\Dk=\infty$ on that path, contradicting Lemma~\ref{l1deltakconv}.}
\end{proof}

\section{Complexity analysis}\label{sec:complexity}
\SSaddnew{The complexity results {\bl derived below} follow from Section~\ref{sec:SupermartFram} once Assumptions~\ref{Assumption:general1} and~\ref{Assumption:cr} are verified, and the verification is common to Algorithms~\ref{algoF} and~\ref{algoR}: it uses the accuracy conclusion of Section~\ref{sec:pointwiseAcc} or~\ref{sec:uniformAcc} and Theorem~\ref{convDeltaTheorem}, which hold for both. The iteration bounds (Theorems~\ref{theor:IterComplex} and~\ref{theor:ExpIterComplex}) are therefore identical for the two variants: $\mathcal O(\varepsilon^{-2})$ almost surely along deterministic $\varepsilon_\ell\downarrow0$, with a constant $R_\Delta$ of finite expectation, and $\mathcal O(\varepsilon^{-2})$ in expectation. The sample complexity bounds (Theorems~\ref{Theorem:WorkComplexity} and~\ref{Theorem:ExpectedWork}) differ only through the oracle count {\bl $\Wk$} of~\eqref{eq:Weps}, $(q+2)\nkrandom$ for \ASTRAR\ and $(q+1)\nkrandom+2\nskrandom$ for \ASTRAF, since both sample sizes obey the same tail, mean and eventual bounds (Lemma~\ref{lemma:Nk}): they are $\mathcal O(\varepsilon^{-4}(\log(1/\varepsilon))^2)$ almost surely, $\mathcal O(\varepsilon^{-4}(1+\log(1/\varepsilon))^{1+\epsLanbda})$ in expectation under Assumption~\ref{Assumption:moments}, and $\mathcal O(\varepsilon^{-4})$ in expectation with a constant inflation (Corollary~\ref{cor:constantInflation}).}

Given $\varepsilon>0$, we consider the family $\accolade{\Tepsilon}_{\varepsilon>0}$ of random variables $\Tepsilon$ defined by
\begin{equation}\label{stoppingTime}
\taurandom_\varepsilon:=\inf\accolade{k\in\N:\normii{\nabla f(\Xzerok)}\leq\varepsilon}.
\end{equation}
Then 
$\Tepsilon$ {\bl is a possibly infinite} stopping time with respect to the filtration $\accoladekinN{\falgebra}$, {\bl whose almost-sure finiteness already guaranteed by Theorem~\ref{theorem32}, will also follow from \SSaddnew{Theorem~\ref{theor:GenIterCompl}} once Assumption~\ref{Assumption:general1} is verified. $\Tepsilon$} trivially satisfies~\eqref{eq:MonotoneTimes} so that the family $\accolade{\Tepsilon}_{\varepsilon>0}$ is monotone. Moreover, $\Dk$ is $\falgebra$-measurable for all~$k\in \N$, hence the $\sigma$-algebras $\gkalgebra$ of Section~\ref{sec:SupermartFram} can be chosen as $\gkalgebra=\falgebra$ \SSaddnew{and $\Delta_0=\delta_0$ is an input of Algorithm~\ref{algoAdapTR}, hence deterministic, as {\bl \SSaddnew{Theorem~\ref{theor:GenIterCompl}}} requires.}
 
\subsection{Iteration complexity}

\SSaddnew{Remark~\ref{Remark:AssumptVerif} verifies Assumption~\ref{Assumption:general1} except its forcing condition~\eqref{eq5:AssumpGen}, which Lemma~\ref{VilleLemmaResult} supplies; Assumption~\ref{Assumption:cr} holds by Remark~\ref{Remark:Parameters1}.}

\SSaddnew{\begin{remark}\label{Remark:AssumptVerif}
We verify Assumption~\ref{Assumption:general1} with $\gkalgebra=\falgebra$, $\Stepsize_k=\Dk$, $\alpha_0=\delta_0$, $\dmax=\deltmax$ and $\Dscr_k:=\AkEvent\cap\ik\cap\jk$ for both algorithms (Remark~\ref{rem:theorem32R}); accordingly we write $S_\Delta,R_\Delta,C_\Delta,\kappa_\Delta,\kappa_{\Delta,\lambda}$ for $S_\Stepsize,R_\Stepsize,C_\Stepsize,\kappa_\Stepsize,\kappa_{\Stepsize,\lambda}$. $\Dk$ is $\falgebra$-measurable, $\delta_0$ is an input of Algorithm~\ref{algoAdapTR}, and $\Dscr_k\in\mathscr{F}_k$ because $\AkEvent$, $\ik$ and $\jk$ are determined by $\Xzerok$, $\QkRandom$, $\Mhatk$, $\Dk$ and the estimates $\Fbarokcur,\Fbarsk$. Condition~(ii) holds since $0<\Dk\leq\deltmax$ and $\Delta_{k+1}\geq\gamma^{-1}\Dk$ at each iteration of Algorithm~\ref{algoF} or~\ref{algoR}; (iv) is~\eqref{fourpttwelve}; (v) is immediate from~\eqref{stoppingTime}. {\bl Regarding}~(i), $$\prob{\Dscr_k|\falgebra}\geq\pi_{\min}^2(1-\beta_Q)=:\pi_\star$$ for \ASTRAF\ by the tower property and~\eqref{tower1} (see also~\cite[(3.4)]{DzaWildSub2022}), and $$\prob{\Dscr_k|\falgebra}\geq\pi_{\min}(1-\beta_Q)=:\pi_\star$$ for \ASTRAR\ by Remark~\ref{rem:theorem32R}; (iii) holds with $\varepsilon_0:=\min\{1,\gamma^{-k_0}\delta_0/\kappa_a\}$ by Lemma \ref{VilleLemmaResult} below, since for $\varepsilon\le\varepsilon_0$ one has $\kappa_a\varepsilon\le\gamma^{-k_0}\delta_0<\deltmax/\gamma$ ($\delta_0<\deltmax$, $k_0\ge1$), which is also the restriction Assumption \ref{Assumption:general1}(iii) imposes.
$\pi_\star>\frac12$ and Assumption~\ref{Assumption:cr} hold by Remark~\ref{Remark:Parameters1}. 
\end{remark}}
 
\begin{lemma}\label{VilleLemmaResult}
Let all the assumptions made in Theorem~\ref{theorem32} hold. Choose $\kappa_a>0$ such that $0<\kappa_a\leq \frac{\alpha_Q}{2\SSaddnew{\kappa_{\rm suc}}+\kmge}$. 
Then for every \SSaddnew{$\varepsilon\in(0,\varepsilon_0]$ and} $k\geq k_0$,~\eqref{eq5:AssumpGen} holds.

\end{lemma}
 
\begin{proof}
\SSaddnew{Fix $\varepsilon\in(0,\varepsilon_0]$. Then }
$\kappa_a\varepsilon<\frac{\delta_{\max}}{\gamma}$. Assume that $\mathscr{E}_{\kappa_a,\varepsilon}:=\accolade{k_0\leq k<\taurandom_{\varepsilon/2},\Dk<\kappa_a\varepsilon}\cap\mathscr{D}_k$ holds, where $\mathscr{D}_k:=\AkEvent\cap\ik\cap\jk$. 
This means that $\normii{\nabla f(\Xzerok)}>\frac{\varepsilon}{2}$ on $\accolade{k<\taurandom_{\varepsilon/2}}$, $\normii{\QkRandomTranspose\nabla f(\Xzerok)}\geq \alpha_Q\normii{\nabla f(\Xzerok)}>\alpha_Q\frac{\varepsilon}{2}$ on $\accolade{k<\taurandom_{\varepsilon/2}}\cap\AkEvent$, and $\normii{\gkhat}\geq \normii{\QkRandomTranspose\nabla f(\Xzerok)}-\frac{1}{2}\kmge\Dk>\alpha_Q\frac{\varepsilon}{2}-\frac{1}{2}\kmge\Dk$ on $\accolade{k<\taurandom_{\varepsilon/2}}\cap\AkEvent\cap\ik$, hence because $\Dk<\kappa_a\varepsilon$,
\begin{equation*}
\normii{\gkhat}>\left(\frac{\alpha_Q}{2\kappa_a}-\frac{1}{2}\kmge\right)\Dk\geq \SSaddnew{\kappa_{\rm suc}}\Dk \quad\text{or equivalently}\quad \Dk<\min\accolade{\mu,\frac{1}{\SSaddnew{\kappah}},\frac{\kfcd(1-\eta)}{2(\kfde+\kmfe)}}\normii{\gkhat}.
\end{equation*}
Since $\jk$ occurs, the estimates are good ($\kfde$-weakly accurate), therefore Lemma~\ref{lemma4p7} guarantees that iteration~$k$ is successful, hence $\Delta_{k+1}=\gamma\Dk$ given that $\Dk<\kappa_a\varepsilon<\deltmax/\gamma$ (meaning that the trust-region radius is not capped), which completes the proof.
\end{proof}
 
\begin{theorem}[\SSaddnew{{\bl almost sure} iteration complexity}]\label{theor:IterComplex}
Let all the assumptions made in Lemma~\ref{VilleLemmaResult}
hold and consider the almost surely finite random variable $R_\Delta:=1+\left(1+\frac{1}{2\SSaddnew{\pi_\star}-1}\right)\frac{S_\Delta}{\kappa_a^2}$, with $\E{R_\Delta}<\infty$.
For every deterministic sequence of $\vepsell\to 0$, 
\begin{equation}\label{complRes1}
\prob{\taurandom_{\vepsell}\vepsell^2>R_\Delta\ \ i.o.}=0.
\end{equation}
In other words, along every deterministic threshold sequence of $\vepsell\downarrow 0$, $\taurandom_{\vepsell}=\mathcal{O}(\vepsell^{-2})$ with an integrable random constant.
\end{theorem}
 
\begin{proof}
The theorem's assumptions guarantee in particular that Assumption~\ref{Assumption:general1} holds (see Remarks~\ref{Remark:AssumptVerif} and~\ref{Remark:Parameters1}, and the paragraph below~\eqref{stoppingTime}), hence the result follows from Theorem~\ref{theor:IterComplex:Gen}, which completes the proof.
\end{proof}

\SSaddnew{\begin{theorem}[expected iteration complexity]\label{theor:ExpIterComplex}
Let all the assumptions made in Lemma~\ref{VilleLemmaResult} hold. Then 
with $\varepsilon_0$ as in Remark~\ref{Remark:AssumptVerif} and $\kappa_T$ as in Theorem~\ref{theor:GenIterCompl}, where $\E{S_\Delta}\le k_0\deltmax^2+\E{\Phi_{k_0}}/\varrho$ by \eqref{fourpttwelve}, 
$\E{\taurandom_\varepsilon}\le\kappa_T\varepsilon^{-2}$ for all $\varepsilon\in(0,\varepsilon_0]$ 
for both \ASTRAF\ and \ASTRAR.
\end{theorem}
\begin{proof}
Assumption~\ref{Assumption:general1} holds by Remarks~\ref{Remark:AssumptVerif} and~\ref{Remark:Parameters1} and the paragraph below~\eqref{stoppingTime}, and $\Delta_0=\delta_0$ is an input of the algorithm, so {\bl \SSaddnew{Theorem~\ref{theor:GenIterCompl}}} applies.
\end{proof}}

\subsection{Work complexity}
 \SSaddnew{Let ${\bl \Wk }$ denote the number of oracle calls made at iteration~$k$, namely ${\bl \Wk }=(q+1)\nkrandom+2\nskrandom$ for \ASTRAF\ and ${\bl \Wk }=(q+2)N_k$ for \ASTRAR\ so that ${\bl \Wk }\le w\max\{N_k,N^{0,s}_k\}$ with $w:=q+3$ for \ASTRAF\ and $w:=q+2$ for \ASTRAR. Define for every $\varepsilon>0$
\begin{equation}\label{eq:Weps}
W_\varepsilon:=\sum_{k=0}^{\taurandom_\varepsilon-1}{\bl \Wk } .
\end{equation}
With equal sample sizes $\nskrandom=\nkrandom$, \ASTRAR\ saves a fraction $1/(q+3)$ of the oracle calls of \ASTRAF, one fifth at $p=1$ and one eleventh at $p=4$; in general the two sampling rules stop at different sample sizes.
The rate in $\varepsilon$ is the same for both.
}

\SSadd{
\SSaddnew{The remainder of this subsection bounds the work, almost surely and in expectation.
}
It is the only place where a moment condition beyond Assumption~\ref{Assumpt:31} is used, and it is stated here rather than globally for the reason recorded in Remark~\ref{rem:momentScope}.}

\SSaddnew{\begin{assumption}\label{Assumption:moments}
There is a constant $\kappa_{m4}<\infty$ such that $\sup_{\x\in\rd}\E{\abs{F(\x;\Theta)-f(\x)}^4}\le\kappa_{m4}$.
\end{assumption}}

\SSaddnew{\begin{remark}\label{rem:momentScope}
Assumption~\ref{Assumption:moments} is used only in Lemma~\ref{lemma:Nk}, hence only in the sample complexity bounds. Convergence and iteration complexity use no condition on the marginal distribution of $F(\x;\Theta)$ beyond the second-moment bound of Assumption~\ref{Assumpt:31}; the Bernstein condition of Assumption~\ref{Assumption:33}, which concerns paired increments rather than marginals, is required throughout. 
The distinction matters under CRN. Consider $F(\x;\Theta)=f(\x)+Z+R(\x;\Theta)$, where $Z$ has finite variance but no fourth moment and $R$ is light-tailed with increments proportional to $\normii{\y-\x}$. Then $Z$ cancels in~\eqref{eq:DIncrement}, so Assumptions~\ref{Assumpt:31} and~\ref{Assumption:33} hold although $\E{\abs{F(\x;\Theta)-f(\x)}^4}=\infty$, and a global moment condition would exclude such problems from the convergence analysis for no reason. The sample size~\eqref{nikDefEq}, however, is driven by the empirical variance of $F$, which includes $Z$; a uniform fourth moment is a convenient sufficient condition for the deviation bound~\eqref{eq:Nktail}, and we do not claim it is necessary.
\end{remark}}

\SSaddnew{
Conditionally on $\mathscr F^Q_{k-1}$, the rule~\eqref{nikDefEq} is a sequential stopping rule of the type introduced in~\cite{ChowRobbins1965}; the next lemma gives the nonasymptotic bounds we need.

\begin{lemma}[adaptive sample sizes]\label{lemma:Nk}
Let Assumptions~\ref{Assumpt:31} and~\ref{Assumption:moments} hold, and set $\Lambda_k:=\lambda_k\kappaas^{-2}\Dk^{-2}$, $\kappa_\sigma:=\max\{1,\sigma_0^2,2v_F\}$ and $\kappa_\Lambda:=\ubar\lambda(\log2)^{1+\epsLanbda}\kappaas^{-2}\deltmax^{-2}$, so that $\Lambda_k\ge\kappa_\Lambda$ for every~$k$. Then, almost surely,
\begin{equation}\label{eq:Nktail}
\prob{\nkrandom>n\,|\,\fqalgebra}\le4(q+1)\kappa_{m4}\,\Lambda_k^2\,n^{-3}\qquad\text{for every integer }n\ge\kappa_\sigma\Lambda_k,
\end{equation}
and consequently, with $\kappa_N:=\kappaas^{-2}\bigl(\kappa_\sigma+(1+8(q+1)\kappa_{m4}\kappa_\sigma^{-2})/\kappa_\Lambda\bigr)$:
\begin{enumerate}
\item[\textup{(a)}] $\E{\nkrandom|\fqalgebra}\le\kappa_N\lambda_k\Dk^{-2}$ almost surely, for every $k$;
\item[\textup{(b)}] if $\E{\sum_{k\in\N}\Dk^2}<\infty$, there is an almost surely finite $K_W$ such that $\nkrandom\le\kappa_N\lambda_k\Dk^{-2}$ for all $k\ge K_W$.
\end{enumerate}
For \ASTRAF, \eqref{eq:Nktail}, \textup{(a)} and \textup{(b)} hold for $\nskrandom$ as well, with $\fmqalgebra$ in place of $\fqalgebra$ in~\eqref{eq:Nktail} and~\textup{(a)} with $K_W$ the larger of the two indices (one for $\nkrandom$ and one for $\nskrandom$).
\end{lemma}
\begin{proof}
Fix $k$ and write $E^i_\ell:=F(\Xki;\Theta_{k,\ell})-f(\Xki)$, so that $\hat\sigma_F^2(\siirandom_k,n;\Theta_k)\le\frac1n\sum_{\ell\le n}(E^i_\ell)^2$. The criterion in~\eqref{nikDefEq} holds at $n$ if and only if $\max\{\sigma_0^2,\hat\sigma_F^2(\siirandom_k,n;\Theta_k):i\in\intbracket{0,q}\}\le n/\Lambda_k$. Let $n\ge\kappa_\sigma\Lambda_k$. If $\nkrandom>n$, the criterion fails at $n$ while $\sigma_0^2\le n/\Lambda_k$, so $\frac1n\sum_{\ell\le n}(E^i_\ell)^2>n/\Lambda_k$ for some~$i$, and then $$\frac1n\sum_{\ell\le n}(E^i_\ell)^2-\sigma_F^2(\Xki)>n/\Lambda_k-v_F\ge n/(2\Lambda_k).$$ Conditionally on $\fqalgebra$, by Lemma~\ref{Lemma:24}, the $E^i_\ell$, $\ell\ge1$, are i.i.d.\ with mean zero, variance $\sigma_F^2(\Xki)$ and $\var{(E^i_1)^2|\fqalgebra}\le\E{(E^i_1)^4|\fqalgebra}\le \kappa_{m4}$, so Chebyshev's inequality gives $$\prob{\frac1n\sum_{\ell\le n}(E^i_\ell)^2-\sigma_F^2(\Xki)>n/(2\Lambda_k)\,|\,\fqalgebra}\le4\kappa_{m4}\Lambda_k^2n^{-3},$$ and a union bound over the $q+1$ design points gives~\eqref{eq:Nktail}.

For~(a), let $n_k:=\lceil\kappa_\sigma\Lambda_k\rceil$ and use $\sum_{n\ge n_k}n^{-3}\le2n_k^{-2}\le2(\kappa_\sigma\Lambda_k)^{-2}$:
\[
\E{\nkrandom|\fqalgebra}=\sum_{n\ge0}\prob{\nkrandom>n|\fqalgebra}\le n_k+8(q+1)\kappa_{m4}\kappa_\sigma^{-2}\le\Lambda_k\Bigl(\kappa_\sigma+\frac{1+8(q+1)\kappa_{m4}\kappa_\sigma^{-2}}{\kappa_\Lambda}\Bigr)=\kappa_N\lambda_k\Dk^{-2}.
\]

For~(b), the tail bound \eqref{eq:Nktail} holds for every deterministic $n\ge\kappa_\sigma\Lambda_k$ and hence, $\Lambda_k$ being $\fqalgebra$-measurable, at $n=n_k$, giving 
$$\prob{\nkrandom>n_k|\fqalgebra}\le4(q+1)\kappa_{m4}\kappa_\sigma^{-3}\Lambda_k^{-1}=4(q+1)\kappa_{m4}\kappa_\sigma^{-3}\kappaas^2\,\Dk^2/\lambda_k.$$ Taking expectations and summing over $k$, $$\sum_k\prob{\nkrandom>n_k}\le4(q+1)\kappa_{m4}\kappa_\sigma^{-3}\kappaas^2\,\E{\sum_k\Dk^2}/(\ubar\lambda(\log2)^{1+\epsLanbda})<\infty,$$ so the Borel-Cantelli lemma gives an almost surely finite $K_W$ with $\nkrandom\le n_k$ for all $k\ge K_W$, and $n_k\le\kappa_\sigma\Lambda_k+1\le\kappa_N\lambda_k\Dk^{-2}$.

For $\nskrandom$, condition on $\fmqalgebra$, with respect to which $\Xzerok$ and $\X^s_k$ are measurable and of which $\Theta^{0,s}_k$ is independent. The same argument with the two points $\Xzerok,\X^s_k$ in place of the $q+1$ design points gives~\eqref{eq:Nktail} and~(a), and $\prob{\nskrandom>n_k}=\E{\prob{\nskrandom>n_k|\fmqalgebra}}$ gives~(b).
\end{proof}}

\begin{theorem}\label{Theorem:WorkComplexity}
Let all the assumptions made in Lemma~\ref{VilleLemmaResult} \SSadd{and Assumption~\ref{Assumption:moments}} hold. 
Consider \SSaddnew{$R_\Delta$ as defined in Theorem~\ref{theor:IterComplex}}  
and define $\Cdelta:=\kappa_a^{-2}\left(4\ceil{R_\Delta}+5\right)<\infty$. Let $\accoladekinN{\lambda_k}$ be the deterministic sequence defined for all $k\in\N$ by
\begin{equation*}
\lambda_k=\ubar{\lambda}(\log(k+2))^{1+\epsilon_\lambda},\qquad \epsLanbda\in(0,1),\ \, \SSaddnew{\ubar{\lambda}> 0}.
\end{equation*}
\SSaddnew{Let $W_\varepsilon$ be given by~\eqref{eq:Weps}.}
Then for every deterministic sequence of $\vepsell\downarrow 0$, there exists an almost surely finite random variable $R_W$ depending on $K_W, \kappa_N, C_\Delta, \SSaddnew{R_\Delta,  \underline\lambda, \epsLanbda,} q$ and $k_0$ ($K_W$ as in Lemma~\ref{lemma:Nk}(b)) 
such that 
\begin{equation*}
\prob{W_{\vepsell}>R_W\vepsell^{-4}\left(\log(1/\vepsell)\right)^2\ \, i.o.}=0.
\end{equation*}
Equivalently, $W_{\vepsell}=\mathcal{O}\left(\vepsell^{-4}\left(\log(1/\vepsell)\right)^2\right)$ eventually almost surely along every deterministic sequence of $\vepsell\downarrow 0$.
\end{theorem}

\begin{proof}
The theorem's assumptions guarantee in particular that Assumptions~\ref{Assumption:general1} and~\ref{Assumption:cr} hold (see Remarks~\ref{Remark:AssumptVerif} and~\ref{Remark:Parameters1}, and the paragraph below~\eqref{stoppingTime}). 
\SSaddnew{Lemma~\ref{lemma:Nk}(b) gives ${\bl \Wk }\le w\kappa_N\lambda_k\Dk^{-2}$ for $k\ge K_W$, which is~\eqref{eq:WkBound1} with $\underline K_W=K_W$ and $N_W=w\kappa_N$.} 
Thus, the result immediately follows from Theorem~\ref{Theorem:WorkComplexity:Gen} by noticing that $W_\varepsilon=\sum_{k=0}^{\taurandom_\varepsilon-1}\Wk$, which completes the proof.
\end{proof}

\SSaddnew{\begin{theorem}[expected work complexity]\label{Theorem:ExpectedWork}
Let the assumptions of Theorem~\ref{Theorem:WorkComplexity} hold, and suppose in addition that $\E{S_\Delta(\log(S_\Delta+3))^{1+\epsLanbda}}<\infty$ (it holds, for instance, whenever $\E{S_\Delta^{1+\delta}}<\infty$ for some $\delta>0$). Let $W_\varepsilon$ be given by~\eqref{eq:Weps}, 
and let $\varepsilon_0$ be as in Theorem~\ref{theor:ExpIterComplex}. Then for every $\varepsilon\in(0,\varepsilon_0/2]$,
\begin{equation}\label{eq:EWald}
\E{W_\varepsilon}\ \le\ w\,\kappa_N\,\E{\sum_{k=0}^{\taurandom_\varepsilon-1}\lambda_k\Dk^{-2}}\ \le\ w\,\kappa_N\,\kappa_{\Delta,\lambda}\,\varepsilon^{-4}\bigl(1+\log(1/\varepsilon)\bigr)^{1+\epsLanbda},
\end{equation}
where $\kappa_N$ is the constant of Lemma~\ref{lemma:Nk} and $\kappa_{\Delta,\lambda}$ that of Lemma~\ref{lemma:ExpInverseRadius}(b).
\end{theorem}
\begin{proof}
By Lemma~\ref{lemma:Nk}(a) and the tower property, $\E{{\bl \Wk }|\falgebra}\le w\kappa_N\lambda_k\Dk^{-2}$ for both algorithms; Assumptions~\ref{Assumption:general1} and~\ref{Assumption:cr} hold by Remarks~\ref{Remark:AssumptVerif} and~\ref{Remark:Parameters1}, so Theorem~\ref{Theorem:ExpectedWork:Gen} applies with $\kappa_W=w\kappa_N$. 
\end{proof}
\begin{corollary}[constant inflation]\label{cor:constantInflation}
Let $\lambda_k\equiv\ubar\lambda$ (with $\kappa_\Lambda:=\underline\lambda\kappaas^{-2}\deltmax^{-2}$ in Lemma~\ref{lemma:Nk}), with $4p\,e^{-\chifde\ubar\lambda}\le1-\pi_{\min}$ for \ASTRAF, and with $\ubar\lambda\ge2p/r_F$ and $e^{-c_F\ubar\lambda}\le\min\{1-\pi_{\min},\SSaddnew{\varpi_\star}\}$ for \ASTRAR. Under the remaining assumptions of Theorem~\ref{Theorem:ExpectedWork}, without the moment condition on $S_\Delta$, $\E{W_\varepsilon}\le w\,\kappa_N\,\ubar\lambda\,\kappa_\Delta\,\varepsilon^{-4}$ for every $\varepsilon\in(0,\varepsilon_0/2]$.
\end{corollary}
\begin{proof}
By Theorems~\ref{Theorem:22} and~\ref{Theorem:Uniform}, the probability bounds used in the proof of Theorem~\ref{convDeltaTheorem} hold for every $k\ge0$, so that proof and Remark~\ref{Remark:AssumptVerif} apply; Lemma~\ref{lemma:Nk}(a) gives $\E{{\bl \Wk }|\falgebra}\le w\kappa_N\ubar\lambda\Dk^{-2}$, and the second part of Theorem~\ref{Theorem:ExpectedWork:Gen} with $\kappa_W=w\kappa_N$ gives the bound.
\end{proof}
}
 
\SSaddnew{\section{Concluding remarks}\label{sec:conclusion}
We have analyzed two variants of an adaptive sampling trust-region method that builds its models in random subspaces and uses common random numbers, and obtained for both the first $\tilde{\mathcal O}(\varepsilon^{-4})$ sample complexity within a STORM-type analysis, almost surely along deterministic sequences and in expectation. The reuse variant \ASTRAR, which evaluates the trial point on the replications that built the model, is the variant one would implement and was out of reach of the existing supermartingale analyses; what brings it within reach is an accuracy statement that holds uniformly over the trust region at the adaptive sample size (Theorem~\ref{Theorem:Uniform}), obtained by chaining a time-uniform Bernstein inequality. Its single accuracy event replaces the two of \ASTRAF, so reuse comes with one sampling stage per iteration instead of two, a larger admissible sketch failure probability, and a better dependence on $\pi_\star$ ($c_r$ smaller, $2\pi_\star-1$ larger), other parameters being equal, at the price of continuous sample paths. Both variants reach the same rates because the same three properties, a per-iteration probability bound on a favorable event, square-summable radii and a forcing property, are all that the supermartingale argument of Section~\ref{sec:SupermartFram} consumes. The dependence of the constants on $d$ and $p$, and the growth of the subspace dimension that a pathwise argument admits, are treated in Part~II. Whether the logarithmic factor and the moment condition on $S_\Delta$ in Theorem~\ref{Theorem:ExpectedWork} can be removed under the logarithmic inflation is open, as is the extension of the uniform-accuracy device to other methods that reuse samples within an iteration, such as line-search and direct-search methods.}

\subsection*{Acknowledgments}
This material was based upon work supported by National Science Foundation Grant CMMI-2226347, OAC-2410949 and Office of Naval Research Grant N000142412398.

\bibliography{KJD_bibliography-2}
\bibliographystyle{abbrvnat} 
 
\newpage
\SSadd{\appendix}

\section{Proofs for Section~\ref{sec:SupermartFram}}\label{app:proofsSec2}
\SSaddnew{The framework of Section~\ref{sec:SupermartFram} replaces the renewal--reward argument of~\cite{blanchet2016convergence} by a direct count. Its statements are in the main text; the proofs follow in the order of the statements. Three of them, those of Theorem~\ref{theor:GenIterCompl}, Lemma~\ref{VilleLemmaResult:Gen} and Lemma~\ref{lemma:ExpInverseRadius}(b), rest on the same device: on the logarithmic scale, an index before $T_{\varepsilon/2}$ at which $A_k$ is small and $\Dscr_k$ occurs contributes $+1$ and every other index at least $-1$, so the number of small indices is bounded by the number of large ones, up to the fluctuation of the count of favorable events about its conditional mean, which is a bounded-increment martingale.}

\subsection{Proof of Theorem~\ref{theor:GenIterCompl}}\label{app:B-theor:GenIterCompl}
\begin{proof}[Proof of Theorem~\ref{theor:GenIterCompl}]
Fix $\varepsilon\in(0,\varepsilon_0]$, write $T:=\taurandom_{\varepsilon/2}$ and, for $\SSaddnew{\bar k}\ge k_0$, $T_\SSaddnew{\bar k}:=T\wedge \SSaddnew{\bar k}$. Let
{\bl 
\begin{align*}
\SSaddnew{\Gamma^-_{\bar k}}&:=\abs{\accolade{k_0\le k<T_\SSaddnew{\bar k}:\Stepsize_k<a_\varepsilon}},\qquad \SSaddnew{\Gamma^+_{\bar k}}:=\abs{\accolade{k_0\le k<T_\SSaddnew{\bar k}:\Stepsize_k\geq a_\varepsilon}},\quad\text{and}\\
\Gamma_\SSaddnew{\bar k}&:=\abs{\accolade{k_0\le k<T_\SSaddnew{\bar k}:\Stepsize_k<a_\varepsilon,\ \Dscr_k\text{ occurs}}},
\end{align*}
}all zero when $T_\SSaddnew{\bar k}\le k_0$. With {\bl $R_k:=\log_\gamma(\Stepsize_k/a_\varepsilon)$, Assumption~\ref{Assumption:general1}-\emph{(ii)}} gives $R_{k+1}-R_k\ge-1$ for every $k$, {\bl and~\eqref{eq:sept} gives} $R_{k+1}-R_k=1$ at the $\Gamma_\SSaddnew{\bar k}$ indices at which $\Dscr_k$ occurs. Since $\Stepsize_{T_\SSaddnew{\bar k}}\le\dmax$ and $\Stepsize_{k_0}\ge\gamma^{-k_0}\alpha_0$, telescoping over $k_0\le k<T_\SSaddnew{\bar k}$ gives
\begin{equation}\label{eq:ExpIterTelescope}
2\Gamma_\SSaddnew{\bar k}-\SSaddnew{\Gamma^-_{\bar k}}-\SSaddnew{\Gamma^+_{\bar k}}\le\kappa_L:=k_0+\log_\gamma(\dmax/\alpha_0).
\end{equation}
The indicator $C_k:=\mathds 1_{\{k<T_\SSaddnew{\bar k}\}}{\bl \mathds 1_{\{\Stepsize_k<a_\varepsilon\}}}$ is $\gkalgebra$-measurable, so {\bl by Assumption~\ref{Assumption:general1}-\emph{(i)},}
\[
\E{\Gamma_\SSaddnew{\bar k}}=\sum_{k=k_0}^{\SSaddnew{\bar k}-1}\E{C_k\,\prob{\Dscr_k|\gkalgebra}}\ge\pi_\star\,\E{\SSaddnew{\Gamma^-_{\bar k}}}.
\]
All counts are at most $\SSaddnew{\bar k}$, so taking expectations in~\eqref{eq:ExpIterTelescope} gives $(2\pi_\star-1)\E{\SSaddnew{\Gamma^-_{\bar k}}}\le\E{\SSaddnew{\Gamma^+_{\bar k}}}+\kappa_L$. Since {\bl $h$ is nondecreasing,  $\SSaddnew{\Gamma^+_{\bar k}}\leq\SSaddnew{S_{\Stepsize,h}/}h(a_\varepsilon)$} and $T_\SSaddnew{\bar k}\le k_0+\SSaddnew{\Gamma^-_{\bar k}}+\SSaddnew{\Gamma^+_{\bar k}}$. {\bl Therefore
\begin{equation*}
\E{T_{\bar k}}\leq k_0+
\frac{\kappa_L}{2\pi_\star-1}+ \frac{2\pi_\star}{2\pi_\star-1}
\frac{\E{S_{A,h}}}{h(a_\varepsilon)},
\end{equation*}\SSaddnew{and the right-hand side does not depend on $\bar k$.} 
}
{\bl Since $T_\SSaddnew{\bar k}\uparrow \taurandom_{\varepsilon/2}$ as $\SSaddnew{\bar k}\to\infty$, regardless of whether $\taurandom_{\varepsilon/2}$ is initially finite, the} monotone convergence {\bl theorem gives 
\SSaddnew{$\E{T_{\varepsilon/2}}=\lim_{\bar k\to\infty}\E{T_{\bar k}}$ is at most the right-hand side of~\eqref{eq:ExpIterGen}, which is finite; }
hence $\taurandom_{\varepsilon/2}<\infty$ almost surely. Assumption~\ref{Assumption:general1}-\emph{(v)} gives $T_{\varepsilon_0}\leq  \taurandom_\varepsilon\le\taurandom_{\varepsilon/2}$ for every $\varepsilon\in(0,\varepsilon_0]$, and $\Tepsilon\leq T_{\varepsilon_0}$ whenever $\varepsilon>\varepsilon_0,$ therefore $\Tepsilon<\infty$ almost surely for all $\varepsilon>0$, which completes the proof.
}

\SSaddnew{For the particular case, \eqref{eq:sept} is Assumption~\ref{Assumption:general1}-\emph{(iii)} with $a_\varepsilon=\kappa_a\varepsilon$, and $h(t)=t^2$ gives $h(a_\varepsilon)=\kappa_a^2\varepsilon^2$; since $\varepsilon\le\varepsilon_0\le1$, the first two terms of~\eqref{eq:ExpIterGen} are at most $(k_0+\kappa_L/(2\pi_\star-1))\varepsilon^{-2}$, which gives~\eqref{eq:ExpIter}.}
\end{proof}

\subsection{Proof of Lemma~\ref{VilleLemmaResult:Gen}}\label{app:B-VilleLemmaResult:Gen}
\SSaddnew{\begin{proof}[Proof of Lemma~\ref{VilleLemmaResult:Gen}]
Write $\bar\varepsilon_j:=2^{-j}$ and consider $j$ with $\bar\varepsilon_j\le\varepsilon_0$, for which Assumption~\ref{Assumption:general1}-\emph{(iii)} gives $\mathscr{E}_{\kappa_a,\bar\varepsilon_j}\subseteq\accolade{\Stepsize_{k+1}=\gamma\Stepsize_k}$ for every $k\ge k_0$. Define
\begin{align*}
    \SSaddnew{\Gamma^-_j}& :=\abs{\accolade{k_0\leq k<\taurandom_{\bar\varepsilon_j/2}:\Stepsize_k<\kappa_a\bar\varepsilon_j}},\\
    \SSaddnew{\Gamma^+_j} & :=\abs{\accolade{k_0\leq k<\taurandom_{\bar\varepsilon_j/2}:\Stepsize_k\geq\kappa_a\bar\varepsilon_j}},\\
    \Gamma_j&:=\abs{\accolade{k_0\leq k<\taurandom_{\bar\varepsilon_j/2}:\Stepsize_k<\kappa_a\bar\varepsilon_j\ \text{and}\ \Dscr_k\ \text{occurs}}},
\end{align*}
all finite since $\bigcap_{j\ge1}\{\taurandom_{\bar\varepsilon_j/2}<\infty\}$ is almost sure {\bl thanks to \SSaddnew{Theorem~\ref{theor:GenIterCompl}}}, and $R_k^{(j)}:=\log_\gamma\bigl(\Stepsize_k/(\kappa_a\bar\varepsilon_j)\bigr)$. 

If $\taurandom_{\bar\varepsilon_j/2}\le k_0$ the three counts are zero and~\eqref{NGammaLog} is trivial since $A_{k_0}\le\dmax$; otherwise telescoping over $k_0\le k<\taurandom_{\bar\varepsilon_j/2}$ gives
\begin{equation*}
R_{\taurandom_{\bar\varepsilon_j/2}}^{(j)}-R_{k_0}^{(j)}=\sum_{k=k_0}^{\taurandom_{\bar\varepsilon_j/2}-1}\left(R_{k+1}^{(j)}-R_k^{(j)}\right)=\log_\gamma\left(\frac{\Stepsize_{\taurandom_{\bar\varepsilon_j/2}}}{\Stepsize_{k_0}}\right)\leq \log_\gamma\left(\frac{\dmax}{\Stepsize_{k_0}}\right).
\end{equation*}
For all~$k\geq k_0$, the inequality $\Stepsize_{k+1}\geq\gamma^{-1}\Stepsize_k$ from~\eqref{eq4:AssumpGen} ensures that $R_{k+1}^{(j)}-R_k^{(j)}\geq -1$
while on $\mathscr{E}_{\kappa_a,\SSadd{\bar\varepsilon_j}}$, we have
$\Stepsize_{k+1}=\gamma\Stepsize_k$ and therefore $R_{k+1}^{(j)}-R_k^{(j)}=1$. 
This means that among the $\SSaddnew{\Gamma^-_j}$ indices for which $\Stepsize_k<\kappa_a \bar\varepsilon_j$, 
the $\Gamma_j$ indices at which $\Dscr_k$ occurs have increments $+1$,
while the remaining $\SSaddnew{\Gamma^-_j} -\Gamma_j$ corresponds to increments bounded below by~$-1$. Likewise, among the $\SSaddnew{\Gamma^+_j}$ indices for which $\Stepsize_k\geq\kappa_a \bar\varepsilon_j$, every increment is bounded below by~$-1$. Therefore, 
\begin{equation}\label{NGammaLog}
\log_\gamma\left(\frac{\dmax}{\Stepsize_{k_0}}\right)\geq R_{\taurandom_{\bar\varepsilon_j/2}}^{(j)}-R_{k_0}^{(j)}\geq \Gamma_j-(\SSaddnew{\Gamma^-_j} -\Gamma_j)-\SSaddnew{\Gamma^+_j}=2\Gamma_j-\SSaddnew{\Gamma^-_j} -\SSaddnew{\Gamma^+_j}.
\end{equation}

We now compare $\Gamma_j$ with $\SSaddnew{\Gamma^-_j}$. Let $C_k^{(j)}:=\mathds 1_{\{k<\taurandom_{\bar\varepsilon_j/2}\}}\mathds 1_{\{\Stepsize_k<\kappa_a\bar\varepsilon_j\}}$, which is $\gkalgebra$-measurable, let $\Pi_k:=\prob{\Dscr_k|\gkalgebra}$, so that $\Pi_k\ge\pi_\star$ for $k\ge k_0$ by~\eqref{eq3:AssumpGen}, and set
\[
O_j:=\sum_{k=k_0}^{\infty}C_k^{(j)}\bigl(\mathds 1_{\Dscr_k}-\Pi_k\bigr),
\]
a sum with finitely many nonzero terms. Since $\sum_{k\ge k_0}C^{(j)}_k\mathds 1_{\Dscr_k}=\Gamma_j$ and $\sum_{k\ge k_0}C^{(j)}_k\Pi_k\ge\pi_\star\SSaddnew{\Gamma^-_j}$, we have $\Gamma_j\ge\pi_\star\SSaddnew{\Gamma^-_j}+O_j$, and~\eqref{NGammaLog} gives
\begin{equation}\label{alephpm}
(2\pi_\star-1)\SSaddnew{\Gamma^-_j}\le\SSaddnew{\Gamma^+_j}+\log_\gamma\Bigl(\frac{\dmax}{\Stepsize_{k_0}}\Bigr)+2\abs{O_j}.
\end{equation}
The summands of $O_j$ are martingale differences with respect to $\{\mathscr G_k\}_{k\ge k_0}$, bounded by one, because $C^{(j)}_k$ and $\Pi_k$ are $\gkalgebra$-measurable while $\Dscr_k\in\mathscr G_k$. Hence, for every $\bar k>k_0$,
\[
\E{\Bigl(\sum_{k=k_0}^{\bar k-1}C^{(j)}_k(\mathds 1_{\Dscr_k}-\Pi_k)\Bigr)^2}=\E{\sum_{k=k_0}^{\bar k-1}C^{(j)}_k\Pi_k(1-\Pi_k)}\le\E{\taurandom_{\bar\varepsilon_j/2}}\le\kappa_T\bar\varepsilon_j^{-2}
\]
by {\bl \SSaddnew{Theorem~\ref{theor:GenIterCompl}}}, and Fatou's lemma gives $\E{O_j^2}\le\kappa_T\bar\varepsilon_j^{-2}$. For every $\delta>0$, Chebyshev's inequality gives $$\sum_{j\ge1}\prob{\bar\varepsilon_j^2\abs{O_j}>\delta}\le\delta^{-2}\kappa_T\sum_{j\ge1}4^{-j}<\infty,$$ so the Borel-Cantelli lemma, applied with $\delta=1/m$ for every $m\in\N_+$, gives $\bar\varepsilon_j^2O_j\to0$ almost surely. Multiplying~\eqref{alephpm} by $\bar\varepsilon_j^2$ and using $\SSaddnew{\Gamma^+_j}(\kappa_a\bar\varepsilon_j)^2\le S_\Stepsize$ yields, almost surely,
\begin{equation}\label{limsupjNj2}
\underset{j\to\infty}{\limsup}\ \bar\varepsilon_j^2\SSaddnew{\Gamma^-_j}\le\frac{S_\Stepsize}{(2\pi_\star-1)\kappa_a^2}.
\end{equation}
Finally, for each~$\ell$, choose $j(\ell)$ such that $\bar\varepsilon_{j(\ell)}/2<\vepsell\le\bar\varepsilon_{j(\ell)}$, so that $\taurandom_{\vepsell}\leq \taurandom_{\bar\varepsilon_{j(\ell)}/2}$ and $\Stepsize_k<\kappa_a \bar\varepsilon_{j(\ell)}$ whenever $\Stepsize_k<\kappa_a\vepsell$, then observe that
\begin{equation}\label{eqsixty}
\abs{\accolade{0\leq k<\taurandom_{\vepsell}:\Stepsize_k<\kappa_a\vepsell}}\leq k_0+\SSaddnew{\Gamma^-_{j(\ell)}}.
\end{equation} 
Multiplying each side of~\eqref{eqsixty} by $\vepsell^2$, using $\vepsell^2\leq \bar\varepsilon_{j(\ell)}^2$ and letting $\ell\to\infty$, yield 
\begin{equation*}
\underset{\ell\to\infty}{\limsup}\ \vepsell^2\abs{\accolade{0\leq k<\taurandom_{\vepsell}:\Stepsize_k<\kappa_a\vepsell}}\leq\underset{\ell\to\infty}{\limsup}\ \bar\varepsilon_{j(\ell)}^2\SSaddnew{\Gamma^-_{j(\ell)}}\leq \frac{1}{(2\pi_\star-1)\kappa_a^2}\sum_{k=0}^{\infty}\Stepsize_k^2,
\end{equation*}
where the last inequality follows from~\eqref{limsupjNj2} and the fact that $j(\ell)\to\infty$ as $\ell\to\infty$. This completes the proof.
\end{proof}}

\subsection{Proof of Theorem~\ref{theor:IterComplex:Gen}}\label{app:B-theor:IterComplex:Gen}
\begin{proof}[Proof of Theorem~\ref{theor:IterComplex:Gen}]
Fix an arbitrary deterministic sequence $\vepsell\to 0$. {\bl By \SSaddnew{Theorem~\ref{theor:GenIterCompl}} and Assumption~\ref{Assumption:general1}-\emph{(v)}}, $\taurandom_{\vepsell}<\infty$ {\bl almost surely for every~$\ell$,}
and we have
\begin{equation}\label{tauvepsel1}
\taurandom_{\vepsell} = \abs{\accolade{0\leq k<\taurandom_{\vepsell}:\Stepsize_k<\kappa_a\vepsell}}+\abs{\accolade{0\leq k<\taurandom_{\vepsell}:\Stepsize_k\geq\kappa_a\vepsell}}.
\end{equation}
Recalling~$S_\Stepsize$ of Assumption~\ref{Assumption:general1}-\emph{(iv)}, we note that $\abs{\accolade{0\leq k<\taurandom_{\vepsell}:\Stepsize_k\geq\kappa_a\vepsell}}(\kappa_a\vepsell)^2\leq S_\Stepsize$ which, combined with~\eqref{tauvepsel1} and Lemma~\ref{VilleLemmaResult:Gen}, yields, almost surely
\begin{equation*}
\underset{\ell\to\infty}{\limsup}\ \taurandom_{\vepsell}\vepsell^2\leq \frac{S_\Stepsize}{(2\SSaddnew{\pi_\star}-1)\kappa_a^2}+ \frac{S_\Stepsize}{\kappa_a^2}=:R_\Stepsize-1,
\end{equation*}
where the random variable $R_\Stepsize:=1+\left(1+\frac{1}{2\SSaddnew{\pi_\star}-1}\right)\frac{S_\Stepsize}{\kappa_a^2}<\infty$ almost surely, with $\E{R_\Stepsize}<\infty$ thanks to~\eqref{eq:SDeltaFinite}. Thus, for almost every~$\omega$, there exists $L_0(\omega)\in\N$ such that $\taurandom_{\vepsell}(\omega)\vepsell^2\leq R_\Stepsize(\omega)$ for all $\ell\geq L_0(\omega)$. Equivalently, this means that~\eqref{complRes1:Gen} holds, which completes the proof.
\end{proof}

\subsection{Proof of Lemma~\ref{lemma:br}}\label{app:B-lemma:br}
\begin{proof}
[Proof of Lemma~\ref{lemma:br}]
The proof of the first result considers two separate cases, the first of which will be broken into two subcases. \SSaddnew{Since $j$ is admissible, $\kappa_a\bar\varepsilon_j\le\kappa_a\varepsilon_0\le\gamma^{-k_0}\alpha_0\le\Stepsize_{k_0}$.
}

{\textbf{Case~1:} $k<\taurandom_{\SSadd{\bar\varepsilon_j/2}}$ and $\Stepsizek<\kappa_a \SSadd{\bar\varepsilon_j}$. 
\\
{\textit{(Subcase~$1$)}} If the event $\mathscr{D}_k$ occurs, then the event $\mathscr{E}_{\kappa_a,\SSadd{\bar\varepsilon_j}}$ 
holds and Assumption~\ref{Assumption:general1}-\emph{(iii)} ensures that $\Stepsize_{k+1}=\gamma\Stepsizek$. Therefore,
\begin{equation*}
\ZkPlusOnej=\max\accolade{1,\gamma^{-2}\left(\frac{\kappa_a \SSadd{\bar\varepsilon_j}}{\Stepsizek}\right)^2}\leq \max\accolade{1,\gamma^{-2}\Zkj},
\end{equation*}
hence, taking the $r$th power and using the inequality $\max\accolade{1,a}\leq 1+a$, yields
\begin{equation}\label{ineq:iDk}
\left(\ZkPlusOnej\right)^r\mathds{1}_{\mathscr{D}_k}\leq \left(1+\gamma^{-2r}\left(\Zkj\right)^r\right)\mathds{1}_{\mathscr{D}_k}.
\end{equation}
{\textit{(Subcase~$2$)}} If the event $\mathscr{D}_k$ does not occur, since it always holds that $\Stepsize_{k+1}\geq \gamma^{-1}\Stepsizek$, then 
\begin{equation*}
\ZkPlusOnej\leq \max\accolade{1,\gamma^{2}\left(\frac{\kappa_a \SSadd{\bar\varepsilon_j}}{\Stepsizek}\right)^2}\leq \gamma^2\Zkj,
\end{equation*}
which implies that
\begin{equation}\label{ineq:iDkc}
\left(\ZkPlusOnej\right)^r\mathds{1}_{\overline{\mathscr{D}}_k}\leq\gamma^{2r}\left(\Zkj\right)^r\mathds{1}_{\overline{\mathscr{D}}_k}.
\end{equation}
Combining~\eqref{ineq:iDk} and~\eqref{ineq:iDkc} yields
\begin{equation}
\left(\ZkPlusOnej\right)^r\leq \left(1+\gamma^{-2r}\left(\Zkj\right)^r\right)\mathds{1}_{\mathscr{D}_k}+\gamma^{2r}\left(\Zkj\right)^r\mathds{1}_{\overline{\mathscr{D}}_k}.
\end{equation}
Taking conditional expectation with respect to $\gkalgebra$ and using that $\Zkj$ is $\gkalgebra$-measurable, gives, almost surely
\begin{equation}\label{ineq:Zk}
\begin{split}
\E{\left(\ZkPlusOnej\right)^r|\gkalgebra}&\leq \Pi_k\left(1+\gamma^{-2r}\left(\Zkj\right)^r\right)+(1-\Pi_k)\gamma^{2r}\left(\Zkj\right)^r\\
&=\Pi_k+\left(\Pi_k\gamma^{-2r}+(1-\Pi_k)\gamma^{2r}\right)\left(\Zkj\right)^r,
\end{split}
\end{equation}
\SSaddnew{where $\Pi_k:=\prob{\Dscr_k\mid\gkalgebra}$ satisfies $\Pi_k\ge\pi_\star$ by Assumption~\ref{Assumption:general1}-(i).}
We note that the function $t\mapsto\gamma^{-2r}t+\gamma^{2r}(1-t)$ whose derivative equals $\gamma^{-2r}-\gamma^{2r}<0$ is decreasing on $(\pi_\star,1)$, hence $\Pi_k\gamma^{-2r}+(1-\Pi_k)\gamma^{2r}\leq \pi_\star\gamma^{-2r}+(1-\pi_\star)\gamma^{2r}=c_r$. Thus, it follows from~\eqref{ineq:Zk} that
\begin{equation}\label{ineq:Zkcase1}
\E{\left.\left(\ZkPlusOnej\right)^r\right\rvert\gkalgebra}\leq c_r\left(\Zkj\right)^r+1.
\end{equation}

{\textbf{Case~2:} $k<\taurandom_{\SSadd{\bar\varepsilon_j/2}}$ and $\Stepsizek\geq\kappa_a \SSadd{\bar\varepsilon_j}$.}\\
We note that $\Zkj=1$ \SSaddnew{on $k\geq k_0$ belonging to  Case 2, }
and the inequalities $\Stepsize_{k+1}\geq \gamma^{-1}\Stepsizek\geq \gamma^{-1}\kappa_a \SSadd{\bar\varepsilon_j}$ yield $\frac{\kappa_a \SSadd{\bar\varepsilon_j}}{\Stepsize_{k+1}}\leq \gamma$, therefore
\begin{equation*}
\left(\ZkPlusOnej\right)^r\leq \gamma^{2r}\leq c_r\left(\Zkj\right)^r+\gamma^{2r},
\end{equation*}
which implies
\begin{equation}\label{ineq:Zkcase2}
\E{\left.\left(\ZkPlusOnej\right)^r\right\rvert \SSaddnew{\gkalgebra}}\leq c_r\left(\Zkj\right)^r+\gamma^{2r}.
\end{equation}

Combining both cases' results~\eqref{ineq:Zkcase1} and~\eqref{ineq:Zkcase2} yields
\begin{equation}\label{ineq:Zkcase}
\E{\left.\left(\ZkPlusOnej\right)^r\mathds{1}_{\accolade{k<\taurandom_{\SSadd{\bar\varepsilon_j/2}}}}\right\rvert\gkalgebra}\leq 
\SSaddnew{\bigl(c_r(\Zkj)^r+b_r\bigr)\mathds 1_{\{k<\taurandom_{\bar\varepsilon_j/2}\}}}
\end{equation}
\SSaddnew{since both cases are on $\{k<\taurandom_{\bar\varepsilon_j/2}\}\in\gkalgebra$.} 
We note from $\accolade{k+1<\taurandom_{\SSadd{\bar\varepsilon_j/2}}}\subset\accolade{k<\taurandom_{\SSadd{\bar\varepsilon_j/2}}}$ that
\begin{equation*}
\left(\XkPlusOnej\right)^r=\left(\ZkPlusOnej\right)^r\mathds{1}_{\accolade{k+1<\taurandom_{\SSadd{\bar\varepsilon_j/2}}}}\leq \left(\ZkPlusOnej\right)^r\mathds{1}_{\accolade{k<\taurandom_{\SSadd{\bar\varepsilon_j/2}}}}
\end{equation*}
which, combined with~\eqref{ineq:Zkcase}, yields~\eqref{ineq:br}, that is, 
\SSaddnew{\[\E{\left.\left(\XkPlusOnej\right)^r\right\rvert\gkalgebra}\leq \bigl(c_r(\Zkj)^r+b_r\bigr)\mathds{1}_{\accolade{k<\taurandom_{\bar\varepsilon_j/2}}}=c_r\left(\Xkj\right)^r+b_r\mathds{1}_{\accolade{k<\taurandom_{\bar\varepsilon_j/2}}}.\]}
From~\eqref{ineq:br}, by defining $\kkj:=\E{\left(\Xkj\right)^r}$ for all $k\geq k_0$, we obtain the recursion $\kkPlusOnej\leq c_r\kkj+b_r,$
which after iterating from~$k_0$ to~$k$, yields
\begin{equation}\label{ineq:recursion}
\kkj\leq c_r^{k-k_0}\SSaddnew{\varkappa_{k_0}^{(j)}}+b_r\sum_{i=0}^{k-k_0-1}c_r^i.
\end{equation}
\SSaddnew{By admissibility, $Z^{(j)}_{k_0}=1$, so $\varkappa^{(j)}_{k_0}\le1$.} 
Since $c_r<1$, it follows from~\eqref{ineq:recursion} that for all $k\geq k_0$,
\begin{equation*}
\kkj\leq c_r^{k-k_0}+b_r\frac{1-c_r^{k-k_0}}{1-c_r}\leq 1+\frac{b_r}{1-c_r},
\end{equation*}
which proves~\eqref{ineq:ar}. 

In order to show~\eqref{ineq:Jensen}, we use Jensen inequality applied to the convex function $t\mapsto t^r$ to obtain
\begin{equation*}
\left(\frac{1}{l-k_0}\sum_{k=k_0}^{l-1}\Xkj\right)^r\leq \frac{1}{l-k_0}\sum_{k=k_0}^{l-1}\left(\Xkj\right)^r.
\end{equation*}
Thus, multiplying each side by $(l-k_0)^r\leq l^r$, then taking expectation and using~\eqref{ineq:ar}, yields~\eqref{ineq:Jensen},
which completes the proof.}
\end{proof}

\subsection{Proof of Lemma~\ref{lemma:InverseRadius}}\label{app:B-lemma:InverseRadius}
\begin{proof}
[Proof of Lemma~\ref{lemma:InverseRadius}]
We fix~$\nu\in\N_+$ and apply Markov inequality using Lemma~\ref{lemma:br} especially~\eqref{ineq:Jensen} with $l=l_{j,\nu}:=\ceil{\nu 2^{2(j+1)}}$ to obtain, for all~$j$ sufficiently large
\begin{equation}\label{eq:newarbound}
\prob{\sum_{k=k_0}^{l_{j,\nu}-1}\Xkj>l_{j,\nu}j^{1-\epsLanbda}}\leq \frac{{\E{\left(\sum_{k=k_0}^{l_{j,\nu}-1}\Xkj\right)^r}}}{l_{j,\nu}^rj^{(1-\epsLanbda)r}}\leq \frac{a_r}{j^{(1-\epsLanbda)r}},
\end{equation}
where $\Xkj:=\Zkj\mathds{1}_{\accolade{k<\taurandom_{\SSadd{\bar\varepsilon_j/2}}}}$ with $\Zkj:=\max\accolade{1,\left(\frac{\kappa_a \SSadd{\bar\varepsilon_j}}{\Stepsizek}\right)^2}$. We note that $\sum_{j=1}^{\infty}a_rj^{-(1-\epsLanbda)r}<\infty$ since $(1-\epsLanbda)r>1$, therefore it follows from the Borel--Cantelli first lemma that almost surely, for all sufficiently large~$j$, $\sum_{k=k_0}^{l_{j,\nu}-1}\Xkj\leq l_{j,\nu}j^{1-\epsLanbda}$. This means that for each fixed $\nu\geq 1$, there exists an almost sure event $\Omega_\nu$ such that for every $\omega\in\Omega_\nu$, there exists $J_\nu(\omega)\in\N$ such that $\sum_{k=k_0}^{l_{j,\nu}-1}\Xkj(\omega)\leq l_{j,\nu}j^{1-\epsLanbda}$ for all $j\geq J_\nu(\omega)$. In the remainder of the proof, we fix $\omega$ in 
\begin{equation*}
\tilde{\Omega}:=\bigcap_{\nu=1}^\infty\Omega_\nu\cap\Omega_{IC}\cap\accolade{R_\Stepsize<\infty},
\end{equation*}
which is almost sure as a countable intersection of almost sure events, where 
\begin{equation*}
\Omega_{IC}:=\accolade{\omega\in\Omega:\exists J_{IC}(\omega)<\infty\ \text{such that}\ \taurandom_{\SSadd{\bar\varepsilon_j/2}}(\omega)\leq R_{\Stepsize}(\omega) 2^{2(j+1)}\ \text{for all}\ j\geq J_{IC}(\omega)}
\end{equation*}
is almost sure thanks to Theorem~\ref{theor:IterComplex:Gen} and $\accolade{R_\Stepsize<\infty}$ is almost sure thanks to~\eqref{eq:SDeltaFinite}. We choose in particular $\nu_\Stepsize:=\ceil{R_{\Stepsize}(\omega)}+1$ so that there exists $J_{\nu_\Stepsize}(\omega)\in\N$ such that 
\begin{equation*}
\sum_{k=k_0}^{l_{j,\nu_\Stepsize}-1}\Xkj(\omega)\leq l_{j,\nu_\Stepsize}j^{1-\epsLanbda}\quad\text{for all}\ j\geq J_{\nu_\Stepsize}(\omega),
\end{equation*}
where $l_{j,\nu_\Stepsize}:=\ceil{\nu_\Stepsize 2^{2(j+1)}}$. Because $\omega\in \Omega_{IC}$, we have for all $j\geq J_{IC}(\omega)$ that
\begin{equation}
\taurandom_{\SSadd{\bar\varepsilon_j/2}}(\omega)\leq R_{\Stepsize}(\omega) 2^{2(j+1)}\leq \nu_\Stepsize 2^{2(j+1)}\leq l_{j,\nu_\Stepsize}\leq (4\nu_\Stepsize+1)2^{2j}.
\end{equation}
Therefore, for all $j\geq J(\omega):=\max\accolade{J_{\nu_\Stepsize}(\omega), J_{IC}(\omega)}$,
\begin{equation}
\sum_{k=k_0}^{\taurandom_{\SSadd{\bar\varepsilon_j/2}}(\omega)-1}\Zkj(\omega)=\sum_{k=k_0}^{l_{j,\nu_\Stepsize}-1}\Xkj(\omega)\leq l_{j,\nu_\Stepsize}j^{1-\epsLanbda}\leq (4\nu_\Stepsize+1)2^{2j}j^{1-\epsLanbda},
\end{equation}
which, combined with the inequality $\Stepsizek^{-2}(\omega)\leq (\kappa_a \SSadd{\bar\varepsilon_j})^{-2}\Zkj(\omega)=\kappa_a^{-2}2^{2j}\Zkj(\omega)$ resulting from the definition of~$\Zkj$, yields 
\begin{equation*}
\sum_{k=k_0}^{\taurandom_{\SSadd{\bar\varepsilon_j/2}}(\omega)-1}\Stepsizek^{-2}(\omega)\leq \kappa_a^{-2}(4\nu_\Stepsize+1)2^{4j}j^{1-\epsLanbda},
\end{equation*}
which completes the proof.
\end{proof}

\subsection{Proof of Theorem~\ref{Theorem:WorkComplexity:Gen}}\label{app:B-Theorem:WorkComplexity:Gen}
\begin{proof}[Proof of Theorem~\ref{Theorem:WorkComplexity:Gen}]
We consider the events $\accolade{R_\Stepsize <\infty}$ satisfying $\prob{R_\Stepsize <\infty}=1$ thanks to~\eqref{eq:SDeltaFinite}, 
\begin{equation*}
\tilde{\Omega}_{IC}:=\accolade{\omega\in\Omega:\exists L(\omega)\in\N\ \text{such that}\ \forall\ell\geq L(\omega), \taurandom_{\vepsell}(\omega)\leq R_\Stepsize(\omega)\vepsell^{-2}}
\end{equation*}
which is almost sure thanks to Theorem~\ref{theor:IterComplex:Gen}, 
\begin{equation*}
\Omega_W:=\accolade{\omega\in\Omega:\exists \SSaddnew{\ubarlong{K}}_W(\omega), N_W(\omega)\ \text{such that}\ \forall k\geq \ubarlong{K}_W(\omega),\Wk(\omega)\leq N_W(\omega)\lambda_k\Stepsize_k^{-2}(\omega)}
\end{equation*}
which is almost sure by hypothesis, and 
\begin{equation*}
\Omega_{1/\Stepsize}:=\accolade{\omega\in\Omega:\exists J(\omega)\in\N\ \text{such that}\ \forall j\geq J(\omega), \sum_{k=k_0}^{\taurandom_{\SSadd{\bar\varepsilon_j/2}}(\omega)-1}\Stepsize_k^{-2}(\omega)\leq \Cdelta(\omega) 2^{4j}j^{1-\epsilon_\lambda}}
\end{equation*}
which is almost sure thanks to Lemma~\ref{lemma:InverseRadius}. In the remainder of the proof, we work on the almost sure event
\begin{equation*}
\breve{\Omega}:=\tilde{\Omega}_{IC}\cap\Omega_W\cap\Omega_{1/\Stepsize}\cap\accolade{R_\Stepsize <\infty}.
\end{equation*}
We first observe that $\chi_W:=\ds{\sum_{k=0}^{K_W-1}}\Wk$, with $K_W:=k_0\vee \ubarlong{K}_W$, is finite almost surely, 
hence for all sufficiently large~$\ell$, 
\begin{equation}\label{ineq:W1}
\begin{split}
W_{\vepsell}&\leq \chi_W+N_W\sum_{k=K_W}^{\taurandom_{\vepsell}-1}\lambda_k\Stepsize_k^{-2}\leq \chi_W+N_W\sum_{k=k_0}^{\taurandom_{\vepsell}-1}\lambda_k\Stepsize_k^{-2}\leq \chi_W+N_W\lambda_{\taurandom_{\vepsell}}\sum_{k=k_0}^{\taurandom_{\vepsell}-1}\Stepsize_k^{-2},
\end{split}
\end{equation}
where the last inequality follows from $\lambda_k\leq \lambda_{\taurandom_{\vepsell}}$ for all sufficiently large~$\ell$ and all $k<\taurandom_{\vepsell}$, given that the sequence of $\lambda_k$ is eventually nondecreasing. It follows from the inequality $\taurandom_{\vepsell}\leq R_{\Stepsize}\vepsell^{-2}$ (for all sufficiently large~$\ell$) 
\SSaddnew{and $\varepsilon_\ell\le1$ that $$\log(\taurandom_{\varepsilon_\ell}+2)\le\log(R_A\varepsilon_\ell^{-2}+2)\le\log(R_A+2)+2\log(1/\varepsilon_\ell),$$}
hence since  $R_{\Stepsize}<\infty$ almost surely, there exists an almost surely finite random variable $C_R>0$ such that 
\begin{equation}\label{ineq:W2}
\lambda_{\taurandom_{\vepsell}}=\ubar{\lambda}\left(\log(\taurandom_{\vepsell}+2)\right)^{1+\epsLanbda}\leq C_R\left(\log(1/\vepsell)\right)^{1+\epsLanbda}.
\end{equation}
Now choose $j(\ell)$ such that $\SSadd{\bar\varepsilon_{j(\ell)}/2}<\vepsell\leq \SSadd{\bar\varepsilon_{j(\ell)}}$ so that $\taurandom_{\vepsell}\leq \taurandom_{\SSadd{\bar\varepsilon_{j(\ell)}/2}}$. By Lemma~\ref{lemma:InverseRadius},
\begin{equation}\label{ineq:W3}
\sum_{k=k_0}^{\taurandom_{\vepsell}-1}\Stepsize_k^{-2}\leq \sum_{k=k_0}^{\taurandom_{\SSadd{\bar\varepsilon_{j(\ell)}/2}}-1}\Stepsize_k^{-2}\leq C_\Stepsize 2^{4j(\ell)}j(\ell)^{1-\epsLanbda}\leq C_\Stepsize\vepsell^{-4}\left(1/\log 2\right)^{1-\epsLanbda}\left(\log(1/\vepsell)\right)^{1-\epsLanbda},
\end{equation}
where the last inequality follows from $2^{4j(\ell)}\leq \vepsell^{-4}$ and $j(\ell)\leq \frac{\log(1/\vepsell)}{\log 2}$ since $2^{j(\ell)}\leq 1/\vepsell$. Combining~\eqref{ineq:W1},~\eqref{ineq:W2} and~\eqref{ineq:W3} yields, for~$\ell$ sufficiently large, almost surely
\begin{equation*}
\begin{split}
W_{\vepsell}&\leq \chi_W+N_WC_R\left(\log(1/\vepsell)\right)^{1+\epsLanbda}C_\Stepsize\vepsell^{-4}\left(1/\log 2\right)^{1-\epsLanbda}\left(\log(1/\vepsell)\right)^{1-\epsLanbda}\\
&= \chi_W+C_W\vepsell^{-4}\left(\log(1/\vepsell)\right)^2\leq R_W\vepsell^{-4}\left(\log(1/\vepsell)\right)^2,
\end{split}
\end{equation*}
where $C_W$ and $R_W$ are finite almost surely, $R_W$ depends on $\ubarlong{K}_W, N_W, C_\Stepsize, \epsLanbda, k_0\SSaddnew{, R_\Stepsize, \ubar{\lambda}}$, and the last inequality follows from the fact that $\vepsell^{-4}\left(\log(1/\vepsell)\right)^2\to\infty$ as $\ell\to\infty$. This completes the proof.
\end{proof}

\subsection{Proof of Lemma~\ref{lemma:ExpInverseRadius}}\label{app:B-lemma:ExpInverseRadius}
\begin{proof}[Proof of Lemma~\ref{lemma:ExpInverseRadius}]
Fix $\varepsilon\in(0,\varepsilon_0/2]$ and let $j$ satisfy $\bar\varepsilon_j/2<\varepsilon\le\bar\varepsilon_j$; then $\bar\varepsilon_j<2\varepsilon\le\varepsilon_0$, so $j$ is admissible in Lemma~\ref{lemma:br}. Write $T:=\taurandom_{\bar\varepsilon_j/2}$, $X_k:=\Xkj$, $y_k:=\E{X_k^r}$ and $p_k:=\prob{k<T}$.

\emph{Reduction.} By Assumption~\ref{Assumption:general1}-\emph{(v)}, $\taurandom_\varepsilon\le T$. For $k_0\le k<T$, $$\Stepsize_k^{-2}\le(\kappa_a\bar\varepsilon_j)^{-2}\Zkj\le(\kappa_a\varepsilon)^{-2}X_k^r,$$ because $X_k\in\{0\}\cup[1,\infty)$ and $r\ge1$; for $k<k_0$, $\Stepsize_k\ge\gamma^{-k}\alpha_0$ by~\eqref{eq4:AssumpGen}. Hence, for any nondecreasing deterministic weights $w_k\ge0$, Tonelli's theorem gives
\begin{equation}\label{eq:EIRreduction}
\E{\sum_{k=0}^{\taurandom_\varepsilon-1}w_k\Stepsize_k^{-2}}\le w_{k_0}\frac{k_0\gamma^{2k_0}}{\alpha_0^2}+(\kappa_a\varepsilon)^{-2}\sum_{k\ge k_0}w_ky_k .
\end{equation}

\emph{Recursion.} Taking expectations in~\eqref{ineq:br} gives $y_{k+1}\le c_ry_k+b_rp_k$ for $k\ge k_0$, and $y_{k_0}\le1$ since $j$ is admissible. Hence
\begin{equation}\label{eq:EIRunroll}
y_k\le c_r^{k-k_0}+b_r\sum_{i=k_0}^{k-1}c_r^{k-1-i}p_i\qquad\text{for all }k\ge k_0 .
\end{equation}

\emph{Part (a).} Summing~\eqref{eq:EIRunroll} over $k\ge k_0$ and using $\sum_{i\ge k_0}p_i\le\E{T}\le\kappa_T\bar\varepsilon_j^{-2}\le\kappa_T\varepsilon^{-2}$, by {\bl \SSaddnew{Theorem~\ref{theor:GenIterCompl}}} at $\bar\varepsilon_j\le\varepsilon_0$, gives $$(1-c_r)\sum_{k\ge k_0}y_k\le1+b_r\kappa_T\varepsilon^{-2},$$ which completes the proof for ~(a) with~\eqref{eq:EIRreduction} with $w_k\equiv1$ and $\varepsilon\le1$.

\emph{Part (b).} In this part $\kappa$ denotes a deterministic constant depending only on the quantities listed in the statement, whose value may change from line to line. For $i,i'\ge0$, $i+3+i'\le(i+3)(1+i')$ and $\log(i+3)\ge1$ give $\lambda_{i+1+i'}\le\lambda_{i+1}(1+\log(1+i'))^{1+\epsLanbda}$, and $\sum_{i'\ge0}c_r^{i'}(1+\log(1+i'))^{1+\epsLanbda}<\infty$. Multiplying~\eqref{eq:EIRunroll} by $\lambda_k$, summing over $k\ge k_0$, and using $\sum_{i\ge k_0}\lambda_{i+1}p_i=\E{\sum_{i=k_0}^{T-1}\lambda_{i+1}}\le\E{T\lambda_T}$,
\begin{equation}\label{eq:EIRweighted}
\sum_{k\ge k_0}\lambda_ky_k\le\kappa\bigl(1+\E{T\lambda_T}\bigr)\le\kappa\bigl(1+\E{u(T)}\bigr).
\end{equation}
It remains to bound $\E{u(T)}$. 
The function $u$ is nondecreasing with $u(0)=0$, and $u(t)\le\kappa(1+t^2)$. For $c\ge1$, $ct+3\le c(t+3)$ and $\log(t+3)\ge1$ give $u(ct)\le c(1+\log c)^{1+\epsLanbda}u(t)$, hence $u(t_1+t_2+t_3)\le\kappa\,(u(t_1)+u(t_2)+u(t_3))$; for $c,t\ge0$, $ct+3\le(c+3)(t+3)$ gives $u(ct)\le2^{\epsLanbda}(c\,u(t)+t\,u(c))$. 
Use $T_{\bar k}$, $\SSaddnew{\Gamma^\pm_{\bar k}}$, $\Gamma_{\bar k}$ and $C_k$ from the proof of {\bl Theorem~\ref{theor:GenIterCompl}, with $a_\varepsilon$ replaced with $\kappa_a \bar\varepsilon_j$ }. Set $\Pi_k:=\prob{\Dscr_k|\gkalgebra}\ge\pi_\star$ and $O_{\bar k}:=\sum_{k=k_0}^{\bar k-1}C_k(\mathds 1_{\Dscr_k}-\Pi_k)$. Then $\Gamma_{\bar k}\ge\pi_\star \SSaddnew{\Gamma^-_{\bar k}}+O_{\bar k}$, and~\eqref{eq:ExpIterTelescope} gives
\[
(2\pi_\star-1)\SSaddnew{\Gamma^-_{\bar k}}\le \SSaddnew{\Gamma^+_{\bar k}}+\kappa_L+2\abs{O_{\bar k}}.
\]
Since $C_k,\Pi_k$ are $\gkalgebra$-measurable and $\Dscr_k\in\mathscr G_k$, $O_{\bar k}$ is a sum of martingale differences bounded by one, so
\[
\E{O_{\bar k}^2}=\E{\sum_{k=k_0}^{{\bar k}-1}C_k\Pi_k(1-\Pi_k)}\le\E{\SSaddnew{\Gamma^-_{\bar k}}}\le\kappa_T\varepsilon^{-2}.
\]
With $a:=(\kappa_a\varepsilon)^{-2}$ and $\SSaddnew{\Gamma^+_{\bar k}}\le aS_\Stepsize$,
\[
\E{u(\SSaddnew{\Gamma^+_{\bar k}})}\le2^{\epsLanbda}\bigl(u(a)\E{S_\Stepsize}+a\,\E{u(S_\Stepsize)}\bigr)\le\kappa\,\varepsilon^{-2}(1+\log(1/\varepsilon))^{1+\epsLanbda},
\]
because $a+3\le(\kappa_a^{-2}+3)\varepsilon^{-2}$. Hence
\[
\E{u(T_{\bar k})}\le\kappa\bigl(1+\E{u(\SSaddnew{\Gamma^-_{\bar k}})}+\E{u(\SSaddnew{\Gamma^+_{\bar k}})}\bigr)\le\kappa\bigl(1+\E{u(\SSaddnew{\Gamma^+_{\bar k}})}+\E{O_{\bar k}^2}\bigr)\le\kappa\,\varepsilon^{-2}(1+\log(1/\varepsilon))^{1+\epsLanbda},
\]
and monotone convergence as ${\bar k}\to\infty$ gives the same bound for $\E{u(T)}$. Finally, \eqref{eq:EIRreduction} with $w_k=\lambda_k$ together with~\eqref{eq:EIRweighted} gives~(b).
\end{proof}

\section{Proof of Proposition~\ref{prop:finiteNk}}\label{app:finitenessOriginal}
\SSaddnew{The stopping-time property and the lower bound~\eqref{eq:NkFloorLB} are proved in Section~\ref{sec:estiModels}; what remains is the finiteness of the sample sizes, that is, the finiteness of the Chow--Robbins rule~\cite{ChowRobbins1965} applied conditionally on a random design. It is proved with the conditional strong law of large numbers of~\cite{MaNoZiCondSLLN2005}, whose two definitions and theorem are recalled first. For \ASTRAF\ the order matters: $\nskrandom$ is defined through the trial point, which exists only once $\nkrandom$ is finite, so $\nkrandom$ is treated first.}
\begin{definition}(\cite[Definition~2.1]{MaNoZiCondSLLN2005})
Random variables $Y_1,Y_2,\dots,Y_n$ are $\halgebra$-independent or conditionally independent given~$\halgebra$, if for every choice of Borel sets $\mathscr{A}_1,\mathscr{A}_2,\dots,\mathscr{A}_n$,
\begin{equation*}
\prob{\left.\bigcap_{i=1}^n\accolade{Y_i\in\mathscr{A}_i}\right\rvert\halgebra}=\prod_{i=1}^{n}\prob{Y_i\in\mathscr{A}_i|\halgebra}\qquad\text{almost surely.}
\end{equation*}
\end{definition}
\begin{definition}(\cite[Definition~4.1]{MaNoZiCondSLLN2005})
Two random variables~$Y$ and~$Z$ have the same conditional distributions given~$\halgebra$~if 
\begin{equation*}
\prob{Y\leq a|\halgebra}=\prob{Z\leq a|\halgebra}\qquad\text{almost surely for every}\ a\in\R.
\end{equation*}
$\accolade{Z_\ell}_{\ell\geq 1}$ is a sequence of random variables with the same conditional distribution given $\halgebra$ if, for all~$a\in\R$
\begin{equation*}
\prob{Z_\ell\leq a|\halgebra}=\prob{Z_1\leq a|\halgebra}\qquad\text{almost surely for every}\ \ell\geq 1.
\end{equation*}
\end{definition}

Next is stated the \SSaddnew{version} of Kolmogorov's strong law of large numbers for conditionally identical distributions.
\begin{theorem}\label{Theorem:42}(\cite[Theorem~4.2]{MaNoZiCondSLLN2005})
Let $(\Omega,\mathscr{F},\pr)$ be a probability space, let $\halgebra\subseteq\mathscr{F}$ be a sub-$\sigma$-algebra, and $\accolade{Z_\ell}_{\ell\geq 1}$ be a sequence of $\halgebra$-independent random variables with the same conditional distribution. \SSaddnew{Let}~$Z$ denote a random variable with this common conditional distribution. 
\SSaddnew{If }
$\E{\abs{Z}|\halgebra}<\infty$ almost surely, \SSaddnew{then}
\begin{equation*}
\frac{1}{n}\sum_{\ell=1}^{n}Z_\ell\longrightarrow\E{Z|\halgebra}\qquad\text{almost surely,}\quad \text{as}\ n\to\infty.
\end{equation*}
\end{theorem}

\SSaddnew{\begin{proof}[Proof of Proposition~\ref{prop:finiteNk}]
Fix $k\geq 0$, consider $\Xki=\Xzerok+\QkRandom\siirandom_k$, $i\in\intbracket{0,q}$, and define for $\ell\geq 1$, $\Zkli:=F(\Xki;\Theta_{k,\ell})$, where the random variables $\Theta_{k,\ell}$ are independent copies of~$\Theta$. For $\x\in\rd$ and $\mathscr{A}\in\mathscr{B}(\R)$, let $\zeta(\x;\mathscr{A}):=\prob{F(\x;\Theta)\in\mathscr{A}}$. Since by construction $\Xki$ is $\fqalgebra$-measurable, $\sigma(\Theta_{k,\ell})$ and~$\fqalgebra$ are independent, and $\E{\abs{\mathds{1}_{\accolade{\Zkli\in\mathscr{A}}}}}<\infty$ trivially, it follows from Lemma~\ref{Lemma:24} that almost surely,
\begin{equation*}
\E{\left.\mathds{1}_{\accolade{F(\Xki;\Theta_{k,\ell})\in\mathscr{A}}}\right\rvert\fqalgebra}=\prob{\Zkli\in\mathscr{A}|\fqalgebra}=\zeta(\Xki;\mathscr{A}).
\end{equation*}
Hence, the conditional distribution of~$\Zkli$ given~$\fqalgebra$ is independent of~$\ell$. Moreover, for distinct indices~$\ell_1,\dots,\ell_r$ and Borel sets $\mathscr{A}_1,\dots,\mathscr{A}_r$, since the $\Theta_{k,\ell_j}, j=1,\dots,r$, are mutually independent, then almost surely,
\begin{equation*}
\prob{\left.\bigcap_{j=1}^r\accolade{Z^i_{k,\ell_j}\in\mathscr{A}_j}\right\rvert\fqalgebra}=\prod_{j=1}^r\prob{Z^i_{k,\ell_j}\in\mathscr{A}_j|\fqalgebra}.
\end{equation*}
Thus, $\accolade{\Zkli}_{\ell\geq 1}$ is $\fqalgebra$-independent and identically conditionally distributed. It follows from Theorem~\ref{Theorem:42} that, almost surely,
\begin{equation}\label{FbarConv1}
\Fbarik(\siirandom_k,n;\Theta_k):=\frac{1}{n}\sum_{\ell=1}^{n}F(\Xki;\Theta_{k,\ell})\underset{n\to\infty}{\longrightarrow}\E{\Zkli|\fqalgebra}=\Esp_{\Theta}\left[F(\Xki;\Theta_{k,\ell})\right]=f(\Xki),
\end{equation}
where the second equality follows from Lemma~\ref{Lemma:24} (applying Lemma~\ref{Lemma:24} to $F^+$ and $F^-$ separately, each with finite conditional expectation since $\mathbb{E}_\Theta{|F(\x;\Theta)|}<\infty$ for every $\x$). On the other hand, defining $\mu_2(\x):=\Esp_{\Theta}[F(\x;\Theta)^2]=\sigma_F^2(\x)+(f(\x))^2$ which is finite thanks to Assumption~\ref{Assumpt:31}, we also observe from Theorem~\ref{Theorem:42} that almost surely,
\begin{equation}\label{FbarSqConv1}
\frac{1}{n}\sum_{\ell=1}^{n}F(\Xki;\Theta_{k,\ell})^2\underset{n\to\infty}{\longrightarrow}\mu_2(\Xki),
\end{equation}
given that the sequence $\accolade{(\Zkli)^2}_{\ell\geq 1}$ is again $\fqalgebra$-independent and identically conditionally distributed since $z\mapsto z^2$ is Borel-measurable. Observing that $\hat{\sigma}_F^2(\siirandom_k,n;\Theta_k):= \frac{1}{n}\sum_{\ell=1}^{n}F(\Xki;\Theta_{k,\ell})^2-\left(\Fbarik(\sii_k,n;\Theta_k)\right)^2$ equivalently to~\eqref{sigmaHatFsqEq}, it follows from~\eqref{FbarConv1} and~\eqref{FbarSqConv1} that almost surely,
\begin{equation*}
\hat{\sigma}_F^2(\siirandom_k,n;\Theta_k)\underset{n\to\infty}{\longrightarrow}\sigma_F^2(\Xki),
\end{equation*}
hence by continuity of $t\mapsto\sqrt{t}$, for each $i\in\intbracket{0,q}$, the event $\mathscr{E}_{\sigma,k}^i:=\accolade{\hat{\sigma}_F(\siirandom_k,n;\Theta_k)\underset{n\to\infty}{\longrightarrow}\sigma_F(\Xki)}$ is almost sure. Thus, since $\ds{\bigcap_{i=0}^q}\mathscr{E}_{\sigma,k}^i\subseteq\accolade{\hat{\sigma}_F^{\max}(n;\Theta_k)\underset{n\to\infty}{\longrightarrow}\underset{i\in\intbracket{0,q}}{\max}\sigma_F(\Xki)}$, then almost surely,
\begin{equation}\label{SigmaMaxConv1}
\hat{\sigma}_F^{\max}(n;\Theta_k)\underset{n\to\infty}{\longrightarrow}\underset{i\in\intbracket{0,q}}{\max}\sigma_F(\Xki).
\end{equation}
It follows from~\eqref{SigmaMaxConv1} that $\sigma_0\vee \hat{\sigma}_F^{\max}(n;\Theta_k)$ converges almost surely to a finite random variable, and is therefore eventually bounded almost surely, hence the event $\mathscr{E}_{\hat{\sigma},k}:=\accolade{\frac{\sigma_0\vee \hat{\sigma}_F^{\max}(n;\Theta_k)}{\sqrt{n}}\underset{n\to\infty}{\longrightarrow}0}$ is almost sure. We note that on $\mathscr{E}_{\hat{\sigma},k}$, there exists a finite random integer $N^\star_k$ such that for every $n\geq N^\star_k$, $\frac{\sigma_0\vee \hat{\sigma}_F^{\max}(n;\Theta_k)}{\sqrt{n}}\leq \kappaas\frac{\Dk}{\sqrt{\lambda_k}}$. Thus, the set $\accolade{n\in\N_+:\frac{\sigma_0\vee \hat{\sigma}_F^{\max}(n;\Theta_k)}{\sqrt{n}}\leq \kappaas\frac{\Dk}{\sqrt{\lambda_k}}}$ is nonempty and consequently, $N_k<\infty$ almost surely. We observe that similar arguments apply to the estimates $\Fbarokbrak(\bm{0},\nskrandom;\Theta_k^{0,s})$ and $\Fbarsk(\srandom_k,\nskrandom;\Theta_k^{0,s})$, but with the conditioning $\sigma$-algebra $\fmqalgebra$. Indeed, $\Xzerok$ and $\X_k^s$ are $\fmqalgebra$-measurable and $\sigma(\Theta_k^{0,s})$ and $\fmqalgebra$ are independent. Thus, by defining $Z_{k,\ell}^{[t]}:=F(\X_k^t;\Theta_{k,\ell}^{0,s}), t\in \accolade{0,s}$, the use of Lemma~\ref{Lemma:24} as above, shows that the sequence $\accolade{Z_{k,\ell}^{[t]}}_{\ell\geq 1}$ is $\fmqalgebra$-independent and identically conditionally distributed. Hence, recalling~\eqref{sigm0Fdef} and~\eqref{sigmsFdef}, Theorem~\ref{Theorem:42} applied to the observations $Z_{k,\ell}^{[t]}$ and their squares, gives for each $t\in\accolade{0,s}$, almost surely, 
\begin{equation}\label{sigmatConv1}
\hat{\sigma}_F^{[t]}(n;\Theta_k)\underset{n\to\infty}{\longrightarrow}\sigma_F(\X_k^t).
\end{equation}
It follows from~\eqref{sigmatConv1} that almost surely,
\begin{equation*}
\sigma_0\vee\hat{\sigma}_F^{[0]}(n;\Theta_k)\vee\hat{\sigma}_F^{[s]}(n;\Theta_k)\underset{n\to\infty}{\longrightarrow}\sigma_0\vee\sigma_F(\Xzerok)\vee\sigma_F(\X_k^s)<\infty,
\end{equation*}
hence the stopping condition is also satisfied in~\eqref{n0sDefEq}. Consequently, $\nskrandom<\infty$ almost surely. 
\SSaddnew{Hence $\nkrandom$ and $\nskrandom$ are finite almost surely for every $k$, and the  intersection of these almost sure events is almost sure.}
\end{proof}
}

\section{Proof of Corollary~\ref{theor31bis}}\label{app:pointwise}
\SSaddnew{Corollary~\ref{theor31bis} restates the bound $1-2e^{-\chifde\lambda_k}$ of Theorem~\ref{Theorem:22} at the level $1-k^{-(1+\ef)}$ used in Section~\ref{sec:convAnalysis}; $\kappaf$ is the index beyond which $\chifde\lambda_k$ exceeds $(1+\ef)\log k+\log2$. The same calibration with $c_F$ in place of $\chifde$ gives~\eqref{eq:UniformBeta} for \ASTRAR.}
\begin{proof}[Proof of Corollary~\ref{theor31bis}]
For any $\ef>0$ fixed, we note that since $\frac{(\log(k+2))^{1+\epsLanbda}}{\log k}\longrightarrow\infty$ and $\frac{\log 2}{\log k}\longrightarrow 0$ as $k\to \infty$, there exists $\kappaf\in\N$ sufficiently large, such that for all $k\geq \kappaf$, 
\begin{equation*}
\frac{\chifde\ubar{\lambda}(\log(k+2))^{1+\epsLanbda}}{(1+\ef)\log k}\geq \frac{\log 2}{(1+\ef)\log k}+1.
\end{equation*}
Thus, the quantity $r_k:=2e^{-\chifde\lambda_k}$ from the right-hand side of~\eqref{eq:chifde} satisfies for all $k\geq \kappaf$,
\begin{equation*}
-\log r_k=\chifde\ubar{\lambda}(\log(k+2))^{1+\epsLanbda}-\log 2\geq (1+\ef)\log k,
\end{equation*}
which implies that $r_k\leq \frac{1}{k^{1+\ef}}$, and the proof is completed.
\end{proof}

\section{Proofs for Section~\ref{sec:convAnalysis}}\label{app:proofsSec5}

\SSaddnew{Lemmas~\ref{lemma4p5}--\ref{lemma4p8} are the counterparts for the present sampling scheme of~\cite[Lemmas~4.1--4.4]{DzaWildSub2022}, and their proofs follow those there with the accuracy constants of Section~\ref{sec:Prob}; Lemma~\ref{lemmaIncr} adds to~\cite[Lemma~4.5]{DzaWildSub2022} the bound~\eqref{fourpt6bis}, used on $\mathscr V_k^c$. The \ASTRAF\ case of Theorem~\ref{convDeltaTheorem} follows the three-subcase argument of~\cite{DzaWildSub2022} and uses the preliminaries stated before the \ASTRAR\ proof in Section~\ref{sec:convAnalysis}; it is where the factor $\pi_{\min}^2$ and the $\sigma$-algebra $\fmqalgebra$ enter. The proofs of Theorem~\ref{theorem32} and Lemma~\ref{l1deltakconv} end the appendix.}

\subsection{Proofs of Lemmas~\ref{lemma4p5}--\ref{lemmaIncr}}
\begin{proof}[Proof of Lemma~\ref{lemma4p5}]
We first observe from~\eqref{firstDeltaCond} that $\normii{\gk}\geq\SSaddnew{\kappah}\dk$, which, combined with~\eqref{CauchyDecreaseEq} leads to
\begin{align}\label{CauchyDecreaseEq2}
\mk(\bm{0})-\mk(\sk)\geq\frac{\kfcd}{2}\normii{\gk}\dk.
\end{align}
It holds that
\begin{eqnarray*}
f(\xk+\Qk\sk)-f(\xk)&=&[\textcolor{black}{(f(\xk+\Qk\sk)-\mk(\sk))-(f(\xk)-\mk(\bm{0}))}]+\mk(\sk)-\mk(\bm{0})\\
&\leq&\textcolor{black}{\frac{1}{2}\kmfe\dk^2}-\frac{\kfcd}{2}\normii{\gk}\dk\leq-\frac{1}{4}\kfcd\normii{\gk}\dk,
\end{eqnarray*}
where the first inequality follows from Lemma~\ref{modFuncErrLemma} and~\eqref{CauchyDecreaseEq2}, while the last one results from the fact that $\kmfe\dk\leq \frac{1}{2}\kfcd\normii{\gk}$ thanks to~\eqref{firstDeltaCond}, and the proof is completed.
\end{proof}

\begin{proof}[Proof of Lemma~\ref{lemma4p6}]
We first observe from the subspace-accuracy assumption that
\begin{align}\label{starEq}
\normii{\gk}\geq\normii{\QkTranspose\nabla f(\xk)}-\frac{1}{2}\kmge\dk,
\end{align}
which, combined with $\normii{\QkTranspose\nabla f(\xk)}\geq \max\accolade{\SSaddnew{\kappah}+\frac{1}{2}\kmge, \frac{2\kmfe}{\kfcd}+\frac{1}{2}\kmge}\dk$ (thanks to~\eqref{secDeltaCond}) leads to~\eqref{firstDeltaCond}. Thus, it follows from Lemma~\ref{lemma4p5} that~\eqref{decrLemma4p5} holds. On the other hand,~\eqref{starEq} combined with~\eqref{secDeltaCond} yields
\begin{align}\label{starEqPrime}
\normii{\gk}\geq 2\max\accolade{\frac{\SSaddnew{\kappah}}{2\SSaddnew{\kappah}+\kmge},\frac{2\kmfe}{4\kmfe+\kmge\kfcd}}\normii{\QkTranspose\nabla f(\xk)}.
\end{align}
Combining~\eqref{decrLemma4p5} and~\eqref{starEqPrime} leads to~\eqref{seconImprov}, and the proof is completed.
\end{proof}

\begin{proof}[Proof of Lemma~\ref{lemma4p7}]
\SSaddnew{Since~\eqref{thirdDeltCond} implies~\eqref{firstDeltaCond}, \eqref{CauchyDecreaseEq2} holds, and we can write}
\begin{eqnarray*}
\abs{\hat{\rho}_k-1}&\leq& \frac{\textcolor{black}{\abs{(\fok-\fsk)-(f(\xk)-f(\xk+\Qk\sk))}}}{\mk(\bm{0})-\mk(\sk)}+\frac{\abs{(f(\xk)-\mk(\bm{0}))-(f(\xk+\Qk\sk)-\mk(\sk))}}{\mk(\bm{0})-\mk(\sk)}\\
&\leq& \frac{\textcolor{black}{\kfde\dk^2}+\textcolor{black}{\frac{1}{2}\kmfe\dk^2}}{\frac{1}{2}\kfcd\normii{\gk}\dk}\leq 1-\eta,
\end{eqnarray*}
where the second inequality follows from the $\kfde$-weak accuracy assumption,  Lemma~\ref{modFuncErrLemma} and~\eqref{CauchyDecreaseEq2}, while the last one follows from~\eqref{thirdDeltCond}. Thus, $\hat{\rho}_k\geq \eta$. Moreover,~\eqref{thirdDeltCond} implies that $\mu\normii{\gk}\geq\dk$, whence the $k$th iteration is successful.
\end{proof}

\begin{proof}[Proof of Lemma~\ref{lemma4p8}]
A successful iteration implies that $\mu\normii{\gk}\geq\dk$ and $\hat{\rho}_k\geq\eta$. Consequently,
\begin{eqnarray*}
\fok-\fsk&\geq&\eta(\mk(\bm{0})-\mk(\sk))
\geq\eta\frac{\kfcd}{2}\normii{\gk}\min\accolade{\frac{\normii{\gk}}{\SSaddnew{\max\{\normeshort{\Hk},1\}}},\dk}\geq\frac{\eta}{2\mu}\kfcd\min\accolade{\frac{1}{\SSaddnew{\kappah}\mu},1}\dk^2.
\end{eqnarray*}
Then the proof is completed by observing that  
\begin{eqnarray*}
f(\xk+\Qk\sk)-f(\xk)&=&(f(\xk+\Qk\sk)-\fsk)-(f(\xk)-\fok)+\fsk-\fok\\
&\leq&\kfde\dk^2-\frac{\eta}{2\mu}\kfcd\min\accolade{\frac{1}{\SSaddnew{\kappah}\mu},1}\dk^2=-c_2\dk^2.
\end{eqnarray*}
\end{proof}

{\bl
\begin{proof}[Proof of Lemma~\ref{lemmaIncr}]
The proof of~\eqref{fourpt6} is identical to that of \cite[Lemma~4.5]{DzaWildSub2022}. In order to show~\eqref{fourpt6bis}, we note that since $\nabla f$ is $\kappa_{Lg}$-Lipschitz, then Taylor inequality gives, for every trial step satisfying $\normii{\sk}\leq \dk$,
\begin{equation*}
\begin{split}
f(\xk+\Qk\sk)-f(\xk)&\leq \nabla f(\xk)^\top(\Qk\sk)+\frac{\kappa_{Lg}}{2}\normii{\Qk\sk}^2\leq \normii{\Qk^\top \nabla f(\xk)}\dk+\frac{1}{2}\kappa_{Lg}\Qmax^2\dk^2\\
&< \zeta\dk^2+\frac{1}{2}\kappa_{Lg}\Qmax^2\dk^2<c_3(\zeta)\zeta\dk^2,
\end{split}
\end{equation*}
where the third inequality follows from the assumption $\normii{\Qk^\top \nabla f(\xk)}<\zeta\dk$, and the proof is completed.
\end{proof}
}

\subsection{Proof of Theorem~\ref{convDeltaTheorem} for \ASTRAF}\label{app:proofConvF}
\begin{proof}
The result is proved by adapting the proof of~\cite[Theorem~4.8]{DzaWildSub2022}. \SSaddnew{The choice of $k_0$ in Theorem~\ref{convDeltaTheorem}} 
guarantees that for all $k\geq k_0$, $1-\frac{1}{k^{1+\ef}}>1-\frac{2p}{k^{1+\ef}}\geq \SSaddnew{\pi_{\min}}$ so that from~\eqref{jkEq1} and~\eqref{mgeEquation1}, almost surely,
\begin{equation}\label{aminIkJkEqs}
\prob{\left.\jk\right\rvert\fmqalgebra}\geq \SSaddnew{\pi_{\min}}\qquad\text{and}\qquad \prob{\ik\left\lvert\fqalgebra\right.}\geq \SSaddnew{\pi_{\min}},\quad\text{for all}\ k\geq k_0.
\end{equation}
\SSaddnew{Throughout, $k\ge k_0$, and we use the notation of the preliminaries preceding the proof for \ASTRAR; it remains to prove~\eqref{fourpt15} and~\eqref{fourpt16}.}

{\textbf{Case~1:} $\ivk=1$.} Three subcases are considered next.\\
{\textit{(Subcase~1i)}} Good model $(\iik=1)$ and good estimates $(\ijk=1)$. It follows from~\eqref{fourpt8} that $$\normii{\QkTranspose\nabla f(\xzerok)}\geq 
\SSaddnew{\left(\frac{1}{2}\kmge+\kappa_{\rm suc}\right)}
\dk,$$ implying in particular~\eqref{secDeltaCond}, whence there is a decrease in~$f$ according to~\eqref{seconImprov}. Moreover, by $(\kmge,\Qk)$-\SSaddnew{subspace-accuracy}, 
\[\normii{\gkhat}\geq\normii{\QkTranspose\nabla f(\xzerok)}-\frac{1}{2}\kmge\dk\geq 
\SSaddnew{\kappa_{\rm suc}}
\dk,\]
 thus implying~\eqref{thirdDeltCond}, hence the iteration is successful. Thus, \SSaddnew{$\phikplusone-\phik\le b_{2,k}$ } follows from~\eqref{fourpt13} and~\eqref{seconImprov}. 
From~\eqref{fourpt17} 
\begin{equation}\label{fourpt18}
\ivk\iik\ijk(\Phikplusone-\Phik)\leq \ivk\iik\ijk\SSaddnew{\brandom_{2,k}}.
\end{equation}
{\textit{(Subcase~1ii)}} Bad model $(\iikc=1)$ and good estimates $(\ijk=1)$. Regardless of the iteration type, i.e., successful or not, the change in~$\phik$ can always be bounded as follows, using~\eqref{fourpt13} and~\eqref{fourpt6}:
\begin{equation}\label{fourpt19}
\phikplusone-\phik\leq \nu c_3\normii{\QkTranspose\nabla f(\xzerok)}\dk+(1-\nu)(\gamma^2-1)\dk^2=\SSaddnew{b_{3,k}}.
\end{equation}
It follows from~\eqref{fourpt19} that
\begin{equation}\label{fourpt20}
\ivk\iikc\ijk(\Phikplusone-\Phik)\leq \ivk\iikc\ijk\SSaddnew{\brandom_{3,k}}\leq \ivk\iikc\SSaddnew{\brandom_{3,k}}.
\end{equation}
{\textit{(Subcase~1iii)}} Bad estimates $(\ijkc=1)$. From Subcase~\textit{1ii}, we observe that the following always holds:
\begin{equation}\label{fourpt21}
\ivk\ijkc(\Phikplusone-\Phik)\leq \ivk\ijkc\SSaddnew{\brandom_{3,k}}.
\end{equation}
With the subcases now complete, we first observe, using ~\eqref{aminIkJkEqs}, the fact that $\iik$ is $\fmqalgebra$-measurable, and the inclusion $\fqalgebra\subseteq\fmqalgebra$, that
\begin{equation}\label{tower1}
\E{\iik\ijk|\fqalgebra}=\E{\left.\iik\E{\ijk|\fmqalgebra}\right\rvert\fqalgebra}\geq\SSaddnew{\pi_{\min}}\E{\iik|\fqalgebra}\geq \SSaddnew{\pi_{\min}^2}.
\end{equation}
Then, given that $\ivk, \SSaddnew{\brandom_{2,k}}$ and $\SSaddnew{\brandom_{3,k}}$ are $\fqalgebra$-measurable thanks to Lemma~\ref{lemmafourp7}, we combine~\eqref{fourpt18},~\eqref{fourpt20} and~\eqref{fourpt21}, and then take the expectation with respect to $\fqalgebra$ to obtain:
\begin{equation}\label{fourpt22}
\begin{split}
\E{\ivk(\Phikplusone-\Phik)|\fqalgebra}&\leq \ivk\SSaddnew{\brandom_{2,k}}\E{\iik\ijk|\fqalgebra}+ \ivk\SSaddnew{\brandom_{3,k}}\!\left(\E{\iikc|\fqalgebra}+\E{\ijkc|\fqalgebra} \right)\\
&\leq \ivk\left(\SSaddnew{\pi_{\min}^2}\SSaddnew{\brandom_{2,k}}+2\SSaddnew{\brandom_{3,k}}(1-\SSaddnew{\pi_{\min}})\right)=:\ivk\SSaddnew{\brandom_{23,k}},
\end{split}
\end{equation}
where the second inequality follows from~\eqref{aminIkJkEqs},~\eqref{tower1} and the fact that $\SSaddnew{\brandom_{2,k}}<0$ while $\SSaddnew{\brandom_{3,k}}>0$. We observe from the fact that $2(1-\SSaddnew{\pi_{\min}})+\SSaddnew{\pi_{\min}^2}=1+(1-\SSaddnew{\pi_{\min}})^2\leq 2$, that 
\[\SSaddnew{\brandom_{23,k}}\leq \nu\normii{\QkRandomTranspose\nabla f(\Xzerok)}\Dk(-\SSaddnew{\pi_{\min}^2}c_1+2(1-\SSaddnew{\pi_{\min}})c_3)+2(1-\nu)(\gamma^2-1)\Dk^2.\]
From~\eqref{fourpt9} and the first condition in~\eqref{fourpt10}, it holds that
\[\SSaddnew{\pi_{\min}^2}c_1-2(1-\SSaddnew{\pi_{\min}})c_3\geq\frac{1}{2}c_1\geq 2\frac{2(1-\nu)(\gamma^2-1)}{\nu\zeta}.\]
Hence, since $\normii{\QkRandomTranspose\nabla f(\Xzerok)}\geq \zeta\Dk$, we easily obtain 
\begin{equation}\label{eq:B23kineq}
\SSaddnew{\brandom_{23,k}}\leq-\frac{1}{4}c_1\nu\normii{\QkRandomTranspose\nabla f(\Xzerok)}\Dk\leq -2(1-\nu)(\gamma^2-1)\Dk^2,
\end{equation}
which, together with~\eqref{fourpt22}, yields~\eqref{fourpt15}.
 
{\textbf{Case~2:} $\ivkc=1$.} We consider the event $\SSaddnew{\mathscr L_k}:=\accolade{\mu\normii{\Gkhat}\geq \Dk}$ and then observe that by construction, $\normii{\Gkhat}$ and hence $\SSaddnew{\ilkhat}, \SSaddnew{\ilkhatc}$ are $\fmqalgebra$-measurable. We note that if $\SSaddnew{\mathscr L_k}$ does not hold, then $\mu\normii{\gkhat}<\dk$ and the iteration is unsuccessful, leading to~\eqref{fourpt14}. Consequently, 
\vspace*{-0.2cm}
\begin{equation}\label{GkbarCase2}
\begin{split}
\E{\ivkc\SSaddnew{\ilkhatc}(\Phikun-\Phik)|\fmqalgebra}&\leq -(1-\nu)\left(1-\frac{1}{\gamma^2}\right)\E{\ivkc\SSaddnew{\ilkhatc}\Dk^2|\fmqalgebra}\\
&\leq -\frac{1}{2}\ivkc\SSaddnew{\ilkhatc}\SSaddnew{\pi_{\min}}(1-\nu)\left(1-\frac{1}{\gamma^2}\right)\Dk^2,
\end{split}
\end{equation}
where the last inequality follows from the $\fmqalgebra$-measurability of $\ivkc$ and $\Dk$, and the fact that $\falgebra\subseteq\fqalgebra\subseteq\fmqalgebra$. 
Inspired by~\cite[Theorem~3]{blanchet2016convergence}, which improved on the proof of~\cite[Theorem~4.11, Case~2]{chen2018stochastic} where four subcases were considered, only two are analyzed below. The overall goal is to show that
\vspace*{-0.3cm}
\begin{equation}\label{GkCase2}
\E{\ivkc\SSaddnew{\ilkhat}(\Phikun-\Phik)|\fmqalgebra}\leq -\frac{1}{2}\ivkc\SSaddnew{\ilkhat}\SSaddnew{\pi_{\min}}(1-\nu)\left(1-\frac{1}{\gamma^2}\right)\Dk^2,
\end{equation}
which, combined with~\eqref{GkbarCase2}, yields \scalebox{1}{$\E{\ivkc(\Phikun-\Phik)|\fmqalgebra}\leq -\frac{1}{2}\ivkc\SSaddnew{\pi_{\min}}(1-\nu)\left(1-\frac{1}{\gamma^2}\right)\Dk^2$},
hence
{\small
\begin{equation*}
\begin{split}
\E{\ivkc(\Phikun-\Phik)|\fqalgebra}&=\E{\E{\ivkc(\Phikun-\Phik)|\fmqalgebra}|\fqalgebra}\\
&\leq \E{-\frac{1}{2}\ivkc\SSaddnew{\pi_{\min}}(1-\nu)\left(1-\frac{1}{\gamma^2}\right)\Dk^2 \, |\fqalgebra}=-\frac{1}{2}\ivkc\SSaddnew{\pi_{\min}}(1-\nu)\left(1-\frac{1}{\gamma^2}\right)\Dk^2,
\end{split}
\end{equation*}
}
showing~\eqref{fourpt16}, with the last equality following from the $\fqalgebra$-measurability of $\Dk$ and $\ivkc$.  What remains to be proved is~\eqref{GkCase2} as achieved below, assuming that $\SSaddnew{\mathscr L_k}$ holds.\\
\textit{(Subcase 2i)}
Good estimates ($\ijk=1$) and $\SSaddnew{\ilkhat}=1$. While $f$ is decreased on successful iterations because  of good estimates thanks to Lemma~\ref{lemma4p8} and in which case~\eqref{eq:fourpt5} holds,  the trust region radius is reduced on unsuccessful iterations, leading to~\eqref{fourpt14}. However, combining~\eqref{eq:fourpt5} and~\eqref{fourpt13} leads to
$\phi_{k+1}-\phi_k\leq \left(-\nu c_2+(1-\nu)(\gamma^2-1)\right)\dk^2\leq \SSaddnew{-b_1\dk^2}
,$
where the last inequality follows from~\eqref{fourpt9}, thus showing that~\eqref{fourpt14} holds in any case. Hence,
\begin{equation}
\label{GammakBarGkJk}
\ivkc\SSaddnew{\ilkhat}\ijk(\Phikun-\Phik) \leq -\ivkc\SSaddnew{\ilkhat}\ijk(1-\nu)\left(1-\frac{1}{\gamma^2}\right)\Dk^2.
\end{equation}
\textit{(Subcase 2ii)}
Bad estimates ($\ijkc=1$) and $\SSaddnew{\ilkhat}=1$. In this case, a successful step can lead to an increase in the objective function {\bl according to~\eqref{fourpt6bis}}, which combined with~\eqref{fourpt13} yields
{\bl
\begin{equation}\label{GammakBarGkJkBar}
\ivkc\SSaddnew{\ilkhat}\ijkc(\Phikun-\Phik) \leq \ivkc\SSaddnew{\ilkhat}\ijkc\left(\nu c_3\zeta+(1-\nu)(\gamma^2-1)\right)\Dk^2.
\end{equation}
}
With the subcases complete, then  combining~\eqref{GammakBarGkJk} and~\eqref{GammakBarGkJkBar}, taking expectations with respect to~$\fmqalgebra$, and recalling that $\E{\ijk|\fmqalgebra}\geq\SSaddnew{\pi_{\min}}$, yields
\begin{eqnarray*}
\E{\ivkc\SSaddnew{\ilkhat}(\Phikun-\Phik)|\fmqalgebra}&\leq& \ivkc\SSaddnew{\ilkhat}\left(-\SSaddnew{\pi_{\min}}(1-\nu)(1-1/\gamma^2)\right. \left. +\ (1-\SSaddnew{\pi_{\min}})\left(\nu c_3\zeta+(1-\nu)(\gamma^2-1)\right) \right)\Dk^2\\
&\leq& -\frac{1}{2}\ivkc\SSaddnew{\ilkhat}\SSaddnew{\pi_{\min}}(1-\nu)\left(1-\frac{1}{\gamma^2}\right)\Dk^2,
\end{eqnarray*}
thus proving~\eqref{GkCase2}, where the second inequality follows from the second condition in~\eqref{fourpt10}, and the proof is completed.
\end{proof}

\subsection{Proof of Theorem~\ref{theorem32}}
\begin{proof}[Proof of Theorem~\ref{theorem32}] 
\SSaddnew{We first give the argument for \ASTRAF.} 
Inspired by the proofs of~\cite[Theorem~4.9]{DzaWildSub2022} and~\cite[Theorem~4.16]{chen2018stochastic}, the result is proved by contradiction, conditioned on the event $\accolade{\Dk\to 0}$ which is almost sure, thanks to Theorem~\ref{convDeltaTheorem}. Assume that with nonzero probability there exist {\bl a deterministic} $\varepsilon>0$ and $\bar{\kappa}\in \N$ such that $\normii{\nabla f(\Xzerok)}\geq \varepsilon$ for all $k\geq \bar{\kappa}$.
Let $\accolade{\xzerok}_{k\geq\bar{\kappa}}$ and $\accolade{\dk}_{k\geq\bar{\kappa}}$ be realizations of $\accolade{\Xzerok}_{k\geq\bar{\kappa}}$ and $\accolade{\Dk}_{k\geq\bar{\kappa}}$, respectively, for which $\normii{\nabla f(\xzerok)}\geq \varepsilon$ for all~$k\geq\bar{\kappa}$. Recall $k_0$ from Theorem~\ref{convDeltaTheorem} for which~\eqref{aminIkJkEqs} holds for all $k\geq k_0$. Since $\dk\to 0$, there exists $k_1\in \N$ such that for all $k\geq \bar{\kappa}_{0,1}:=\max\accolade{\bar{\kappa}, k_0, k_1}$,
\begin{equation}\label{deltabEq}
\dk <b:=\alpha_Q\cdot\min\accolade{
\SSaddnew{\frac{\varepsilon}{\frac{1}{2}\kmge+\kappa_{\rm suc}}}
, \frac{\delta_{\max}}{\gamma}}.
\end{equation}
Let $\rrandomk$ be the random variable with realizations $r_k:=\log_\gamma\left(\frac{\dk}{b}\right)$. Then for all $k\geq \bar{\kappa}_{0,1}$, $r_k<0$. The main idea of this proof is to show that such realizations occur only with probability zero, leading to a contradiction. 

{\bl
Recall the events~$\ik$ and~$\jk$ satisfying~\eqref{aminIkJkEqs}, and~$\AkEvent$ satisfying~\eqref{WAMEquation}, and define $\Dscr_k:=\ik\cap\jk\cap\AkEvent$.
For every $k\geq k_0$, it follows from~\eqref{aminIkJkEqs},~\eqref{WAMEquation}, and the tower property that $\prob{\Dscr_k|\falgebra}\geq \pi_{\min}^2(1-\beta_Q)=:\pi_\star>\frac{1}{2}$,
where the last inequality follows from the hypothesis on~$\beta_Q$. Consider now a given realization as above and an iteration $k\ge\bar\kappa_{0,1}$ for which $\Dscr_k$ occurs.}  
It follows from~\eqref{deltabEq} that
\[\normii{\nabla f(\xzerok)}\geq\frac{1}{\alpha_Q}
\SSaddnew{\left(\frac{1}{2}\kmge+\kappa_{\rm suc}\right)}
{\bl \dk},\] which, together with well alignment and subspace-accuracy, yields 
\begin{equation*}
\begin{split}
\normii{\gkhat}&\geq \normii{\QkTranspose \nabla f(\xzerok)}-\frac{1}{2}\kmge\dk\geq \alpha_Q\normii{\nabla f(\xzerok)}-\frac{1}{2}\kmge\dk\geq 
\SSaddnew{\kappa_{\rm suc}}
\dk.
\end{split}
\end{equation*}
Lemma~\ref{lemma4p7} guarantees that the $k$th iteration is successful, hence $\dkun=\gamma\dk$, leading to $r_{k+1}=r_k+1$. For all other outcomes of $\ik$, $\jk$ and $\AkEvent$, 
it always holds that $\dkun\geq \gamma^{-1}\dk$, hence $r_{k+1}\geq r_k-1$. Consequently, {\bl along every such realization and for every $k\ge\bar\kappa_{0,1}$,
\begin{equation}\label{eq:centtrentrois}
r_{k+1}-r_k\geq 2\mathds 1_{\Dscr_k}-1.
\end{equation}
Now define the random walk \(\{\widetilde W_k\}_{k\geq k_0}\) globally by $\widetilde W_k:= \sum_{i=k_0}^{k-1}\bigl(2\mathds 1_{\Dscr_i}-1\bigr)$ for all $k>k_0$ and $\widetilde W_{k_0}:=0$.
For every $k\geq k_0$, $\E{\widetilde W_{k+1}-\widetilde W_k|\falgebra}=2\prob{\Dscr_k|\falgebra}-1\geq 2\pi_\star-1>0$. Since $\Dscr_i\in\mathscr F_i$, the process $\accolade{\widetilde W_k}_{k\geq k_0}$ is adapted to $\accolade{\falgebra}_{k\geq k_0}$ and satisfies $\abs{\widetilde W_k}\leq k-k_0$. Hence, it is an integrable submartingale with respect to this filtration, with increments in $\accolade{-1,1}$. Moreover, $\abs{\widetilde W_{k+1}-\widetilde W_k}=1$ for every $k\ge k_0$ so $\accolade{\widetilde W_k}$ cannot converge to a finite value. Therefore, by~\cite[Theorem~4.4]{chen2018stochastic},
\begin{equation}\label{eq:centtrentequatre} 
\prob{\limsup_{k\to\infty}\widetilde W_k=+\infty}=1.
\end{equation}
Returning to the realization fixed above, summing~\eqref{eq:centtrentrois} from~$\bar\kappa_{0,1}$ to $k-1$ gives, for every $k>\bar\kappa_{0,1}$,
\begin{equation}\label{eq:centtrentecinq}
r_k-r_{\bar\kappa_{0,1}}\geq \widetilde w_k-\widetilde w_{\bar\kappa_{0,1}},
\end{equation}
where~$\widetilde w_k$ denotes the corresponding realization of~$\widetilde W_k$. On the almost sure event in~\eqref{eq:centtrentequatre}, the right-hand side of~\eqref{eq:centtrentecinq} has infinite $\limsup$, and therefore $\limsup_{k\to\infty}r_k=+\infty$. However, by construction, $r_k<0$ for every $k\geq\bar\kappa_{0,1}$, which is a contradiction. Thus, the sequence of realizations satisfying $\normii{\nabla f(\xzerok)}\geq\varepsilon$ for all sufficiently large~$k$ occurs with probability zero, which proves the result. 

For \ASTRAR, the same argument applies with \SSaddnew{$\Dscr_k:=\AkEvent\cap\ik\cap\jk$, where $\jk$ is built from the reuse estimates, and $\pi_\star:=\pi_{\min}(1-\beta_Q)$, by Remark~\ref{rem:theorem32R}}.
}\end{proof}

\subsection{Proof of Lemma~\ref{l1deltakconv}}
\begin{proof}[Proof of Lemma~\ref{l1deltakconv}]
We know from Theorem~\ref{convDeltaTheorem} that $\Dk\to 0$ almost surely. Recalling the constants $k_0$ and $\zeta$ from Theorem~\ref{convDeltaTheorem}, it holds for each realization of \SSadd{Algorithm~\ref{algoF} or~\ref{algoR}} that there exists a finite integer $k_1\geq k_0$ such that $\dk\leq \SSaddnew{\varepsilon}\alpha_Q/\zeta$ for all $k\geq k_1$. {\bl Let $\krandom_1$ be the random variable with realizations~$k_1$, and define $\mathscr{K}_{\SSaddnew{\varepsilon},1}:=\accolade{k\in \mathscr{K}_\SSaddnew{\varepsilon}:k\geq \krandom_1}$ and
$\mathscr{C}_k:=\accolade{\normii{\nabla f(\Xzerok)}>\SSaddnew{\varepsilon}, \Dk\leq \SSaddnew{\varepsilon}\alpha_Q/\zeta}$. Consider $\mathscr{W}_k:=\accolade{\normii{\QkRandomTranspose\nabla f(\Xzerok)}\geq \alpha_Q\normii{\nabla f(\Xzerok)}}$} satisfying $\prob{\mathscr{W}_k|\falgebra}\geq 1-\beta_Q$ almost surely, and assume that $\mathscr{C}_k\cap\mathscr{W}_k$ holds, which means in particular that 
the event $\mathscr{V}_k:=\accolade{\normii{\QkRandomTranspose\nabla f(\Xzerok)}\geq \zeta\Dk}$ holds (thus, $\mathscr{C}_k\cap\mathscr{W}_k\subseteq \mathscr{V}_k$ hence $\ickwk\leq \ivk$). From the proof of Theorem~\ref{convDeltaTheorem}, 
\SSaddnew{(\eqref{fourpt22} and~\eqref{eq:B23kineq} for \ASTRAF, the last inequality of Case~1 for \ASTRAR),}
we know that for all $k\geq k_0$, almost surely, $$\E{\ivk(\Phikplusone-\Phik)|\fqalgebra}\leq \ivk\left(-\frac{1}{4}c_1\nu\normii{\QkRandomTranspose\nabla f(\Xzerok)}\Dk\right),$$which implies that for all $k\geq k_0$, almost surely, 
\begin{equation}\label{eq:V}
\begin{split}
\E{\ivk(\Phikplusone-\Phik)|\fqalgebra}&\leq \ickwk\left(-\frac{1}{4}c_1\nu\normii{\QkRandomTranspose\nabla f(\Xzerok)}\Dk\right)\\
&\leq \ickwk\left(-\frac{1}{4}c_1\nu\alpha_Q\normii{\nabla f(\Xzerok)}\Dk\right)\leq \ickwk\left(-\frac{1}{4}c_1\nu\SSaddnew{\varepsilon}\alpha_Q\Dk\right).
\end{split}
\end{equation}
On the other hand it follows from~\eqref{fourpt16} that trivially, for all $k\geq k_0$, $\E{\ivkc(\Phikplusone-\Phik)|\fqalgebra}\leq \Dk^2$ almost surely, which, combined with~\eqref{eq:V}, yields for all $k\geq k_0$, almost surely
\begin{equation*}
\E{\Phikplusone-\Phik|\fqalgebra}\leq \iwk\ick\left(-\frac{1}{4}c_1\nu\SSaddnew{\varepsilon}\alpha_Q\Dk\right)+\Dk^2,
\end{equation*}
hence for all $k\geq k_0$, almost surely
\begin{equation}\label{fortysix}
\E{\left.\E{\Phikplusone-\Phik|\fqalgebra}\right\rvert\falgebra}\leq -\frac{1}{4}c_1\nu\SSaddnew{\varepsilon}\alpha_Q\Dk\ick\E{\iwk|\falgebra}+\Dk^2,
\end{equation}
where we used the $\falgebra$-measurability of $\ick$ and $\Dk$. Defining $c_\SSaddnew{\varepsilon}:=\frac{1}{4}c_1\nu\alpha_Q(1-\beta_Q)\SSaddnew{\varepsilon}$, it follows from~\eqref{fortysix} that for all $k\geq k_0$, almost surely,
\begin{equation}\label{fortyseven}
\E{\Phikplusone-\Phik|\falgebra}\leq -c_\SSaddnew{\varepsilon}\Dk\ick+\Dk^2,
\end{equation}
hence \begin{equation}\label{fortyeight}
\SSaddnew{c_\SSaddnew{\varepsilon}\E{\sum_{k= k_0}^{\bar k-1}\Dk\ick}\leq \sum_{k= k_0}^{\bar k-1}\E{\Phik-\Phikplusone}+\E{\sum_{k= k_0}^{\bar k -1}\Dk^2}\leq \E{\Phi_{k_0}}+\sum_{k\geq k_0}\E{\Dk^2}<\infty,}
\end{equation}
\SSaddnew{and $\bar k\to\infty$ by monotone convergence, } where we used the Fubini-Tonelli theorem to permute $\Esp$ and $\sum$, the first inequality follows from~\eqref{fortyseven}, the second one results from $\Phik>0$ for all~$k$ and the last one follows from Theorem~\ref{convDeltaTheorem}. Thus, it follows from~\eqref{fortyeight} and the definition of $\mathscr{C}_k$ that
{\bl 
\begin{equation*}
\sum_{k\geq k_0}\Dk\ick<\infty\qquad\text{almost surely}.
\end{equation*}
By the definition of~$K_1$, almost surely, $\Dk\leq\SSaddnew{\varepsilon}\alpha_Q/\zeta$ for every $k\geq K_1$. Consequently, if $k\in \mathscr{K}_{\SSaddnew{\varepsilon},1}$, then $k\geq K_1$ and \(\|\nabla f(X_k^0)\|_2>\SSaddnew{\varepsilon}\), and hence~$\mathscr{C}_k$ occurs. Therefore, almost surely,
\[\sum_{k\in\mathscr K_{\SSaddnew{\varepsilon},1}}\Dk=\sum_{\accolade{k\in\mathscr{K}_{\SSaddnew{\varepsilon}}:k\geq\krandom_1}}\Dk
\leq\sum_{k\geq k_0}\Dk\ick<\infty,\]
}
hence, almost surely,
\begin{equation*}
\sum_{k\in\mathscr{K}_{\SSaddnew{\varepsilon}}}\Dk=\sum_{\accolade{k\in\mathscr{K}_{\SSaddnew{\varepsilon}}:k<\krandom_1}}\Dk+\sum_{\accolade{k\in\mathscr{K}_{\SSaddnew{\varepsilon}}:k\geq\krandom_1}}\Dk<\infty,
\end{equation*}
which completes the proof.
\end{proof}

\newpage

\clearpage
\framebox{\parbox{\linewidth}{
The submitted manuscript has been created by UChicago Argonne, LLC, Operator of 
Argonne National Laboratory (``Argonne''). Argonne, a U.S.\ Department of 
Energy Office of Science laboratory, is operated under Contract No.\ 
DE-AC02-06CH11357. 
The U.S.\ Government retains for itself, and others acting on its behalf, a 
paid-up nonexclusive, irrevocable worldwide license in said article to 
reproduce, prepare derivative works, distribute copies to the public, and 
perform publicly and display publicly, by or on behalf of the Government.  The 
Department of Energy will provide public access to these results of federally 
sponsored research in accordance with the DOE Public Access Plan. 
http://energy.gov/downloads/doe-public-access-plan.}}

\end{document}